\documentclass{amsart}
\usepackage{graphicx} % Required for inserting images
\usepackage[utf8]{inputenc} % allow utf-8 input
\usepackage[T1]{fontenc}    % use 8-bit T1 fonts
\usepackage{hyperref}       % hyperlinks
\usepackage{url}            % simple URL typesetting
\usepackage{booktabs}       % professional-quality tables
\usepackage{amsfonts}       % blackboard math symbols
\usepackage{nicefrac}       % compact symbols for 1/2, etc.
\usepackage{microtype}      % microtypography
\usepackage{xcolor}         % colors
\usepackage{latexsym}
\usepackage{epsfig}
\usepackage{amssymb}
\usepackage{amstext}
\usepackage{amsgen}
\usepackage{amsxtra}
\usepackage{amsgen}
\usepackage{bm}     
\usepackage{algorithm}
\usepackage{algpseudocode} 

\usepackage{amsthm,mathrsfs,mathtools}
\usepackage{bbm}
\usepackage{enumerate,amsmath,subfigure,color}
\usepackage{tikz}
\usepackage{comment}
\usepackage{ulem}
\usepackage{enumitem}

\DeclareGraphicsExtensions{.eps,.pdf,.jpg,.png}

\usepackage[foot]{amsaddr}

\makeatletter
\newcommand{\addresseshere}{%
  \enddoc@text\let\enddoc@text\relax
}
\makeatother

\newtheorem{thm}{Theorem}[section]

\theoremstyle{remark}

\theoremstyle{definition}

\newcommand{\R}{\mathbb{R}}

\renewcommand{\epsilon}{\varepsilon} % epsilon

\DeclareMathOperator*{\argmin}{arg\,min}
\def\<{\langle}
\def\>{\rangle}

\newcommand{\Input}{\Require}
\newcommand{\Output}{\Ensure}
\algrenewcommand\algorithmicrequire{\textbf{Input:}}
\algrenewcommand\algorithmicensure{\textbf{Output:}}

\newcommand{\indic}{\chi}

\title{The learned range test method for the inverse inclusion problem}

\author[S.\ Sun]{Shiwei Sun \!$^{1, 2}$} \email{swsun@seu.edu.cn}
\author[G.\ S.\ Alberti]{Giovanni S.\ Alberti \!$^2$} \email{giovanni.alberti@unige.it}
\address{$^1$ School of Mathematics, Southeast University, Nanjing 210096, P.R. China.}
\address{$^2$ Machine Learning Genoa Center (MaLGa), Department of Mathematics, Department of Excellence 2023–2027, University of Genoa, Italy.}


\date{May 24, 2025}


\begin{document}

\begin{abstract}
We consider the inverse problem consisting of the reconstruction of an inclusion $B$ contained in a bounded domain $\Omega\subset\R^d$ from a single pair of Cauchy data $(u|_{\partial\Omega},\partial_\nu u|_{\partial\Omega})$, where $\Delta u=0$ in $\Omega\setminus\overline B$ and $u=0$ on $\partial B$. We show that the reconstruction algorithm based on the range test, a domain sampling method, can be written as a neural network with a specific architecture. We propose to learn the weights of this network in the framework of supervised learning, and to combine it with a pre-trained classifier, with the purpose of distinguishing the inclusions based on their distance from the boundary. The numerical simulations show that this learned range test method provides accurate and stable reconstructions of polygonal inclusions. Furthermore, the results are superior to those obtained with the standard range test method (without learning) and with an end-to-end fully connected deep neural network, a purely data-driven method. 
\end{abstract}

\maketitle

\section{Introduction}\label{Sec:Introduction}
Electrostatic and thermal imaging are important techniques with various applications in scientific and industrial disciplines, including non-destructive testing and evaluation, remote sensing, ultrasound imaging and so on; see, for instance, \cite{Akduman-2002-18, Hu-2020, Kress-2014, Aparicio-1996-12} and the references therein. Roughly speaking, these techniques aim to recover the unknown boundary impedance or an abnormal inclusion (for example, cavity or crack) inside the conducting medium, from the knowledge of  voltage and current measurements  or of the temperature  on the boundary of the medium. In mathematical terms, it can be modeled as an inverse boundary value problem for the Laplace equation.

In this work, we focus on the reconstruction of a perfectly insulated inclusion inside a homogeneous medium. Specifically, let $\Omega\subset\mathbb {R}^d$ be a bounded domain. Assume that there is a perfectly insulated inclusion $B\subseteq\Omega$. The electric potential $u\in H^{1}(\Omega \setminus \overline{B})$ satisfies
\begin{equation*}
\begin{cases}
\Delta u =0 & \mathrm{in}\ \Omega \setminus \overline{B}, \\
u = g & \mathrm{on}\ \partial \Omega,\\
u = 0 & \mathrm{on}\ \partial B,
\end{cases}
\end{equation*}
with a prescribed boundary condition $g\in H^{\frac{1}{2}}(\partial \Omega)$. Our goal is to recover the inclusion $B$ from a single measurement $\{g,\partial_\nu u|_{\partial\Omega}\}$, which is a nonlinear and ill-posed inverse problem. This is a variation of Calderón's inverse problem \cite{borcea-2002,uhlmann-2009}, modeling electrical impedance tomography, in which the unknown is not a simple inclusion but the full conductivity density, and the data consists of all Cauchy pairs on $\partial\Omega$. 

In recent years, several theoretical and numerical works have provided many insights into this inverse problem; see  $\S$\ref{sub:review} below for a short literature review. In this work, we are particularly interested in a non-iterative method, the range test (RT). This is a qualitative method, or more precisely, a domain sampling method \cite{cakoni2011linear,Cakoni-2014}. This method was initially introduced for the inverse acoustic scattering problem \cite{Potthast-2003-19}, and was later extended to the Oseen problem \cite{Zia-2016-304}, for which the convergence of the RT is studied in detail. Recently, the RT was employed to solve the inverse boundary value problem for the heat equation \cite{sun-2024}. For additional applications of the RT we refer the readers to \cite{Jakubik-2008-58, Potthast-2005-75, Alves-2009-25}. In essence, the RT tackles the corresponding inverse problem by detecting the range of a specific operator. In other words, the philosophy behind the RT is that the presence of an anomaly within the conductor will affect the regularity of the solution and therefore the range of an integral operator. 

The RT has been recently extended to the inverse boundary value problem for the Laplace equation discussed above  \cite{Lin-2021-15, Sun-2023-485}. A sophisticated indicator function based on the range of a boundary integral operator relative to a prescribed test domain allows us to establish whether the test domain contains the inclusion. Hence, one can detect the inclusion by taking the intersection of all the test domains containing the inclusion. The RT method inherits the advantages of qualitative approaches: it can perform the reconstruction with just one boundary measurement and requires minimal {\it a priori} information about the inclusion. However, it has three notable limitations. First, a careful selection of parameters, which often vary with the inclusion, is needed for the numerical implementation of the RT method, making it impractical. Second, the reconstruction it provides is relatively rough. Lastly, the method becomes computationally intensive when applied to large sampling domains.

Over the past decade, deep learning methods have attracted extensive attention because of their salient success in solving inverse problems; see, for example, \cite{Chung-2024,Ye-2018, Tanyu-2023-39,Raissi-2019-378,2023bubba}. Directly solving inverse problems by constructing an end-to-end neural network and ignoring the physical model loses interpretability and often leads to less effective reconstruction methods \cite{2019-arridge-etal-actanumerica,mukherjee-etal-2023}. As a result, the integration of conventional methods with deep learning methods has become a popular direction in recent years, also in the context of qualitative methods. For example, in \cite{Guo-2021-43}, the authors propose a deep learning method for electrical impedance tomography based on the direct sampling method. They employ fully connected neural networks and convolutional neural networks to approximate the indicator function of the inclusion. The study \cite{Le-2023-784} integrates the orthogonality sampling method with deep neural networks to reconstruct the geometry of a penetrable object from scattering measurements generated by a single incident wave. In \cite{Ning-2024-40}, the authors propose a sampling-based deep learning method to tackle the inverse medium scattering problem. Specifically, they integrate the results generated by the direct sampling method with a U-Net to learn the relationship between the indicator function and the true resolution of the target. In \cite{Li-2024-40}, the authors consider the inverse inhomogeneous medium problem using a learning-based iterative algorithm. Briefly speaking, they make use of a deep neural network to learn the {\it a priori} information on the shape of the unknown scatterer, and then refine the reconstruction by using the iteratively regularized Gauss-Newton method. Inspired by the linear sampling method in \cite{Lorenzo-2024}, the work \cite{Meng-2024-SIAM} employs kernel methods for inverse source and scattering problems, which are modeled by integral operators with specific kernels.
Inclusion reconstruction may also be viewed as an imaging problem, for which convolutional neural networks (CNN) are often used because of their capability to capture spatial structures.
Several CNN-based methods have been proposed to solve this kind of inverse problems, see e.g.\ \cite{Cen-2023-493, Wu-2021-7, Ongie-2020}.

In this paper, to overcome the weaknesses of the RT, we propose an RT-based deep learning method, denoted as the learned range test (LRT), to realize the reconstruction of the inclusion with high accuracy and efficiency.
The main contributions of our study are the following. 
\begin{itemize}
    \item We show that the reconstruction algorithm based on the RT described in \cite{Sun-2023-485} can be expressed as a neural network with a specific architecture. The weights of this NN can then be learned with a classical supervised learning approach.
    \item We illustrate that it is impractical to expect that one universal network with this architecture can provide accurate reconstructions for all inclusions. As a remedy, we propose to combine several neural networks with this architecture with a pre-trained classifier, in order to treat in a different way the inclusions that have a different distance from the boundary. 
    \item  Various numerical experiments are conducted to verify the effectiveness of the LRT. The numerical results suggest that the LRT has the capacity to reconstruct inclusions with high accuracy and efficiency from one measurement, especially if compared with the standard RT method or with a fully connected end-to-end deep neural network. 
\end{itemize}

The remainder of this paper is structured as follows. The mathematical formulation of our inverse problem is stated in Section~\ref{Sec:IP}, where a literature review concerning this problem is also included. In Section~\ref{Sec:RT}, we briefly revisit the RT from the theoretical and numerical perspectives. Section~\ref{sec:LRT} is devoted to the LRT method, where the architecture of the neural network and the three-step strategy are described. Various numerical experiments are carried out in Section \ref{Sec:Experiments}. In Section~\ref{Sec:Conclusion}, we provide some concluding remarks. Finally, some details on the computational aspects, with a particular focus on the training of the neural networks of this work, and partial numerical results are discussed in the appendix.

\section{The inverse problem}\label{Sec:IP}
In this section, we first briefly revisit the mathematical model of the inverse problem of our interest. After that, we give a short literature survey on the relative work. 

\subsection{Mathematical formulation}\label{sub:mathIP}

Let $\Omega$ be a bounded domain in $\mathbb {R}^d$ $(d=2,\,3)$ with  $C^2$ boundary $\partial \Omega$, and $\nu$ denote the unit outward normal vector to $\partial \Omega$. Assume that there is a perfectly insulated inclusion (hereinafter referred to as inclusion) $B\subseteq\Omega$ with Lipschitz boundary such that $\Omega \setminus \overline{B}$ is connected. The physical field $u\in H^{1}(\Omega \setminus \overline{B})$ satisfies 
\begin{equation}
\begin{cases}
\Delta u =0 & \mathrm{in}\ \Omega \setminus \overline{B}, \\
u = g & \mathrm{on}\ \partial \Omega,\\
u = 0 & \mathrm{on}\ \partial B,
\end{cases}
\label{u}
\end{equation}
with a prescribed boundary condition $g\in H^{\frac{1}{2}}(\partial \Omega)$. Furthermore, let $\tilde u\in H^{1}(\Omega)$ be the solution to 
\begin{equation}
\begin{cases}
\Delta \tilde u =0 & \mathrm{in}\ \Omega, \\
\tilde u = g & \mathrm{on}\ \partial \Omega.
\end{cases}
\label{ub}
\end{equation}
 We first highlight that the forward problems \eqref{u} and \eqref{ub} are well-posed; see, for instance, \cite{Lin-2021-15, Evans-2010}.
Physically, $\tilde u$ models the state of the homogeneous medium $\Omega$ without the inclusion $B$, which is usually denoted as the background solution. The differences between $u$ and $\tilde u$ might reveal some information on $B$. Setting $\omega: =\omega(B)= u-\tilde u$, we have
\begin{equation}\label{omega}
\begin{cases}
\Delta \omega =0 & \mathrm{in}\ \Omega \setminus \overline{B}, \\
\omega = 0 & \mathrm{on}\ \partial \Omega,\\
\omega = -\tilde u & \mathrm{on}\ \partial B.
\end{cases}
\end{equation}
 Therefore, for a prescribed $g\in H^{\frac{1}{2}}(\partial \Omega)$, we are able to obtain the Neumann data $\partial_{\nu} \omega|_{\partial \Omega} = \partial_{\nu} u|_{\partial \Omega} - \partial_{\nu} \tilde u|_{\partial \Omega}$ by measuring the normal derivative of the physical field on $\partial\Omega$ and by solving the forward problem \eqref{ub}.

In this paper, the inverse problem of our concern is an inverse boundary value problem for the Laplace equation. Specifically, we aim to detect the inclusion $B$ from one boundary measurement $\partial_{\nu}\omega|_{\partial \Omega}$. This is a nonlinear and ill-posed inverse problem, with applications to 
nondestructive testing \cite{Aparicio-1996-12, Akduman-2002-18}.

\subsection{Literature review}\label{sub:review}
Over the past decades, numerous works have been devoted to this inverse problem from the theoretical and numerical points of view. We list some representative literature as follows. From the theoretical perspectives, the uniqueness and conditional stability for the determination of the target from Cauchy data are well studied; see, for instance, \cite{Kress-2014, Akduman-2002-18, Karageorghis-2009-17, Alessandrini-2001-176, Beretta-1999-30, Bukhgeim-1999-15, Haddar-2015-21} and the references therein. Analogous results in the context of the Calderón problem can be found in \cite{1996-seo,1989-friedman-isakov,2022-beretta-francini,2023-alberti-arroyo-santacesaria,Hanke_2024}.

There also exist various numerical schemes developed to solve this inverse problem. These may be divided into two main categories, namely, the iterative methods and the qualitative methods. Just as the name implies, the former approaches aim to realize the numerical reconstruction iteratively. For example,  Kress and Rundell in \cite{Kress-2005-21} convert this inverse boundary problem into a two-by-two system, which consists of an ill-posed linear equation and a nonlinear integration equation. Then, these two equations can be solved by the regularized Newton iteration after linearizing the nonlinear equation. We refer to \cite{Burger-2001-17, Bonnet-2008-24,Bourgeois-2010-4,Cakoni-2007-1, Karageorghis-2009-17} for additional details on iterative methods.

In general, iterative methods yield relatively good reconstructions via the iterative procedure. However, their drawbacks are the heavy computational burden and the need for extensive {\it a priori} information about the target, which might not be available in practice. Therefore, as a remedy, qualitative methods were introduced in order to obtain some limited information of the target. Some popular qualitative methods include the sampling method \cite{Li-2008-30, Chow-2014-30}, the enclosure method \cite{Ikehata-1999-15, Ikehata-2002-18}, the probe method \cite{Ikehata-1998-23} and the MUSIC-type method \cite{Ammari-2008-108}. We also refer the reader to \cite{Kim-2002-62, Cakoni-2014, Kirsch-2007, Potthast-2006-22} for additional qualitative methods and the relative technical details. Compared with the iterative methods, the qualitative methods allow for a good trade-off between accuracy and computational cost. 

\section{Review of the range test}\label{Sec:RT}
The range test (RT) is a qualitative method introduced for the inverse acoustic scattering problem \cite{Potthast-2003-19}. It is capable of determining the convex support of the scatterer by using the far field patterns from one or a few incident waves, when the physical properties of the scatterer, for instance, the shape, size and location, are unknown. This method, in the context of the inverse boundary value problem for the Laplace equation considered in $\S$\ref{sub:mathIP}, was studied in \cite{Lin-2021-15, Sun-2023-485} from the theoretical and numerical points of view.  In this section, we will briefly revisit the theoretical foundations and the numerical algorithm for the use of the RT in solving our inverse problem. 

\subsection{Theoretical analysis of the RT}
We start with a quick summary of the whole method for solving the inverse problem introduced in $\S$\ref{sub:mathIP}. Assuming that the solution $u$ of \eqref{u} cannot be analytically extended across $\partial B$ (because otherwise $B$ would be invisible), we can design an indicator function $\mathcal I(G)$ for a  convex test subdomain $G\Subset \Omega$ such that $\mathcal I(G)$ is finite if and only if $\overline B\subset \overline G$. Therefore,  we can deduce some geometric information about $B$ or, more precisely, of the convex hull of $B$, by taking the intersection of various test domains $G$ for which $\mathcal I(G)<\infty$.

Now, let us explain the details of the RT more explicitly. Let $G\Subset \Omega$ be a bounded convex Lipschitz domain, which will serve as a test domain. Letting $\mathcal K(x,y)$ be the Green function of the Laplace equation in $\Omega$ with Dirichlet boundary condition, we define the single-layer potential operator $\mathcal S\colon H^{-\frac{1}{2}}(\partial G) \to H^{1}(\Omega \setminus \partial G) $  by
\begin{equation*}
  \mathcal S[\psi](x) := \int_{\partial G} \mathcal K(x,\,y)\,\psi(y)\, ds(y),\quad \psi \in H^{-\frac{1}{2}}(\partial G),
\end{equation*}
where $ds(y)$ represents either the line element ($d = 2$) or the surface element ($d = 3$) of $\partial G$. Furthermore, we introduce the operator $\mathcal R\colon H^{-\frac{1}{2}}(\partial G) \to H^{-\frac{1}{2}}(\partial \Omega)$ as the  normal derivative of $\mathcal S$ on $\partial \Omega$, namely,
\begin{equation*}
  \mathcal R[\psi](x) := \partial_\nu \mathcal S[\psi](x) = \int_{\partial G} \frac{\partial}{\nu(x)}\mathcal K(x,\,y)\,\psi(y)\, ds(y).
\end{equation*}
 Then, we consider the following integral equation:
\begin{equation}\label{main-eq}
  \mathcal R[\psi] = \partial_{\nu}\omega.
\end{equation}
The operator $\mathcal R$ is compact and injective  \cite{Lin-2021-15}, and therefore we can solve \eqref{main-eq} by using Tikhonov regularization:
\begin{equation*}
  \psi_\alpha = (\alpha I + \mathcal R^*\mathcal R)^{-1}\mathcal R^*\partial_{\nu}\omega,
\end{equation*}
where $\alpha>0$ indicates the regularization parameter, and $\mathcal R^*\colon H^{-\frac{1}{2}}(\partial \Omega)\to H^{-\frac{1}{2}}(\partial G)$ is the adjoint operator of $\mathcal R$. Denoting the range of $\mathcal R$ by $\operatorname{Im}\mathcal R$, classical regularization theory allows us to obtain the following characterization (see, for instance, \cite{Engl,Nakamura-2015}):
\begin{itemize}
    \item if $\partial_{\nu}\omega \in \operatorname{Im}\mathcal R$, $\lim_{\alpha \to 0} \psi_\alpha = \tilde \psi$ for some $\tilde \psi\in H^{-\frac{1}{2}}(\partial G)$;

    \item if $\partial_{\nu}\omega \not\in \operatorname{Im}\mathcal R$, $\lim_{\alpha \to 0} \|\psi_\alpha\|_{H^{-\frac{1}{2}}(\partial G)} = +\infty$.
\end{itemize}
Furthermore, assuming that $\omega$ cannot be analytically extended across $\partial B$, $\partial_{\nu}\omega$ belongs to $\operatorname{Im}\mathcal R$ if and only if $\overline B\subset \overline G$; see \cite{Lin-2021-15} for more details. Therefore, we can investigate the relationship between $B$ and $G$ by checking the convergence of $\|\psi_\alpha\|_{H^{-\frac{1}{2}}(\partial G)}$ as $\alpha \to 0$. In order to do this, we define for a test domain $G$ the indicator function $\mathcal{I}(G)$ by
\begin{equation*}
  \mathcal I(G) :=
\begin{cases}
\lim\limits_{\alpha \to 0} \Vert \psi_\alpha \Vert_{H^{-\frac{1}{2}}(\partial G)}  & \textrm{if the finite limit exists,} \\
+\infty & \textrm{otherwise}.
\end{cases}
\end{equation*}
Now, we are ready to state the key result on the RT.
\begin{thm}[\cite{Lin-2021-15}]\label{main-thm}
    Assume that $B\Subset \Omega$ is a convex polygon satisfying the distance property:
    \begin{equation}\label{eq:distanceproperty}
        \mathrm{diam}(B)<\mathrm{dist}(B, \partial \Omega),
    \end{equation}
    where $\mathrm{diam}(B)$ represents the diameter of $B$, and $\mathrm{dist}(B, \partial \Omega)$ is the distance between $B$ and $\partial \Omega$. Assuming that the test domain $G$ is a convex polygon, we have:
    \[
\overline{B}\subset \overline{G} \iff \mathcal I(G) <+\infty.
    \]
\end{thm}
Henceforth, for the sake of simplicity, we call a test domain $G$ positive if $\mathcal I(G)<+ \infty$. The above theorem suggests that the target $B$ can be reconstructed by taking the intersection of all positive test domains.

\subsection{Numerical implementations of the RT}\label{Algorithm-RT}
The numerical implementation of the RT was completed in the previous work \cite{Sun-2023-485}, where the numerical algorithm and the experiments are discussed in detail. In this subsection, we revisit the algorithm and some representative numerical results.

First, we approximate  $\mathcal I(G)$ by
\begin{equation}\label{I-RT}
\mathcal I(\alpha, G): = \Vert \mathcal W(\alpha, G)\partial_{\nu}\omega \Vert_{L^{2}(\partial G)}.
\end{equation}
for a sufficiently small $\alpha>0$, where 
$\mathcal W(\alpha, G): =(\alpha I + \mathcal R^*\mathcal R)^{-1}\mathcal R^*$.
Next, a large threshold $\mathcal C$ should be selected to establish whether $\overline{B}\subset \overline{G}$, as in Theorem \ref{main-thm}. Specifically, we consider the following two cases:
\begin{equation}\label{eq:RT2cases}
\begin{cases}
 \mathcal I(\alpha, G) < \mathcal C \implies \overline{B}\subset \overline{G},\\
 \mathcal I(\alpha, G) \geq \mathcal C \implies \overline{B}\not\subset \overline{G}.
\end{cases}
\end{equation}

\begin{figure}%[htbp]
 \centering
 \subfigure[Test domains] 
  {
	\begin{minipage}{3cm}
	\centering
	\includegraphics[scale=0.25]{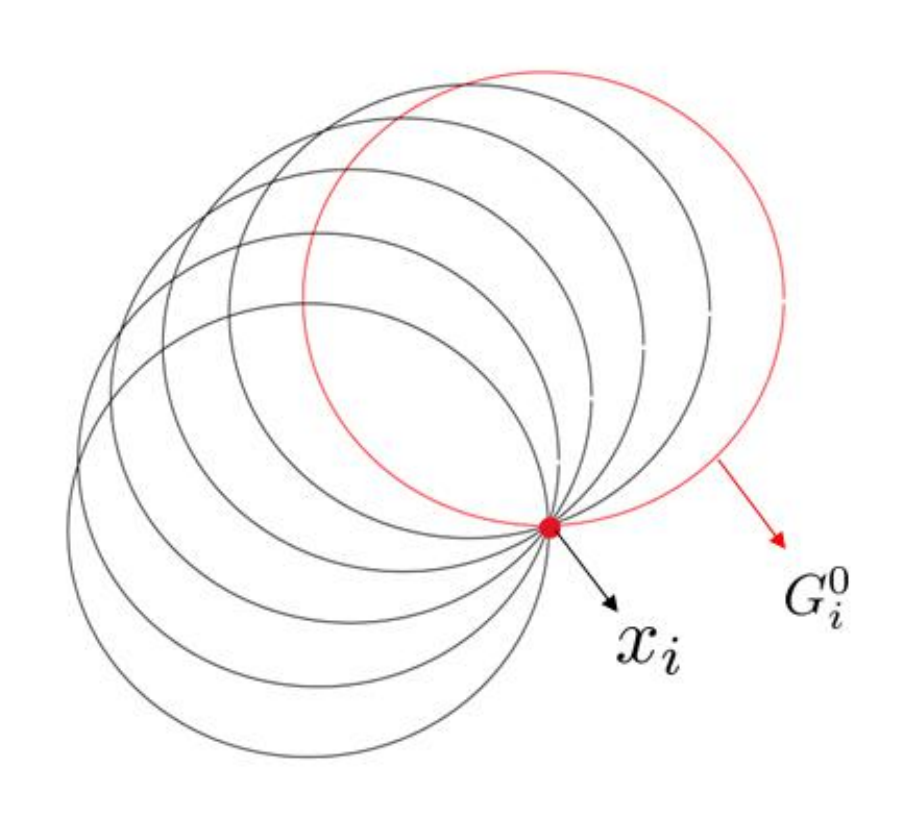}
	\end{minipage}
 \label{Algorithm-illustration1}
  } 
  \hspace{5mm}
  \subfigure[$x_i \in B$]
  {
	\begin{minipage}{3cm}
	\centering
	\includegraphics[scale=0.25]{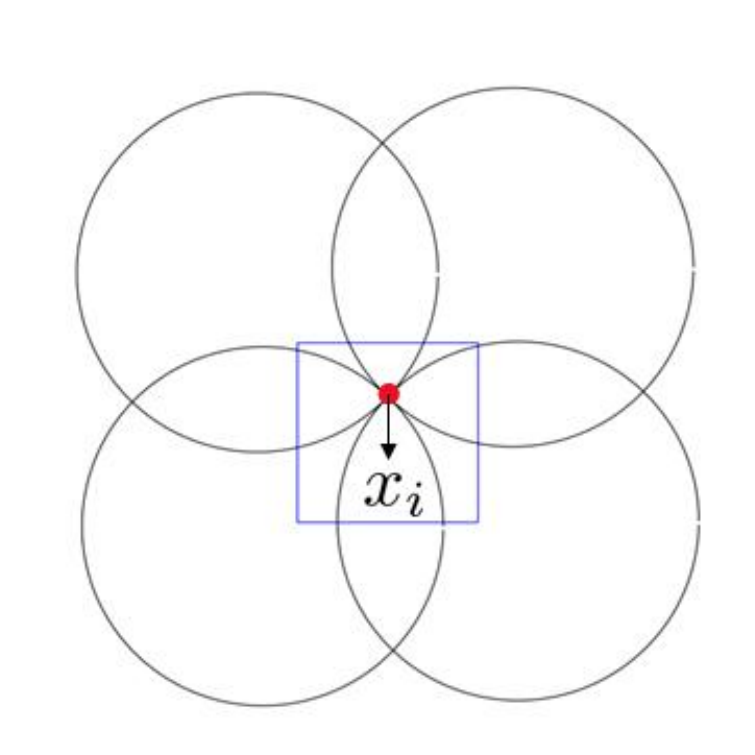}
	\end{minipage}
  \label{Algorithm-illustration2}
  }\hspace{5mm}
   \subfigure[$x_i \not\in B$]
  {
	\begin{minipage}{3cm}
	\centering
	\includegraphics[scale=0.25]{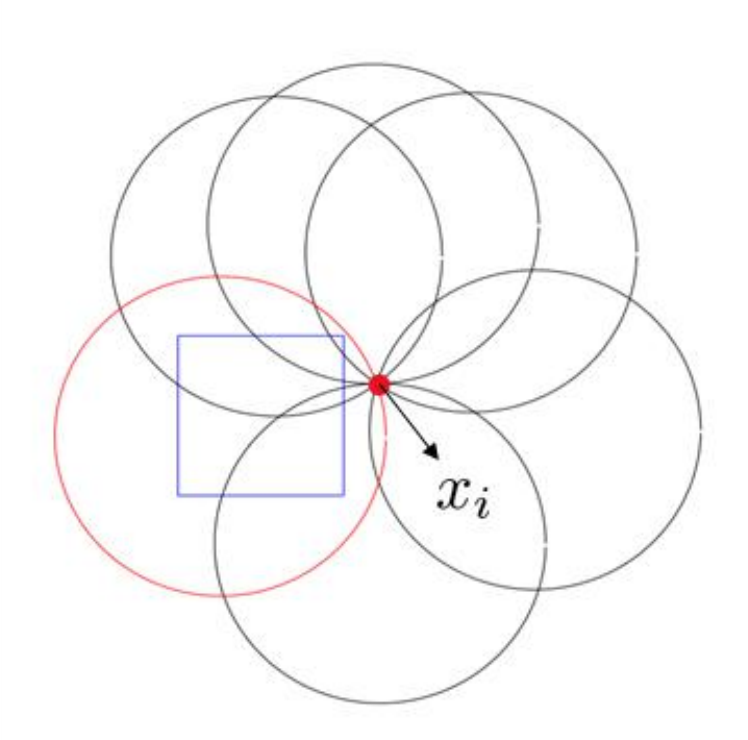}
	\end{minipage}
  \label{Algorithm-illustration3}
  }

 \caption{(a): Some test domains in $P(x_i)$. In (b) and (c), the blue square represents the inclusion $B$. (b): if $x_i \in B$, then all test domains in $P(x_i)$ are not positive. (c): if $x_i \not\in B$, then there exists a positive test domain in $P(x_i)$, for instance, the red one.}
 \label{Algorithm-illustration}
\end{figure}

Let us briefly explain how to generate test domains. Letting $\tilde \Omega$ be a bounded domain such that $\overline{B} \Subset \tilde \Omega \subset \Omega$, we attempt to reconstruct $B$ by searching in $\tilde \Omega$, and therefore $\tilde \Omega$ is called the initial searching area. Consider a grid $\Gamma = \{x_1,\,x_2,\,\cdots,\, x_N\}\subset \tilde \Omega$. For every point $x_i\in \Gamma$, we generate a convex polygonal test domain $G^0_i$ with $x_i \in \partial G^0_i$. Then, we can rotate $G^0_i$ around $x_i$ by different angles to obtain a family of test domains $P(x_i):=\{G_i^j:j=0,\dots,M\}$ associated with $x_i$, see Figure~\ref{Algorithm-illustration1}. 

In view of \eqref{eq:RT2cases}, we have the following two cases:
\begin{itemize}
    \item if $x_i\in \overline B$, then every test domain $G\in P(x_i)$ satisfies $\overline B\not\subset \overline G$, and therefore $\mathcal I(\alpha, G) \geq \mathcal C$ (see Figure~\ref{Algorithm-illustration2});

    \item if $x_i\not\in \overline B$, we can  find  a test domain $G\in P(x_i)$ such that $\overline B\subset \overline G$, assuming that the initial test domain $G^0_i$ is chosen appropriately and that sufficiently many rotations are considered (see Figure~\ref{Algorithm-illustration3}). Thus, $\mathcal I(\alpha, G) < \mathcal C$.
\end{itemize}
Therefore, defining for every $x_i\in \Gamma$ the quantity
\begin{equation}\label{eq:Ii-min}
    I_i:= \min_{G \in P(x_i)} {\mathcal I(\alpha, G)},
\end{equation}
we can use the rules
\begin{equation*}
\begin{cases}
I_i < \mathcal C \implies x_i \notin \overline{B},\\
I_i \geq \mathcal C \implies x_i \in \overline{B},
\end{cases}
\end{equation*}
to obtain an approximation of the inclusion $B$ given by
\begin{equation}\label{eq:fx_i}
f(x_i)=
\begin{cases}
    1 & \text{if $I_i \geq \mathcal C$,}\\
    I_i/\mathcal C & \text{if $I_i < \mathcal C$.}
\end{cases}
\end{equation}

\begin{figure}%[htbp]
 \centering
 \subfigure[Trapezoid ($\delta = 0$)]
  {
	\begin{minipage}{.48\textwidth}
	\centering
	\includegraphics[width=\columnwidth]{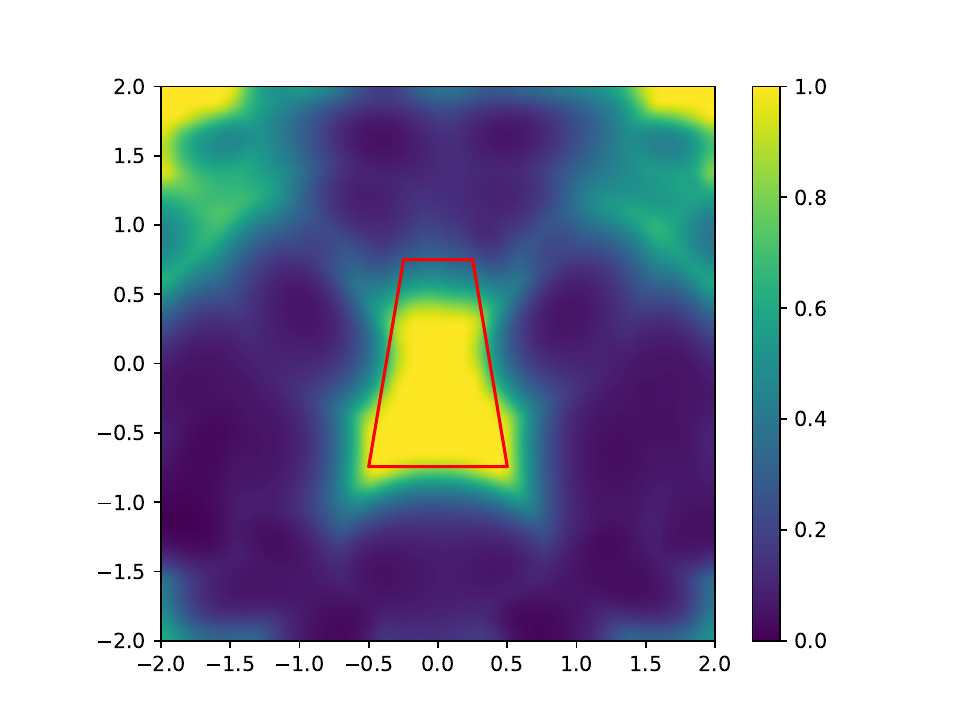}
	\end{minipage}
  }\hfill
  \subfigure[Trapezoid ($\delta = 3\%$)]
  {
	\begin{minipage}{.48\textwidth}
	\centering
	\includegraphics[width=\columnwidth]{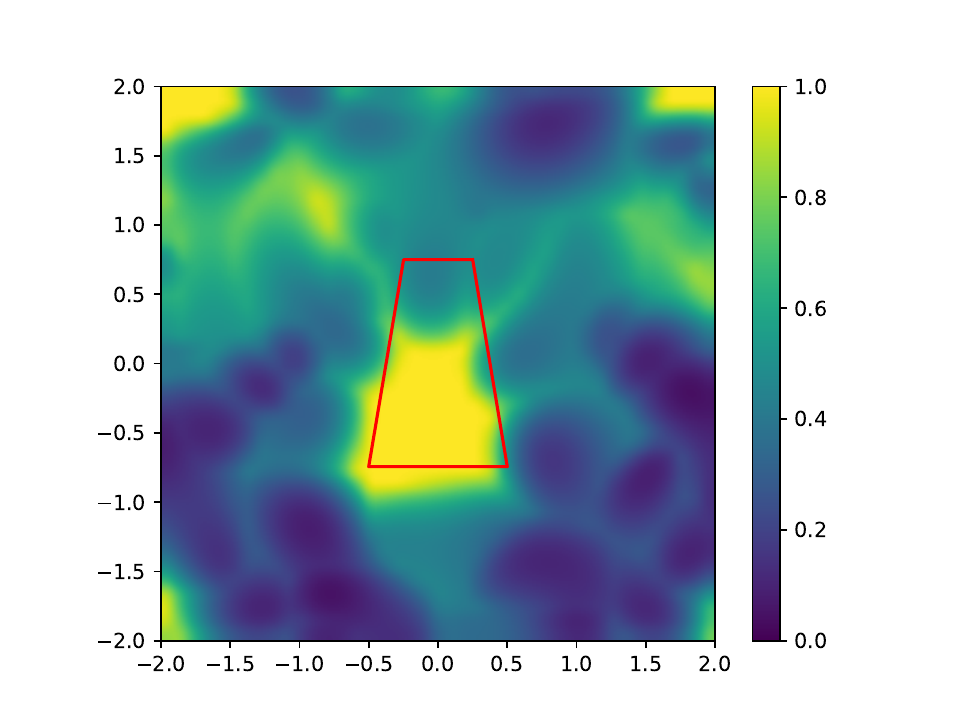}
	\end{minipage}
  }

   \subfigure[Triangle ($\delta = 0$)]
  {
	\begin{minipage}{.48\textwidth}
	\centering
	\includegraphics[width=\columnwidth]{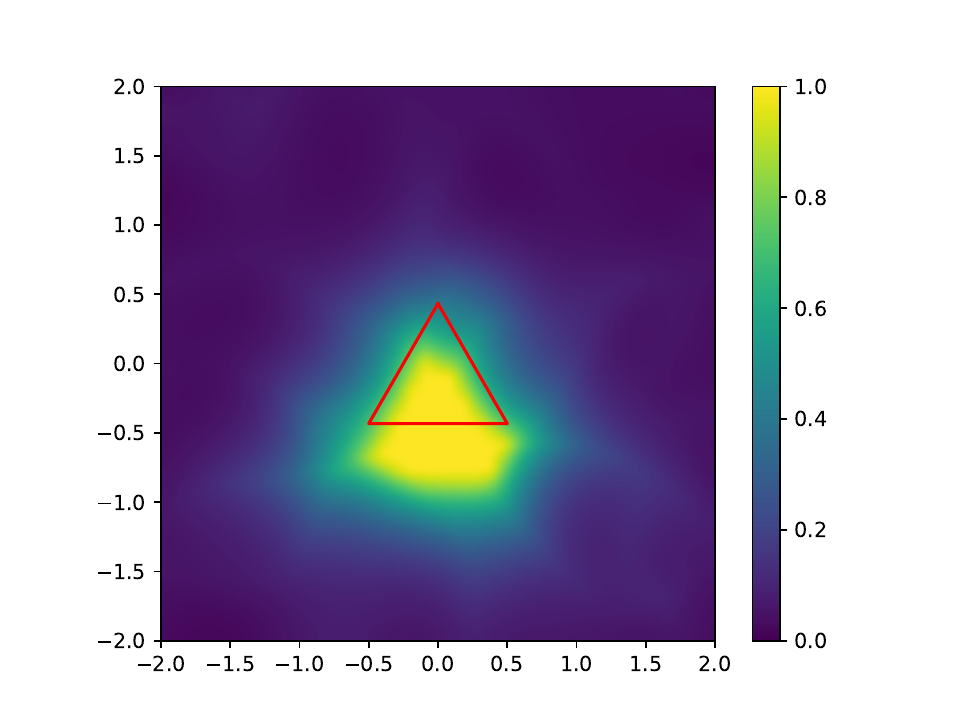}
	\end{minipage}
  }\hfill
  \subfigure[Triangle ($\delta = 3\%$)]
  {
	\begin{minipage}{.48\textwidth}
	\centering
	\includegraphics[width=\columnwidth]{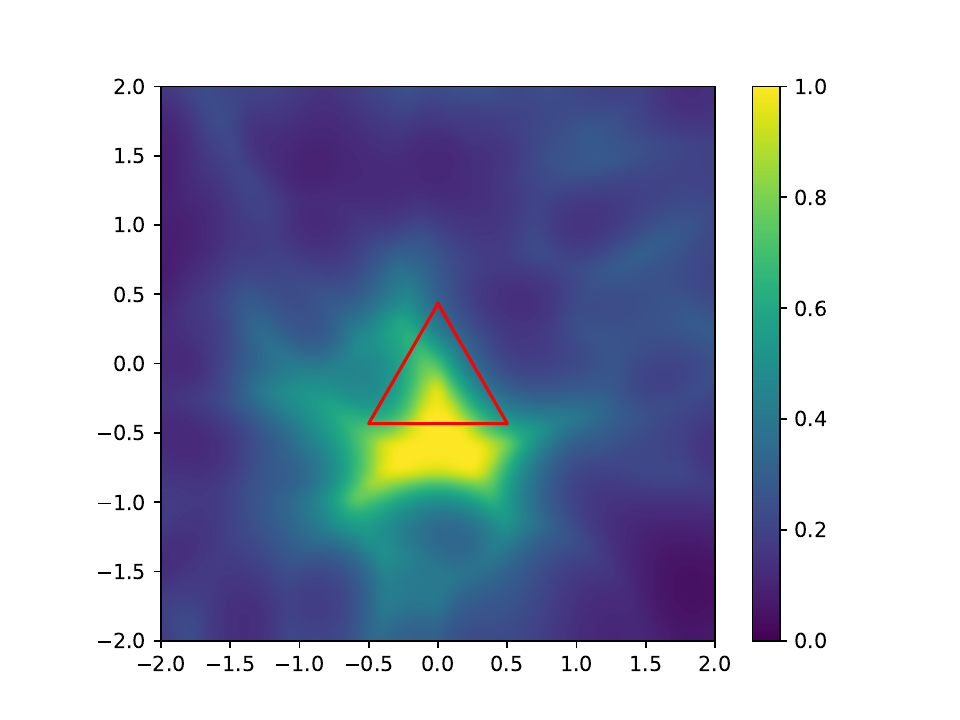}
	\end{minipage}
  }
 \caption{The red polygon represents the target $B$. (a) and (c) are the results generated by using noise-free measurement, whereas (b) and (d) correspond to the measurement with $3\%$ noise.}
 \label{RT-results}
\end{figure}

Next, we exhibit some representative numerical results obtained by using the above algorithm. Here, $d=2$ and $\Omega$ is the disk centered at the origin with radius $10$, and the initial searching area is $\tilde\Omega=[-2, 2]\times[-2,2]$. To test the stability of the algorithm,  random Gaussian noise with noise level $\delta$ is added to the boundary measurements. The numerical results are shown in Figure~\ref{RT-results}, and are taken from \cite{Sun-2023-485}.

This example suggests that the RT can realize a reliable and stable reconstruction of the inclusion from one boundary measurement. However, there are three main limitations of the RT. First, for every (unknown) inclusion $B$, one has to select with caution the parameters (which usually vary with the inclusion), including the grid points, the test domains, the regularization parameter $\alpha$ and the threshold $\mathcal C$, which is impractical. Second, as shown in Figure~\ref{RT-results}, this approach only provides a rough reconstruction of the inclusion. Finally, it struggles with the huge computational burden if the grid  $\Gamma\subset\tilde \Omega$ is large.

\section{The learned range test method}\label{sec:LRT}
In this section, we propose a deep learning method based on the RT, denoted as learned range test (LRT), to realize the reconstruction of the inclusion with high accuracy and efficiency. Indeed,  we show that the whole implementation process of the RT can be incorporated into the framework of deep learning, which allows us to train the neural network to achieve better performance. In addition, we will illustrate that it is impractical to expect one universal neural network (NN) that can reconstruct all inclusions well. As a remedy, a three-step strategy is proposed to address this problem.

\subsection{The RT as a neural network}\label{sub:LRT}

\begin{figure}
    \centering
   \includegraphics[width=\textwidth]{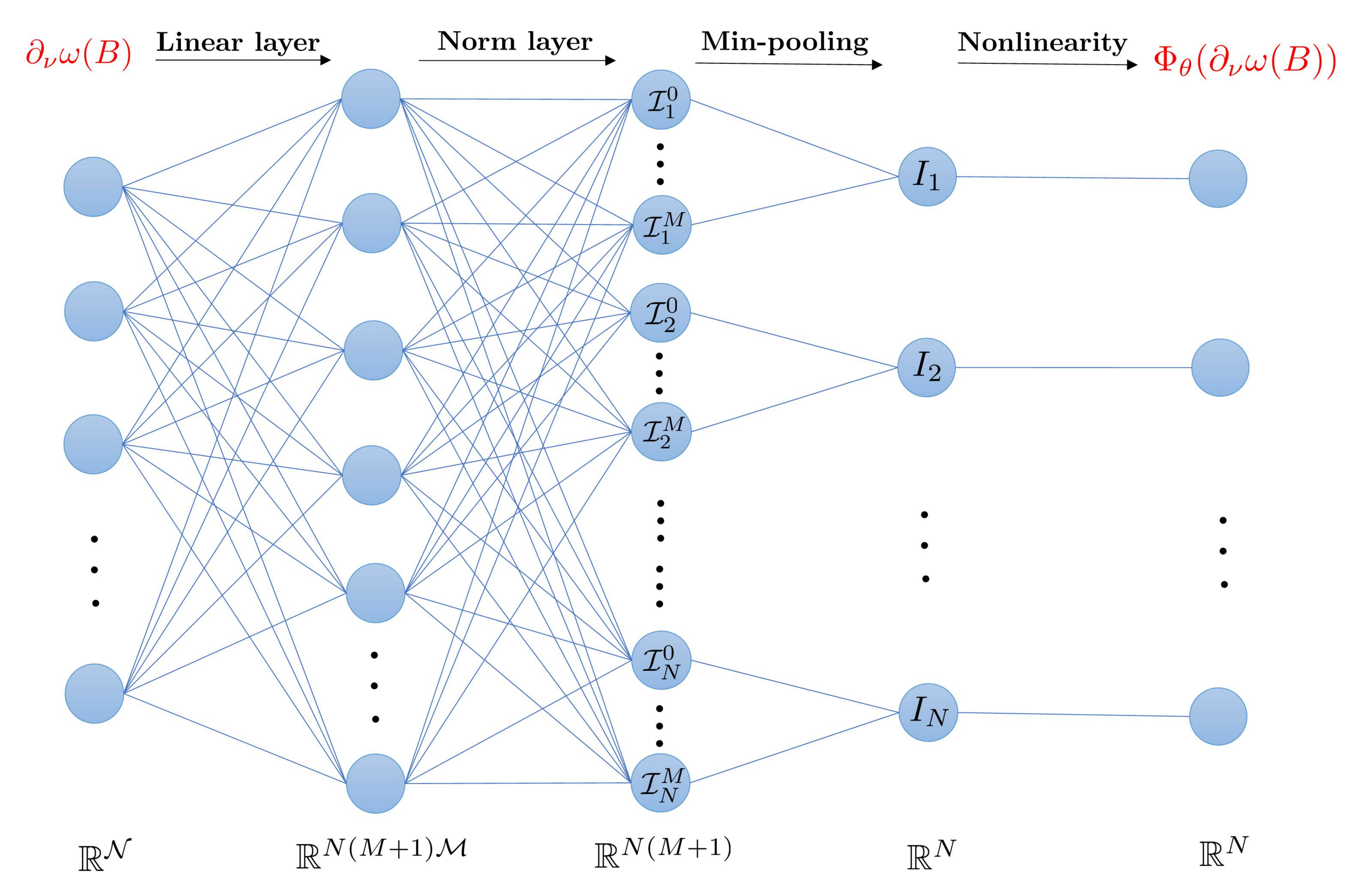}
    \caption{The architecture of the NN relative to the RT method.}
    \label{fig:whole-archi}
\end{figure}

We now show how the algorithm based on the RT and discussed in Section~\ref{Algorithm-RT} may be written as a NN, which is the composition of five steps.
\begin{enumerate}[label=\arabic*)]
    \item  \label{step1} \textit{ The linear layer.} The input of the algorithm (and of the NN) is the measurement $\partial_\nu \omega \in H^{-1/2}(\partial \Omega)$, and in the first  layer we compute the linear map
    \[
\partial_\nu\omega\in H^{-1/2}(\partial \Omega)  \longmapsto \mathcal W(\alpha,G_i^j) \partial_\nu \omega \in L^2(\partial G_i^j)
    \]
    for $i=1,\dots,N$ and $j=0,\dots,M$. 
    \item \label{step2} \textit{The norm layer.} In the second layer we compute $\mathcal I(\alpha,G_i^j)$ (see \eqref{I-RT}) by computing the norm of the output of the first layer:
    \[
\mathcal I(\alpha,G_i^j) = \|\mathcal W(\alpha,G_i^j) \partial_\nu\omega\|_{L^2(\partial G_i^j)}
    \]
    for $i=1,\dots,N$ and $j=0,\dots,M$.
    \item \label{step3}\textit{The min-pooling.} Next, we compute $I_i$ in \eqref{eq:Ii-min} as
    \[
     I_i = \min_{j=0,\dots,M} \mathcal I(\alpha,G_i^j)
    \]
    for $i=1,\dots,N$.
    \item \label{step4} \textit{The bias term.} In view of the threshold $\mathcal{C}$, we compute
    \[
    I_i-\mathcal{C}
    \]
    for $i=1,\dots,N$.
     \item \textit{The nonlinearity.} Finally, we compute the output of the neural network
     \[
     (\sigma(I_i-\mathcal{C}))_{i=1}^N\in\R^N,
     \]
     where $\sigma(x) = \frac{1}{1+e^{-x}}$ is the sigmoid. This is an  approximation of the indicator function of $B$ defined on the grid $\Gamma=\{x_i:i=1,\dots,N\}\subset\tilde\Omega$, as in \eqref{eq:fx_i}. Note that we changed the activation function compared to \eqref{eq:fx_i} because the sigmoid is a much more popular choice.
\end{enumerate}
For simplicity, we wrote the NN directly with infinite-dimensional function spaces. In practice, the network has to be discretized:  the operator $\mathcal{W}(\alpha,G^j_i)$ in step \ref{step1} is approximated by a matrix $W(\alpha,G^j_i)\in \R^{\mathcal{M}\times\mathcal{N}}$, as a linear  map $\R^\mathcal{N}\to\R^\mathcal{M}$. Here, $\mathcal{N}$ denotes the number of scalar measurements of $\partial_\nu\omega$ on $\partial\Omega$ and $\mathbb{R}^\mathcal{M}$ is a descritization of $L^2(\partial G^j_i)$. 
The whole architecture is visualized in Figure~\ref{fig:whole-archi}.

\subsection{The neural network with learned weights}\label{sub:NNlearned}

Inspired by the concept of unrolling \cite{gregor-lecun-2010,2019-arridge-etal-actanumerica} and some related more recent developments  \cite{dehoop-etal-2021,Bubba-2021-14,alberti2024learning}, we propose to use the architecture discussed above, and to learn the weights of the network. More precisely:
\begin{itemize}
    \item the linear maps in step \ref{step1} are replaced by an affine map
    \[
    W\cdot + b\colon \R^\mathcal{N}\to \R^{N(M+1)\mathcal{M}}
    \]
    with $W\in\R^{N(M+1)\mathcal{M}\times\mathcal{N}}$ and $b\in \R^{N(M+1)\mathcal{M}}$;
    \item the constant translation by $\mathcal{C}$ in step \ref{step4} is replaced by a more general bias term
    \[
    I_i \longmapsto I_i - c_i 
    \]
    with $c\in\R^N$.
\end{itemize}

The learnable parameters, namely, the matrix $W\in\R^{N(M+1)\mathcal{M}\times\mathcal{N}}$ and the bias terms $b\in \R^{N(M+1)\mathcal{M}}$ and $c\in\R^N$, are collected into $\theta:=(W,b,c)$. The corresponding NN is denoted by
\[
\Phi_\theta\colon\R^\mathcal{N}\to\R^N.
\]
As will be discussed in Appendix~\ref{app:Training}, the training is performed by minimizing the empirical loss
\begin{equation*}
    \min_\theta \frac{1}{\# \mathcal T} \sum_{B\in\mathcal T} \frac1N\Vert \indic(B)- \Phi_\theta(\partial_\nu\omega(B))\Vert^2,
\end{equation*}
where $\mathcal T$ is a training set of inclusions and  $\indic (B) :=( \mathbbm{1}_{\overline{B}}(x_i))_{i=1,\dots,N}$.

Ideally, we would hope to learn a universal $\Phi_{\theta^*}$ that can be used to reconstruct a polygonal inclusion  as long as it satisfies the distance property mentioned before. Nonetheless, we claim it is an impractical task. Indeed, the norm of $\partial_\nu\omega(B)$ heavily depends on $d(B,\partial\Omega)$ (this is related to the fact that inclusions that are farther from $\partial\Omega$ are harder to reconstruct \cite{alessandrini-scapin-2017,garde-hyvonen-2020}, and have a smaller $\|\partial_\nu\omega(B)\|$), and so a unique NN with the above architecture is unlikely to be able to simultaneously deal with all inclusions. This is consistent to the fact that in the work \cite{Sun-2023-485} the parameters $\alpha$ and $\mathcal C$ were chosen in a suitable way depending on the inclusion.

\subsection{Three-step strategy}\label{S-TSS}

%Assume that we have finished the training process and define the learned parameter set as $\Theta^*: =(W^*, b^*_{LT}, b^*_{c})$. Theoretically, for any inclusion $B_t\in\mathcal B_{\tilde\Omega}$, we know $X_{\Theta^*}(B_t) = S(\Xi-b^*_{c})$ with $\Xi = (I_1,\, I_2,\,\cdots, I_N)$ denoting the input from $L_3$ to $L_4$ in Figure \ref{LRT-Architecture}(c). Thus, there must exist some components of $(\Xi-b^*_{c})$ are less than $0$ and others are nonnegative, otherwise, we cannot obtain any reasonable reconstruction from $X_{\Theta^*}(B_t)$. On the one hand, it is easy to see that $\Xi$ is closely associated with $\Vert W^*\partial_\nu \omega(B_t)+b^*_{LT}\Vert$, or more precisely, $\Vert \partial_\nu \omega(B_t)\Vert$. On the other hand, we roughly observe from the numerical computation that $\Vert \partial_\nu \omega(B_t)\Vert$ will be larger if $B_t$ is closer to $\partial \Omega$, which also changes the magnitude of the values of the components in $\Xi$. Therefore, it is unreasonable to expect that a $\Phi_\Theta$ can reconstruct all inclusions in $\mathcal B_{\tilde\Omega}$ well when the diameter of $\tilde \Omega$ is large.

As a remedy, we propose a three-step strategy consisting of a training step (TS), a classification step (CS) and a reconstruction step (RS). To begin with, we consider $L$ training sets of inclusions%$\mathcal T_\ell,\,\ell = 1,\,2,\,\cdots,L$, where
\begin{equation*}
    \mathcal T_\ell = \{B^{k}_\ell: k=1,2,\cdots, N_{\ell}\},\qquad \ell = 1,\,2,\,\cdots,L.
\end{equation*}
(The details on these datasets are presented below.)
For the above training sets, we assume that:
\begin{itemize}
    \item all the samples $\{B^{k}_\ell\}$ have a similar size;

    \item the samples in the same $\mathcal T_\ell$ have similar distances from $\partial \Omega$ and the samples from different training sets have different distances from $\partial \Omega$, namely
    \begin{equation*}
    \begin{cases}
        |\mathrm{dist}(B_1, \partial \Omega)-\mathrm{dist}(B_2, \partial \Omega)|<\epsilon, & \mathrm{if}\,B_1, B_2 \in \mathcal T_\ell,\,\ell = 1,\,2,\,\cdots,L,\\
        |\mathrm{dist}(B_1, \partial \Omega)-\mathrm{dist}(B_2, \partial \Omega)|\geq\epsilon, & \mathrm{if}\,B_1\in \mathcal T_\ell, B_2\in \mathcal T_k, \ell\neq k,
    \end{cases}
    \end{equation*}
    for some small constant $\varepsilon>0$.
\end{itemize}

\begin{figure}%[htbp]
\centering
\includegraphics[width=\textwidth]{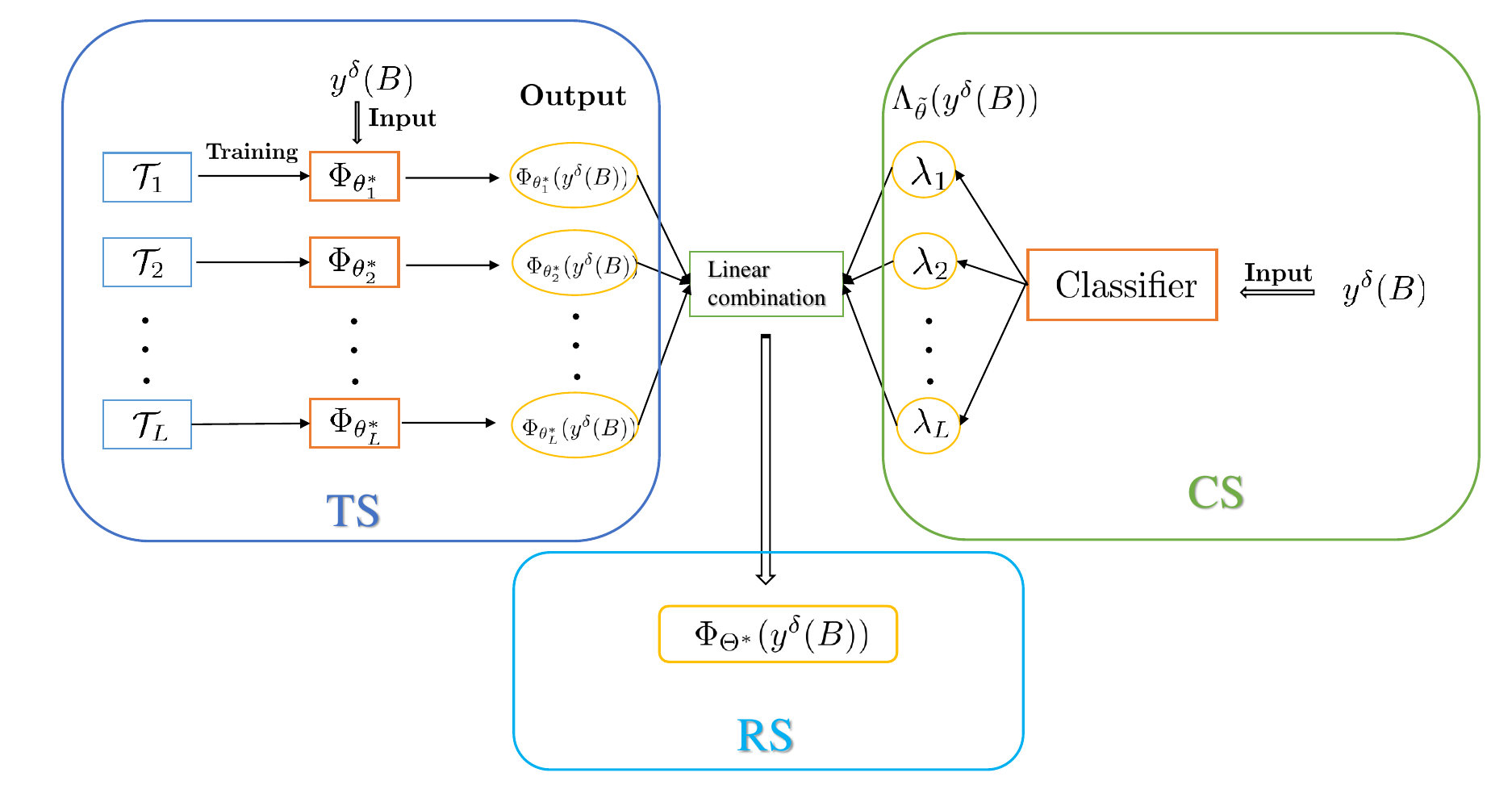}
 \caption{  The diagram of the three-step strategy. Here, $y^\delta(B)$ denotes the noisy version of $\partial_\nu\omega(B)$ with noise level $\delta$, see \eqref{eq:def_noise}. }
 \label{Three-step}
\end{figure}

The three steps of our LRT method are as follows (see Figure \ref{Three-step} for a schematic visualization). 
\begin{itemize}
    \item[(TS)]\label{step-TS} We use $\mathcal T_1, \mathcal T_2, \cdots, \mathcal T_L$ to train  $L$ neural networks $\Phi_{\theta^*_1}, \Phi_{\theta^*_2},\dots, \Phi_{\theta^*_L}$ separately. (The details of the training procedure are presented below.) Each neural network $\Phi_{\theta^*_\ell}$ is expected to perform well with inclusions similar to those in the training set $\mathcal T_\ell$, but not necessarily with those in the others $\mathcal T_k$ for $k\neq \ell$.

    \item[(CS)] \label{step-CS} By using the training set $\mathcal{T}=\mathcal{T}_1\cup\dots\cup\mathcal{T}_L$, we construct a fully connected neural network designed to classify the inclusions according to their distances from $\partial\Omega$. This NN maps
    \[
    \partial_\nu\omega(B)\in\R^\mathcal{N}\longmapsto \Lambda_{\tilde\theta}(\partial_\nu\omega(B)):= (\lambda_1,\,\lambda_2,\,\cdots,\,\lambda_L)\in\R^L,
    \]
    where $\lambda_1+\lambda_2+\cdots+\lambda_L=1$ and $\lambda_\ell\in [0,1]$ is the probability that $B$ belongs to $\mathcal T_\ell$ for $\ell = 1,\,2,\,\cdots, L$. We denote the optimal parameters of this network by $\tilde\theta^*$.

    \item[(RS)] The final reconstruction is obtained by combining the results in (TS) and (CS). More precisely, we set the final reconstruction as
    \begin{equation*}
        \Phi_{\Theta^*}(\partial_\nu\omega(B)) = \lambda_1 \Phi_{\theta_1^*}(\partial_\nu\omega(B))+\cdots+\lambda_L \Phi_{\theta_L^*}(\partial_\nu\omega(B)).
    \end{equation*}
    where $\Theta^*=(\theta_1^*,\dots,\theta_L^*,\tilde\theta^*)$.
    This weighted combination yields a better performance with the inclusions that are on the decision boundaries of the classifier in (CS), compared to the (simpler) alternative of considering only the NN $\Phi_{\theta_\ell^*}$ corresponding to the largest probability $\lambda_\ell$.
    
\end{itemize}

\section{Numerical experiments}\label{Sec:Experiments}
In this section, we carry out some numerical experiments to illustrate that our proposed strategy is effective, accurate for reconstructing the inclusion, and robust to  noise.
\subsection{Setup}
In all the numerical experiments, we restrict ourselves to the case where $d=2$, with the background medium $\Omega$ being a disk centered at the origin with radius $R_\Omega = 10$. The initial searching area $\tilde\Omega = [-2,2]\times [-2,2]$ is divided into a uniform $41\times 41$ Cartesian grid $\Gamma$, so that $N=41^2$. This is the same setup used in \cite{Sun-2023-485}, see also Figure~\ref{RT-results}, which allows us to better compare the results. For every inclusion, the synthetic data $\partial_\nu\omega$ is obtained by solving the forward problems \eqref{u}  and \eqref{ub} with boundary value $g(x,y) = xy$, $(x,y)\in \partial\Omega$. We use $\partial_\nu\omega$ at $\mathcal N=400$ collocation points on $\partial\Omega$ as the measurement. To demonstrate the robustness of our strategy, we consider the following noisy data with the noise level $\delta$:
\begin{equation}\label{eq:def_noise}
    y^\delta(B) = \partial_\nu\omega(B)\left(1+\delta \frac{\psi}{\Vert\psi\Vert}_2\right),
\end{equation}
where $\psi$ is a standard normal random distribution in $\R^\mathcal N$. 
In addition, we let $\mathcal M = 60$ (the number of the collocation points on the boundaries of the test domains $\partial G^j_i$), $M = 24$ (the number of test domains for every point in $\Gamma$).

\subsection{Construction of the dataset}

Following the notation of $\S$\ref{S-TSS}, we consider $L=3$ training sets $\mathcal{T}_1$, $\mathcal{T}_2$ and $\mathcal{T}_3$. For later use, we let
\[
r_0 = 0,\qquad r_1 = 0.4,\qquad r_2=0.8,\qquad r_3 = 1.2. 
\]
For $l=1,2,3$, the training set $\mathcal{T}_l$ consists of 4,000 inclusions. Each of them is randomly constructed as follows. Pick  the following quantities uniformly at random:
\begin{itemize}
    \item $h\in \R^n$ such that $r_{l-1}\le |h| < r_l$;
    \item $r\in [0.4,0.6]$;
    \item $k\in \{3,4,5,6\}$;
    \item and $0\leq \psi_1<\psi_2<\cdots<\psi_k<2\pi$.
\end{itemize}
The inclusion is the convex polygon with $k$ sides obtained by sequentially connecting the points
\begin{equation}\label{eq:Sample}
  h+r(\cos\psi_j, \sin\psi_j),\qquad j=1,\dots,k. 
\end{equation}

This construction guarantees that the inclusions are convex polygons satisfying the distance property \eqref{eq:distanceproperty}. The test sets are generated in a similar way. The only difference between the two procedures is the strategy used to select the parameter $h$. Specifically, to generate the test set,  $h$ is sampled uniformly at random in the interval $[r_0, r_3)$. The purpose of the classification step is precisely to classify the inclusions depending on their distance from the boundary.

\subsection{Training of the networks}\label{sub:trainingnet}

For each $\ell=1,2,3$, we train the network $\Phi_{\theta_\ell}$ introduced in $\S$\ref{sub:NNlearned} by minimizing the loss
\begin{equation}\label{eq:lossthetal}
\mathcal L_\ell (\theta_\ell) = \frac{1}{\# \mathcal T_\ell} \sum_{B\in\mathcal T_\ell}\frac1N \Vert \indic(B)- \Phi_{\theta_\ell}(y^\delta(B))\Vert_2^2.
\end{equation}
This is done by using stochastic gradient descent (see Appendix~\ref{app:Training} for additional details). This requires an initial guess $\theta_\ell^0$, and we propose to use the weights corresponding to the deterministic algorithm discussed in $\S$\ref{sub:LRT}. More precisely, we select the parameters $( \alpha',{\mathcal C}')$ used in the example shown in Figure~\ref{RT-results}, and set $ W' = \left(\mathcal{W}(\alpha', G^0_1), \mathcal{W}(\alpha', G^1_1), \cdots, \mathcal{W}(\alpha', G^M_N)\right)^\top$. The values of $ {\mathcal C}'$ and  of the components of $W'$ are  large, which may lead to the vanishing gradient problem during the training process (see, for instance, \cite{Tan-2019}). Thus, we normalize them and choose
    \[
\theta_\ell^0 = \left(\frac{ W'}{ {\mathcal C'}}, 0, \mathbf{1}\right),
    \]
    where $\mathbf{1}=(1,\dots,1)\in\R^N$.
    Note that this choice is independent of $\ell$. Alternatively, the initial guess can also be chosen at random. We denote the optimal parameter found by $\theta_\ell^*$.

\begin{figure}%[htbp]
\centering
\includegraphics[scale=0.1]{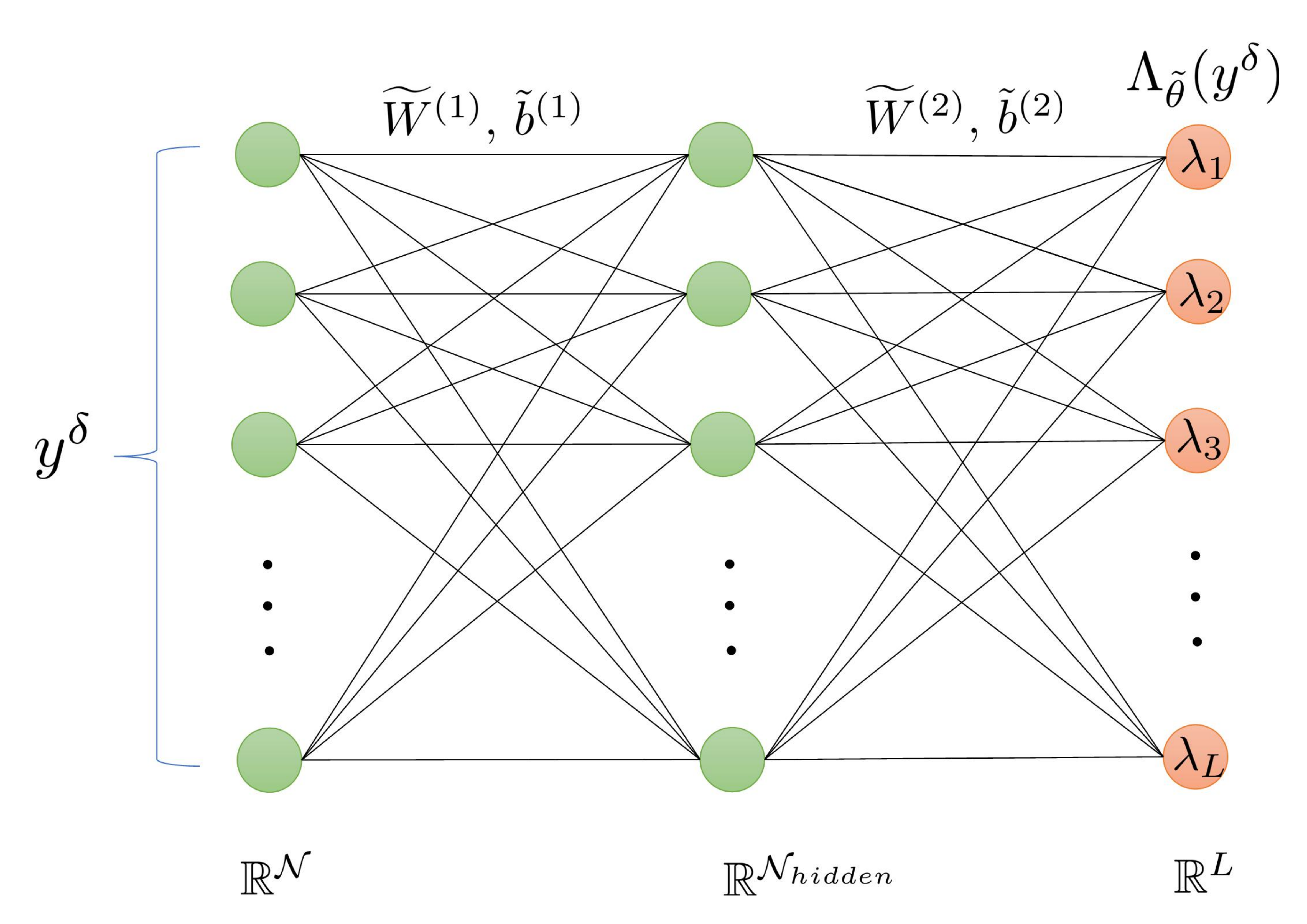}
 \caption{The architecture of the classifier.}
 \label{NN-Class}
\end{figure}

We now discuss the architecture and the training of the classifier introduced in  Section~\ref{S-TSS}, step (CS). It is well known that multi-layer feedforward neural networks have excellent performance in classification tasks; see, for instance, \cite{Murtagh-1991-2}. The network takes as input $y^\delta(B)$ and gives as output a vector $\Lambda_{\tilde\theta}(y^\delta(B))\in\R^3$, whose $l$-th component is an approximation of the probability that the inclusion $B$ is similar to those in the training set $\mathcal{T}_l$.  We construct a fully connected neural network with three layers, with $\mathcal N_{hidden} = 300$  neurons in the hidden layer, see Figure~\ref{NN-Class}.
More precisely, the network $ \Lambda_{\tilde\theta}\colon\R^\mathcal{N}\to\R^3$ is given by
\begin{equation*}
    \Lambda_{\tilde\theta}(y) =\sigma\left(\widetilde{W}^{(2)}\left(\rho\left(\widetilde W^{(1)}y+\tilde b^{(1)}\right)\right)+\tilde b^{(2)}\right),
\end{equation*}
where $\rho$ and $\sigma$ are the  activation functions given by   $\rho(x) = \max(0,x)$ (the rectified linear unit, ReLU, acting pointwise) and 
\(
\sigma(s)_\ell = \frac{e^{s_\ell}}{\sum_{j=1}^3 e^{s_j}}
\)
(the softmax). Here, $\tilde \theta=\bigl(\widetilde W^{(1)},\widetilde W^{(2)},\tilde b^{(1)},\tilde b^{(2)}\bigr)$ collects all the learnable parameters of the network. The training is performed with the training set $\mathcal{T}=\mathcal{T}_1\cup\mathcal{T}_2\cup\mathcal{T}_3$ by minimizing the cross entropy loss, given by 
\begin{equation}\label{eq:lossclass}
\widetilde{\mathcal L}(\tilde\theta) = -\frac{1}{\#\mathcal{T}}\sum_{B\in\mathcal T}\sum_{\ell=1}^3 \eta_\ell(B) \log\bigl(\Lambda_{\tilde\theta}(y^\delta(B))_\ell\bigr),
\end{equation}
where $\eta_\ell(B)=1$ if $B\in\mathcal{T}_\ell$ and $\eta_\ell(B)=0$ otherwise. We denote the optimal parameter by $\tilde\theta^*$. Additional details on the training procedure are presented in Appendix~\ref{app:Training}.

\subsection{Results}
\begin{figure}%[htbp]
\centering
 \subfigure[LRT, $\delta = 0$.]
  {
	\begin{minipage}{.47\textwidth}
 	\centering
	\includegraphics[width=\columnwidth]{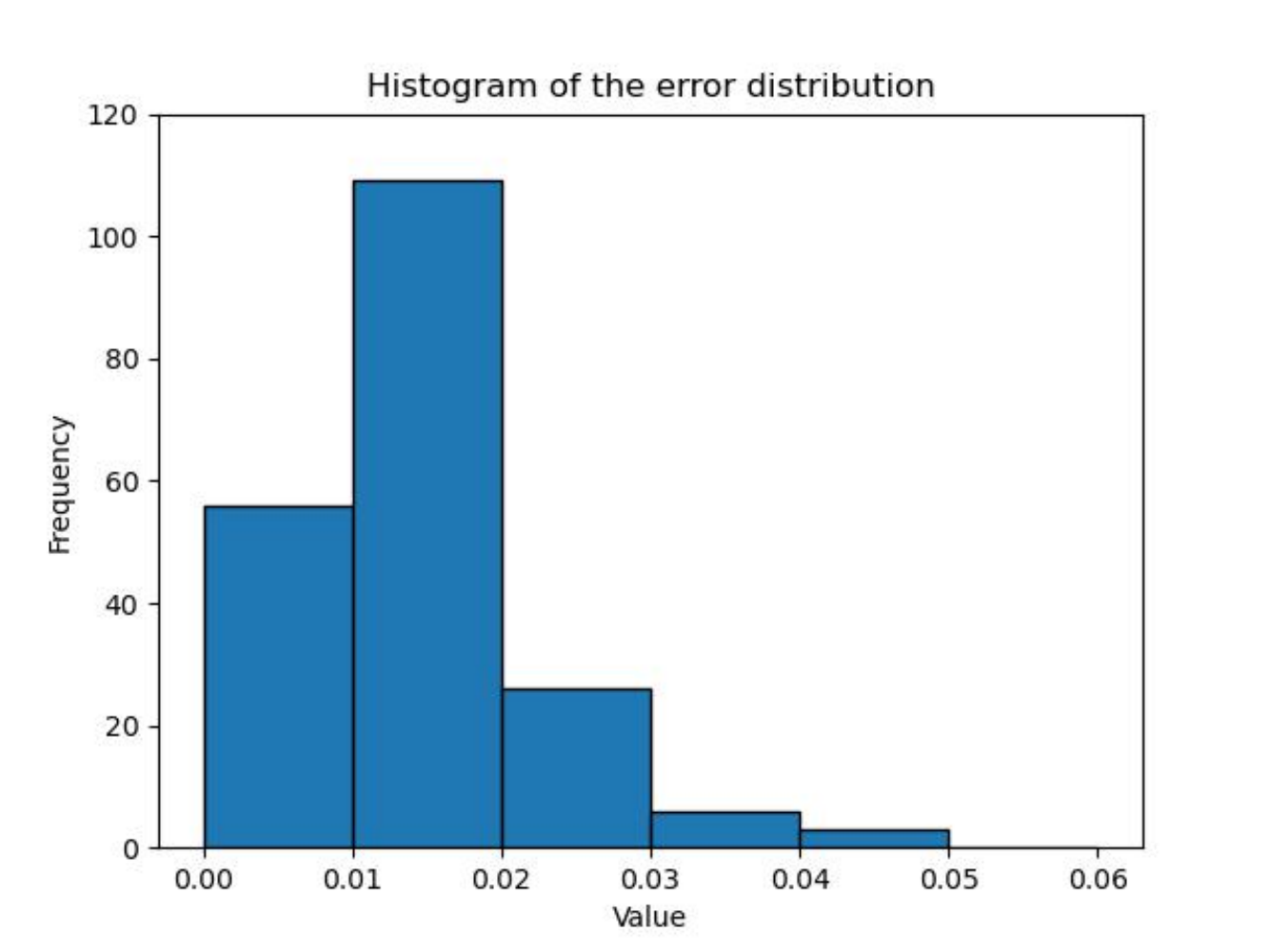}
	\end{minipage}
 \label{Histograma}
  }
   \subfigure[LRT, $\delta = 3\%$.]
  {
	\begin{minipage}{.47\textwidth}
 	\centering
	\includegraphics[width=\columnwidth]{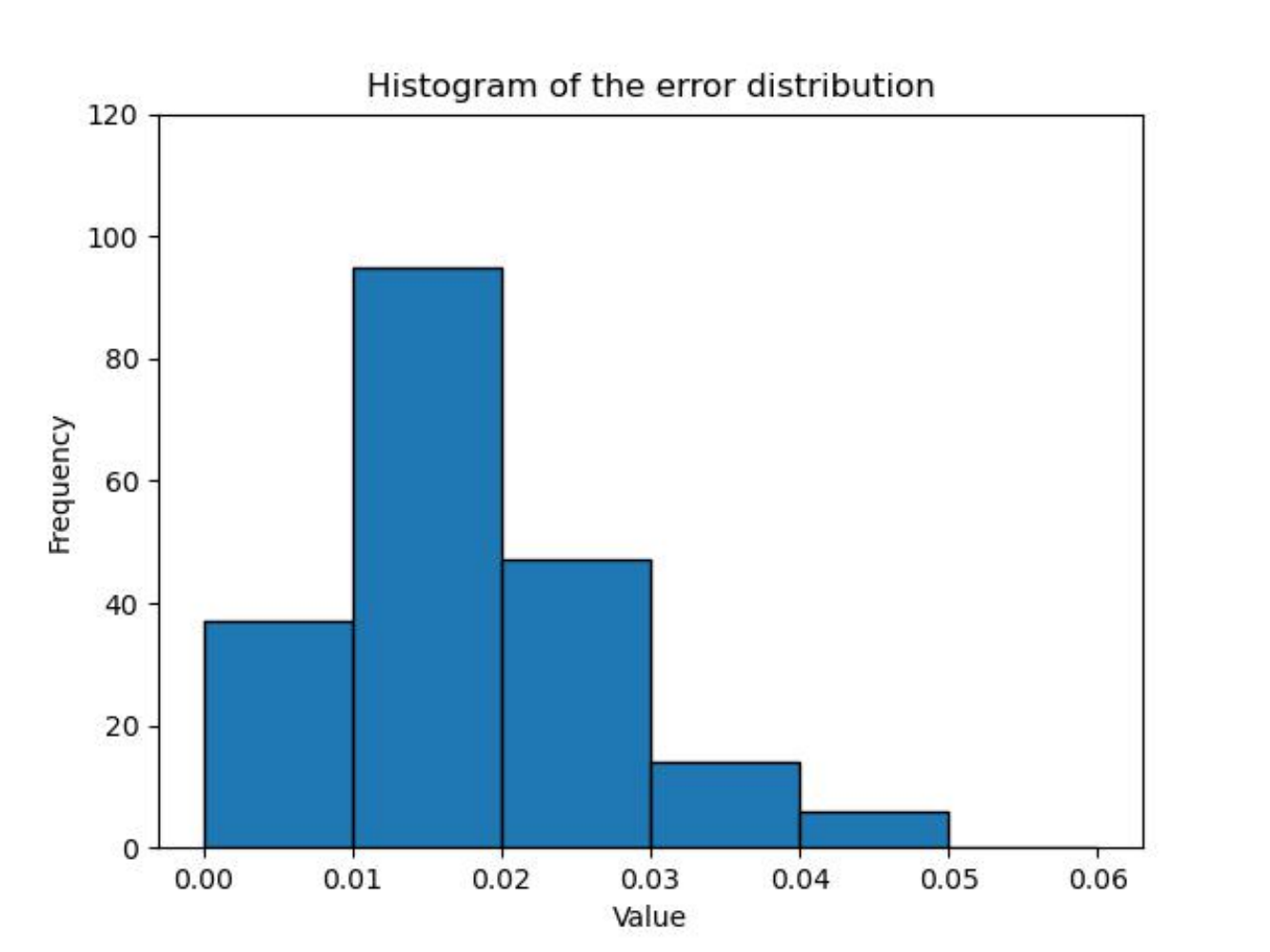}
	\end{minipage}
 \label{Histogramb}
  }
  \subfigure[Deterministic RT, $\delta = 0$.]
  {
	\begin{minipage}{.47\textwidth}
 	\centering
	\includegraphics[width=\columnwidth]{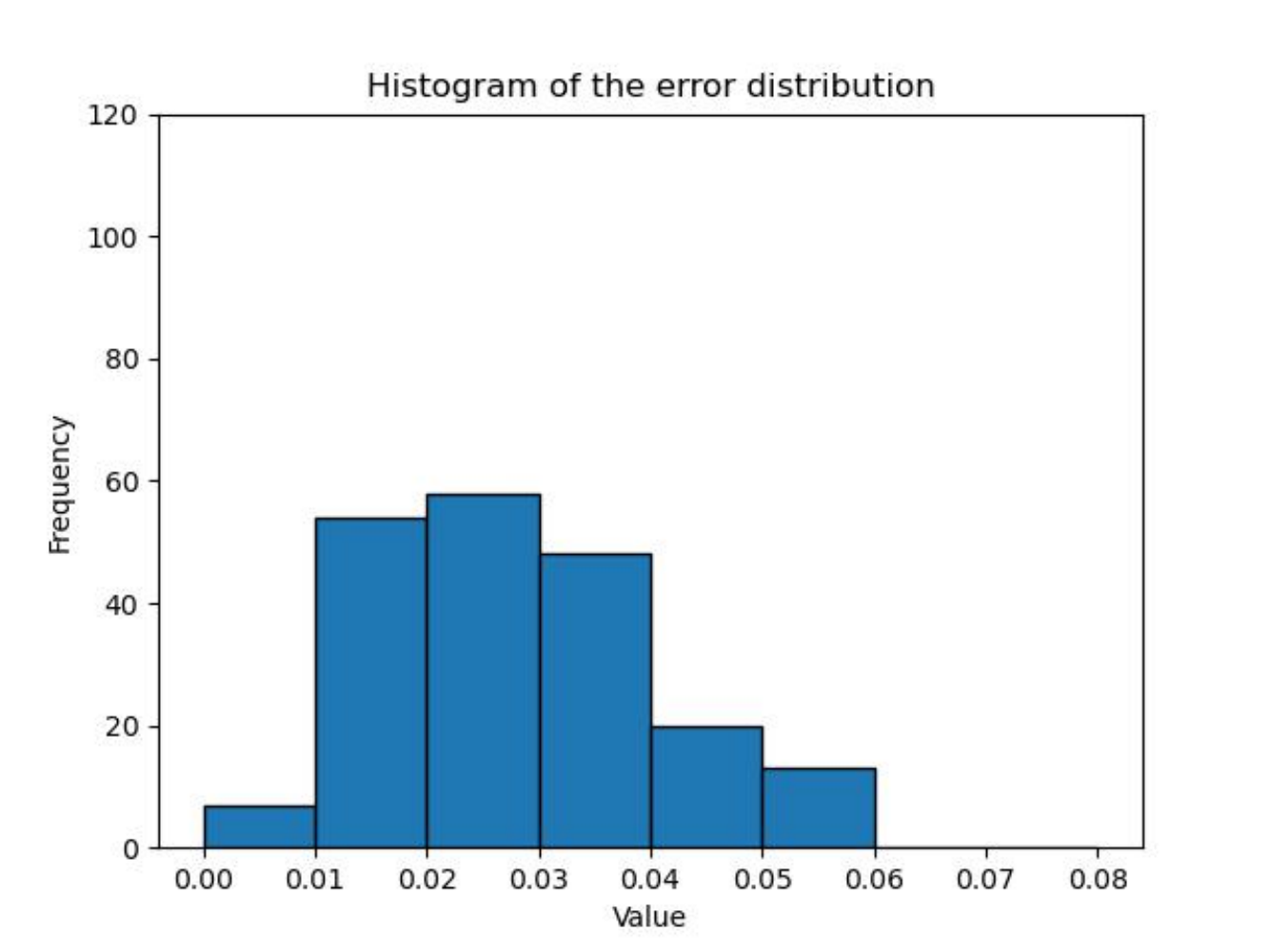}
	\end{minipage}
 \label{Histogramc}
  }
   \subfigure[Deterministic RT, $\delta = 3\%$.]
  {
	\begin{minipage}{.47\textwidth}
 	\centering
	\includegraphics[width=\columnwidth]{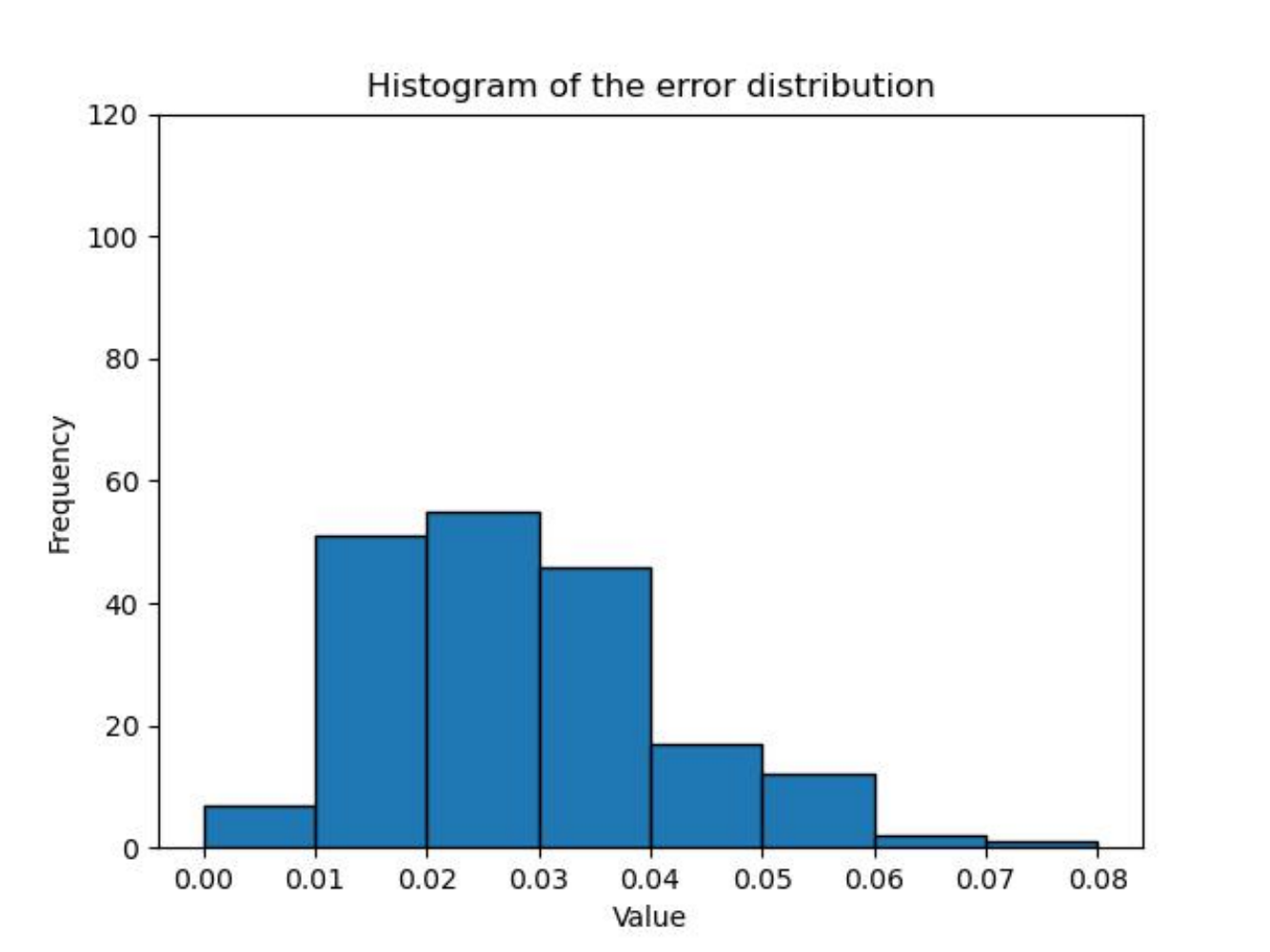}
	\end{minipage}
 \label{Histogramd}
  }
 \caption{The error distribution of the mean squared errors for all the samples in the test set relative to the reconstructions obtained with the LRT and the deterministic RT.} %(a) The error distribution of Scenario 1. (b) The error distribution of Scenario 2. }
 \label{Histogram}
\end{figure}

\begin{figure}%[htbp]
\centering
 \subfigure[0.0077, $0.9168$]
{
	\begin{minipage}{3cm}
 	\centering
	\includegraphics[scale=0.18]{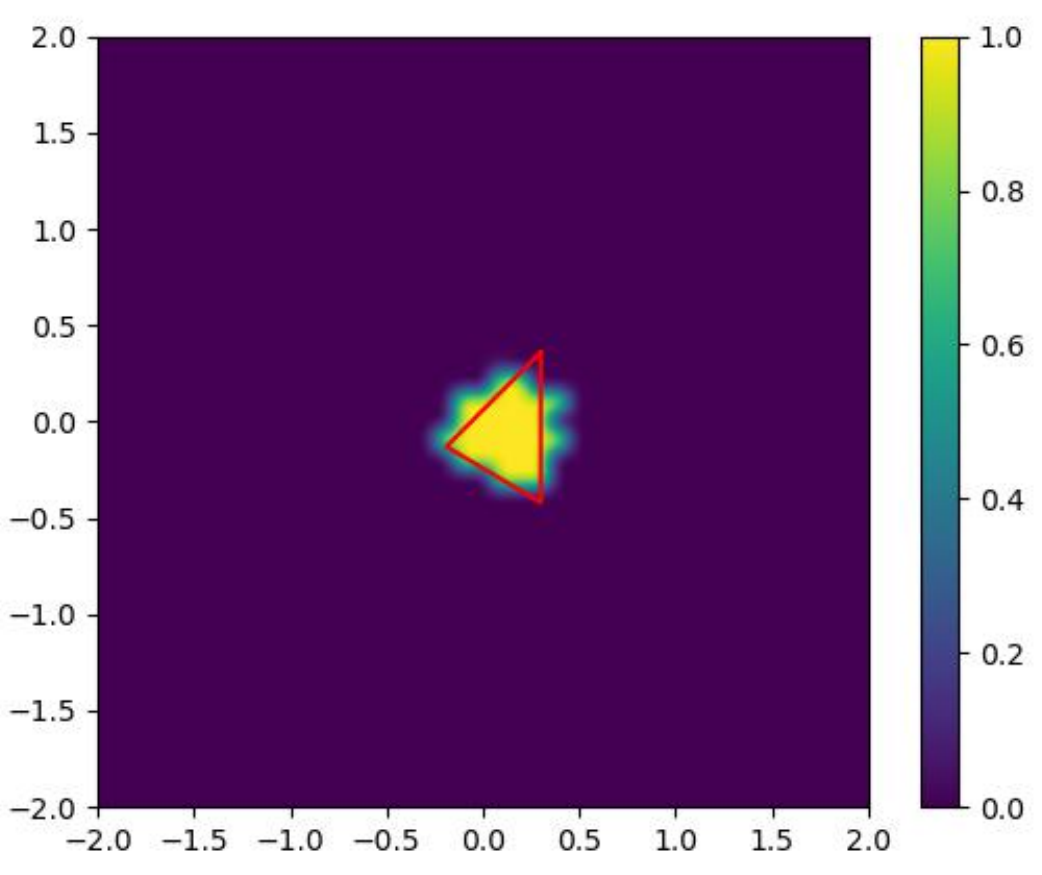}
	\end{minipage}
}
  \subfigure[0.0155, $0.8537$]
  {
	\begin{minipage}{3cm}
 	
	\includegraphics[scale=0.18]{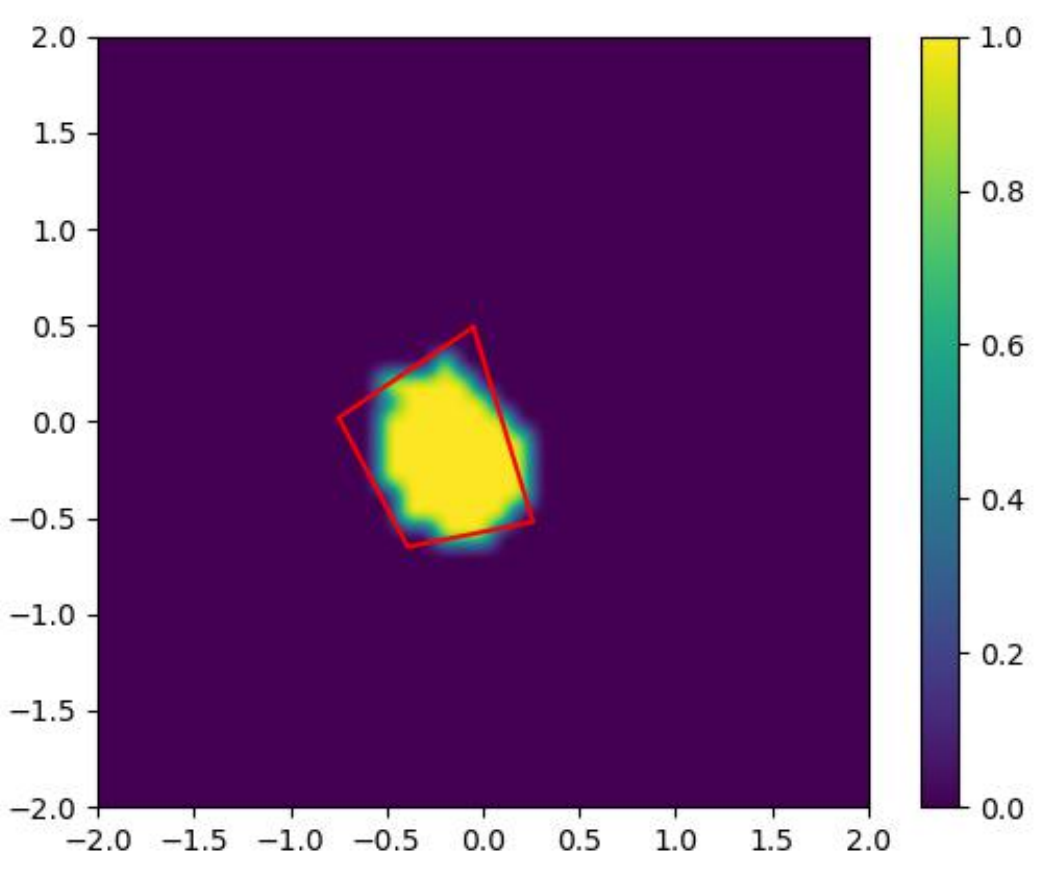}
	\end{minipage}
  }
   \subfigure[0.0077, $0.8956$]
  {
	\begin{minipage}{3cm}
 	\centering
	\includegraphics[scale=0.18]{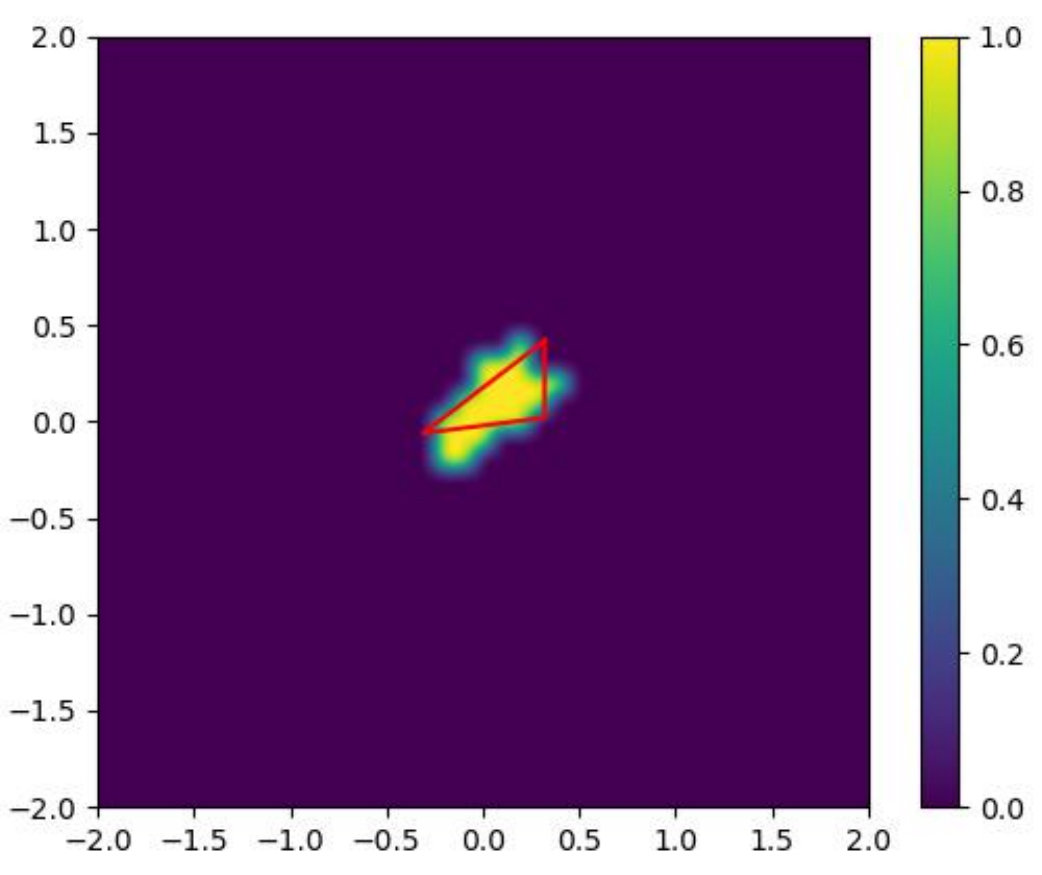}
	\end{minipage}
  \label{Results-Noisec}
  }

   \subfigure[0.0120, $0.8922$]
  {
	\begin{minipage}{3cm}
 	\centering
	\includegraphics[scale=0.18]{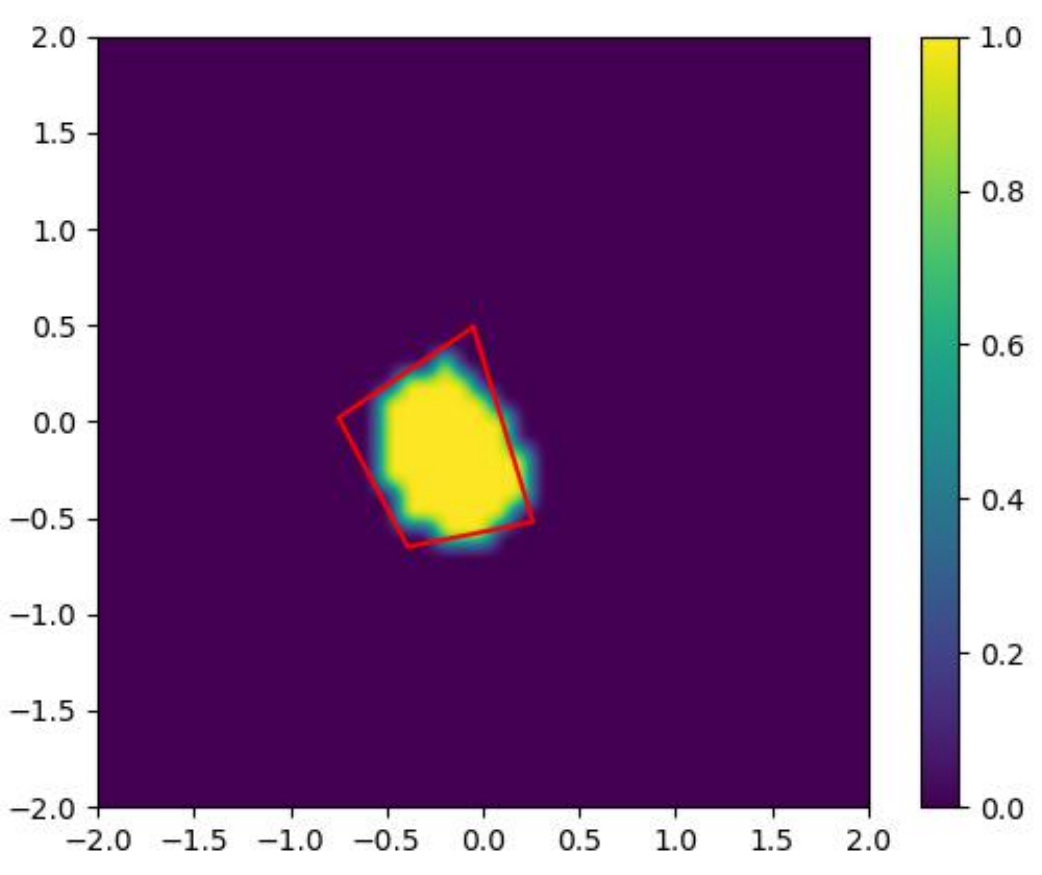}
	\end{minipage}
  }
  \subfigure[0.0160, $0.8485$]
  {
	\begin{minipage}{3cm}
 	
	\includegraphics[scale=0.18]{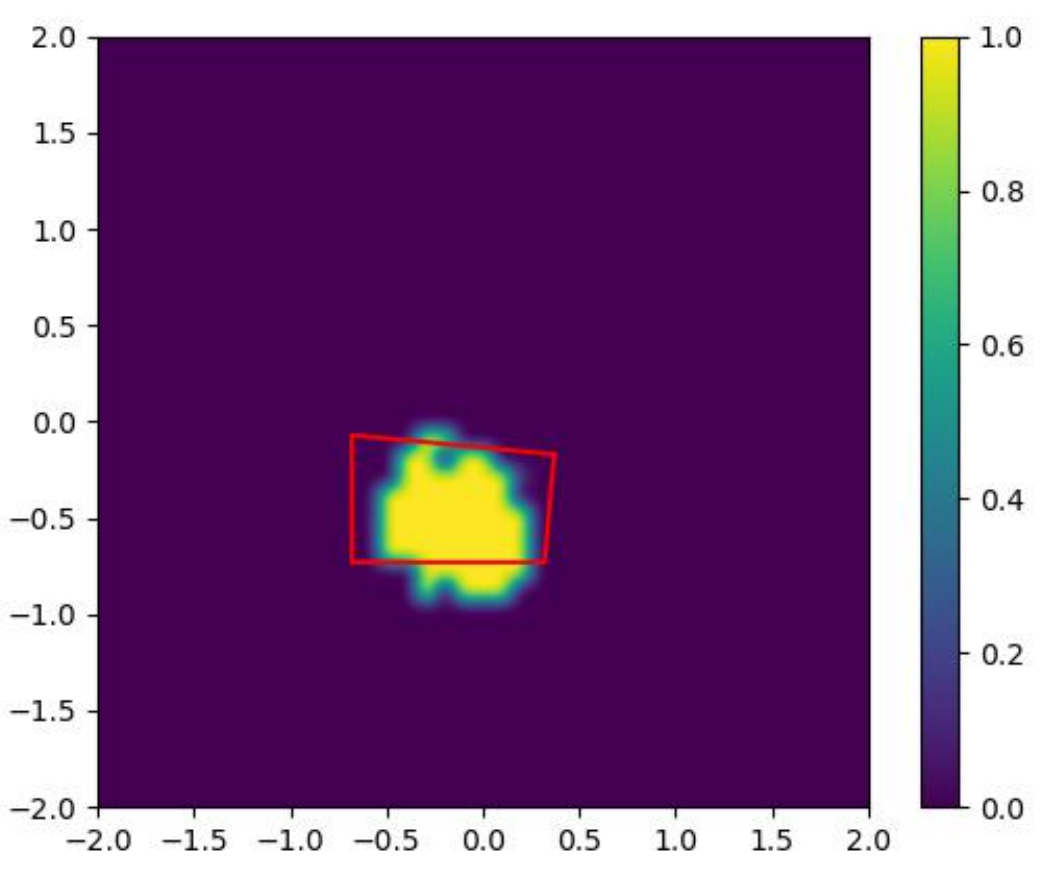}
	\end{minipage}
  }
   \subfigure[0.0054, $0.9365$]
  {
	\begin{minipage}{3cm}
 	\centering
	\includegraphics[scale=0.18]{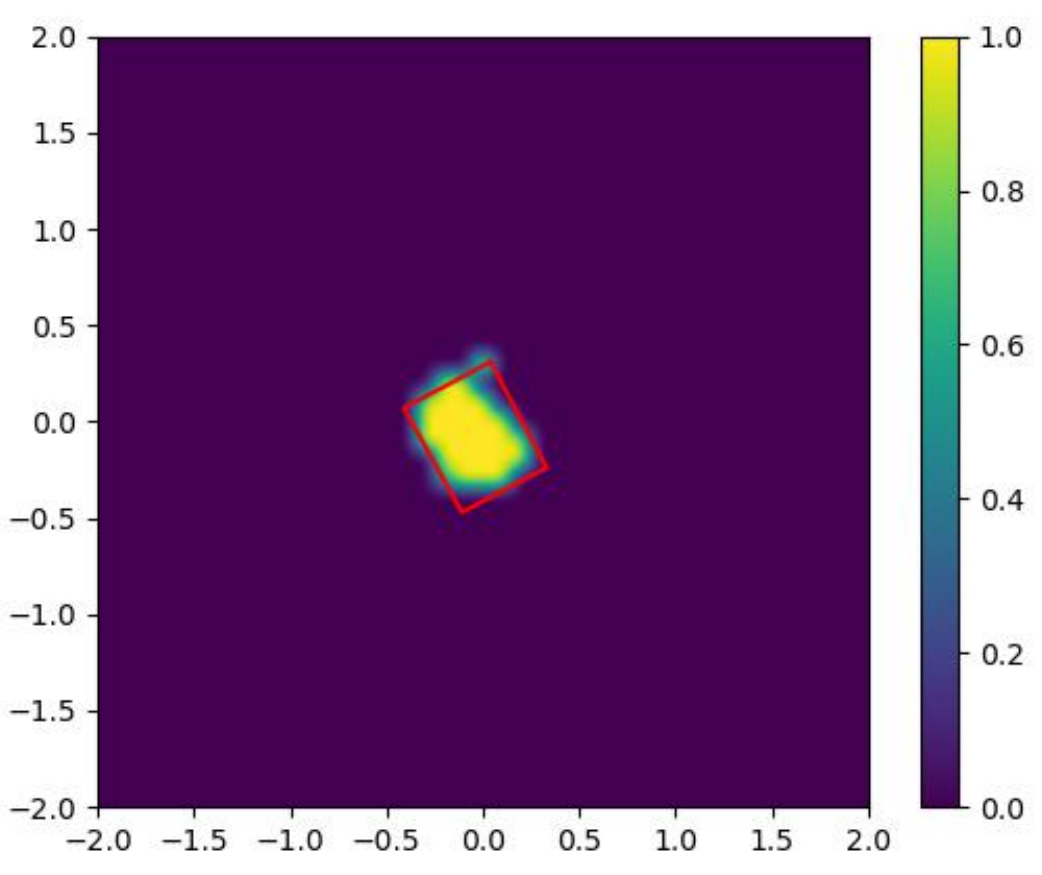}
	\end{minipage}
  }

   \subfigure[0.0095, $0.9037$]
  {
	\begin{minipage}{3cm}
 	\centering
	\includegraphics[scale=0.18]{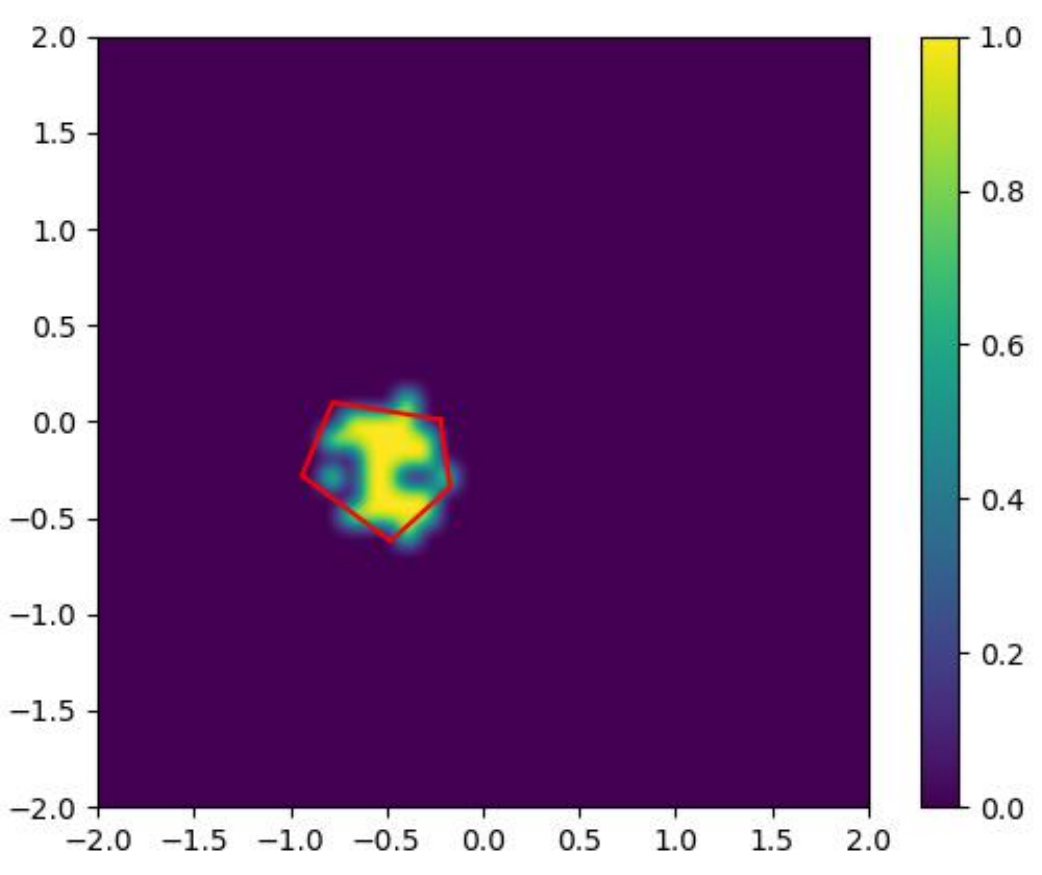}
	\end{minipage}
  }
  \subfigure[0.0149, $0.8444$]
  {
	\begin{minipage}{3cm}
 	
	\includegraphics[scale=0.18]{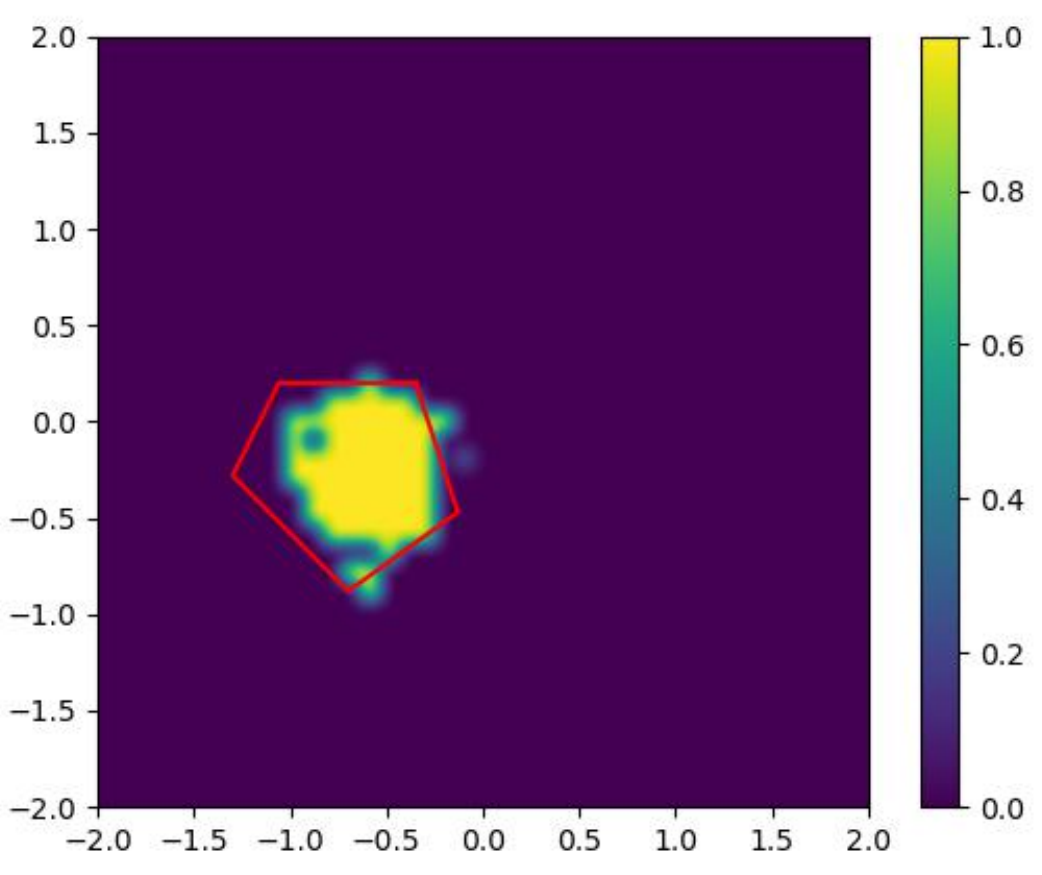}
	\end{minipage}
  }
   \subfigure[0.0224, $0.7857$]
  {
	\begin{minipage}{3cm}
 	\centering
	\includegraphics[scale=0.18]{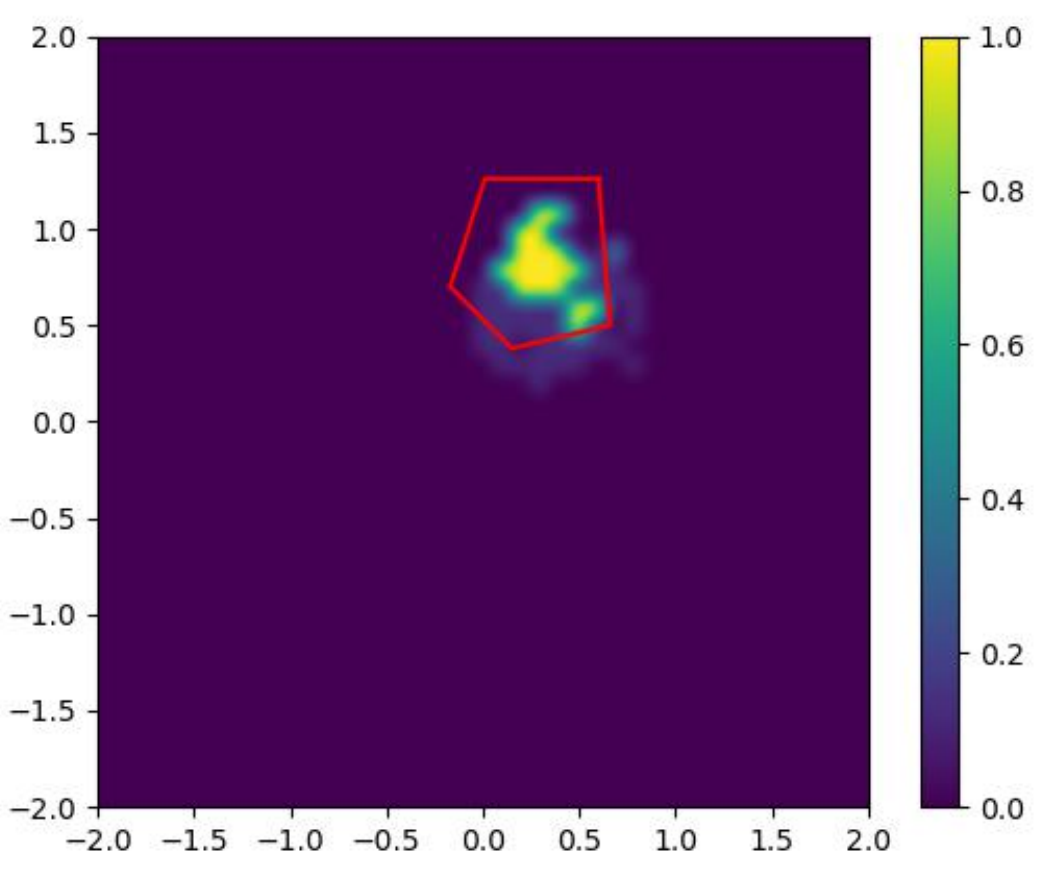}
	\end{minipage}
  }

   \subfigure[0.0089, $0.9021$]
  {
	\begin{minipage}{3cm}
 	\centering
	\includegraphics[scale=0.18]{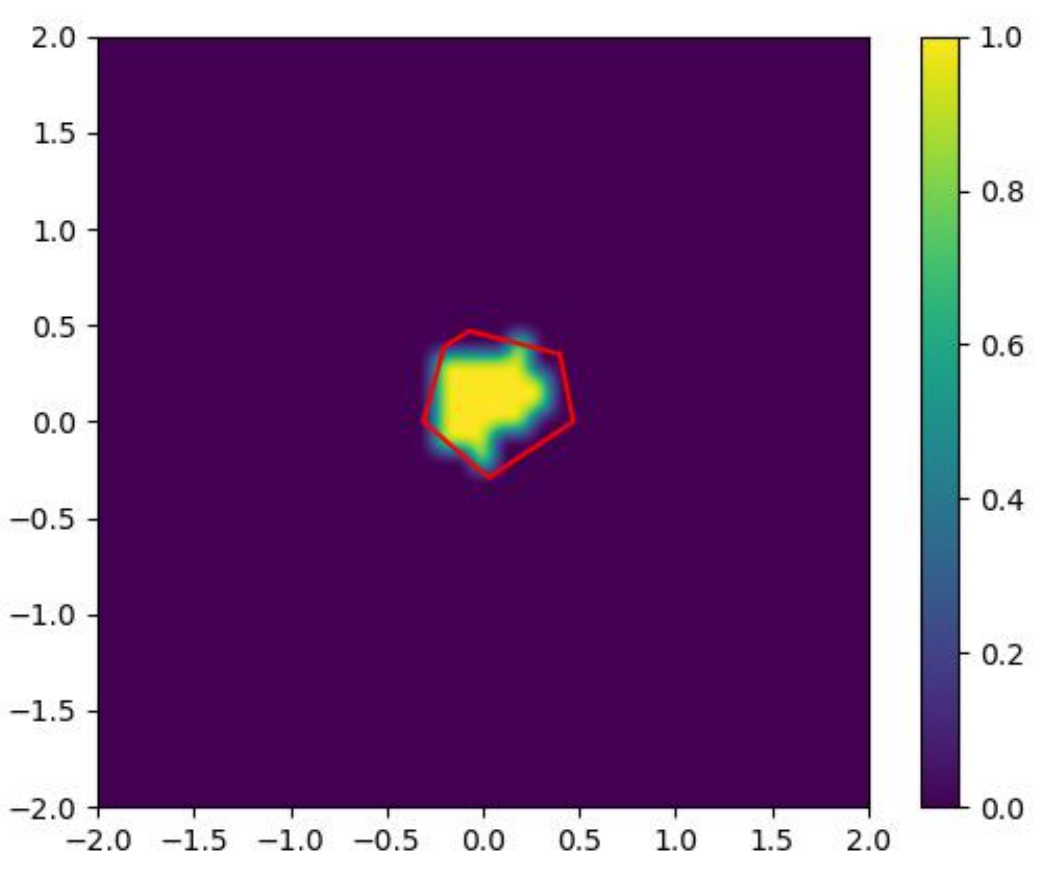}
	\end{minipage}
  }
  \subfigure[0.0232, $0.7851$]
  {
	\begin{minipage}{3cm}
 	
	\includegraphics[scale=0.18]{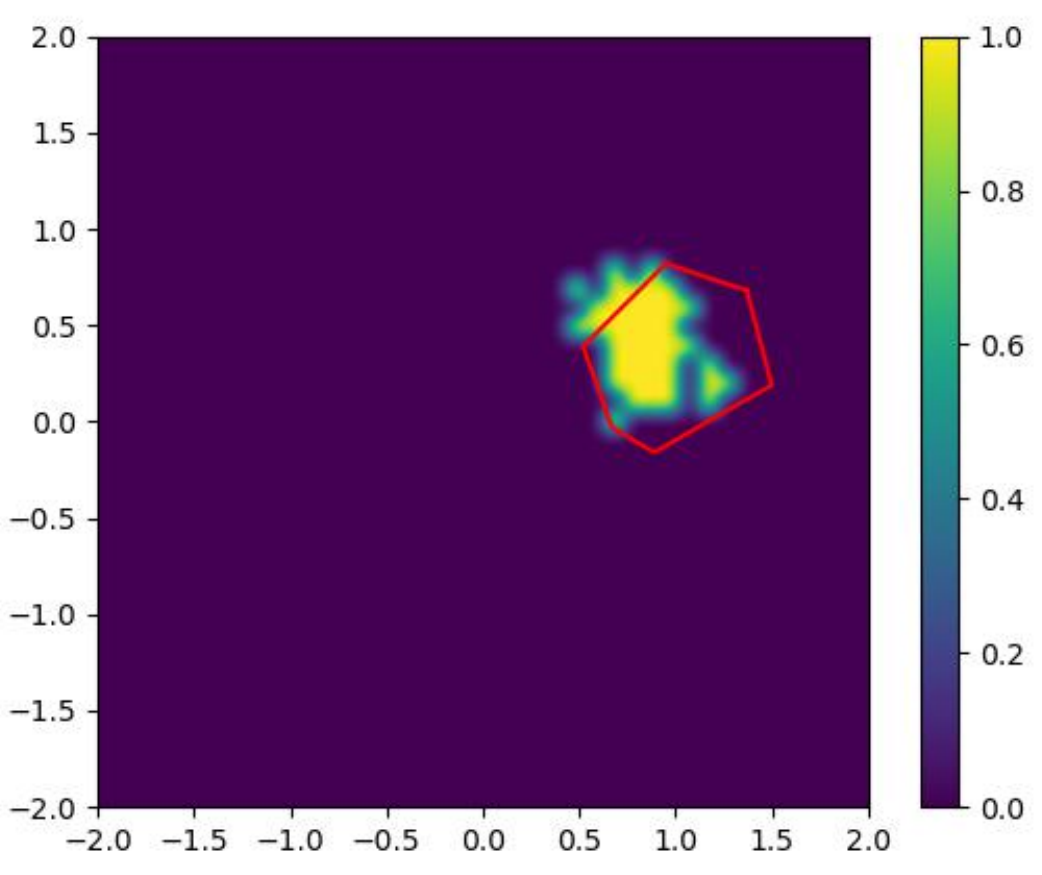}
	\end{minipage}
  }
   \subfigure[0.0202, $0.7926$]
  {
	\begin{minipage}{3cm}
 	\centering
	\includegraphics[scale=0.18]{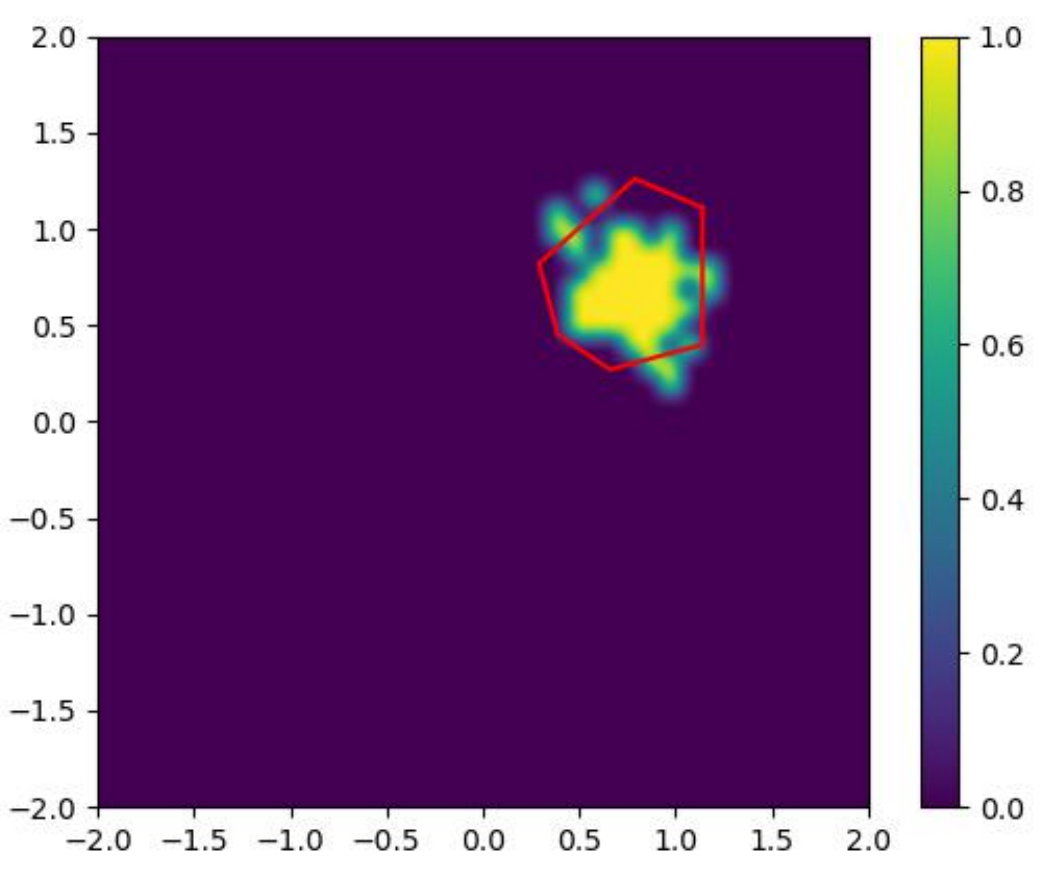}
	\end{minipage}
  }

 \caption{Some samples of the reconstructions of Scenario 2. The MSE (left) and SSIM (right) are recorded below each image.}
 \label{Results-Noise}
\end{figure}

\begin{figure}%[htbp]
\centering
 \subfigure[0.0225, $0.7671$]
  {
	\begin{minipage}{3cm}
 	\centering
	\includegraphics[scale=0.18]{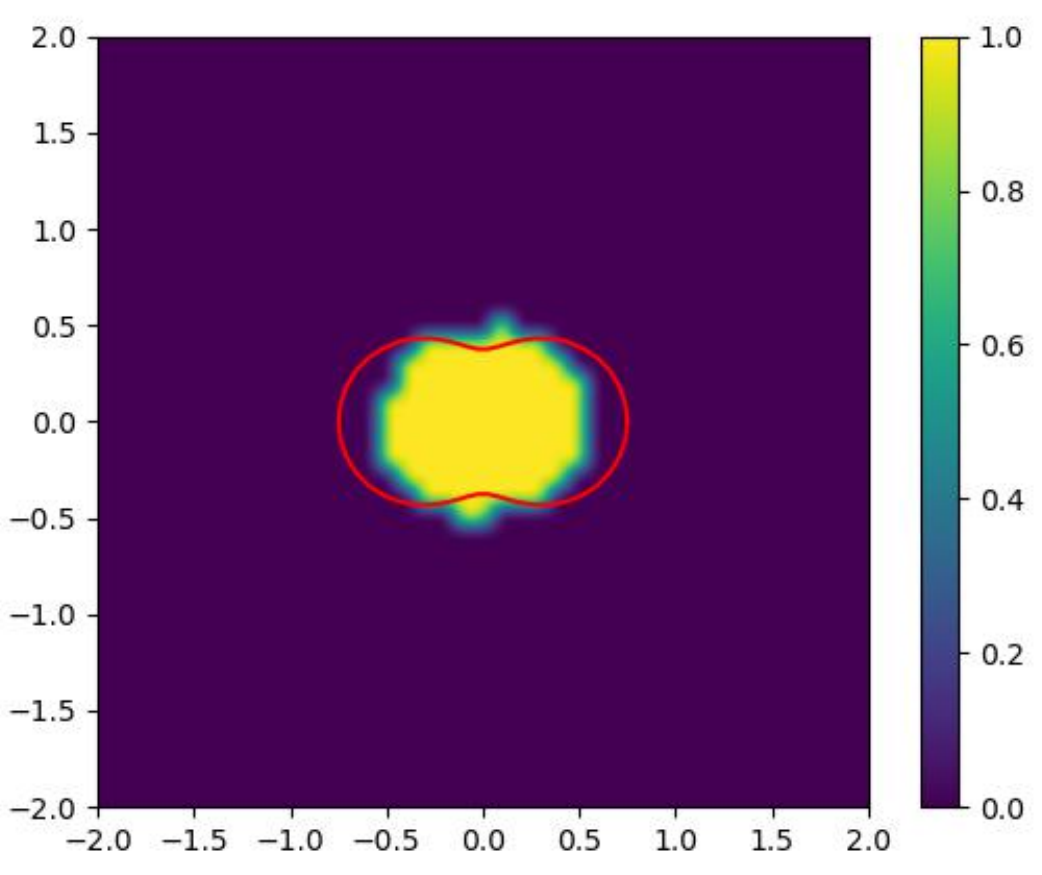}
	\end{minipage}
  }
  \subfigure[0.0105, $0.8850$]
  {
	\begin{minipage}{3cm}
 	
	\includegraphics[scale=0.18]{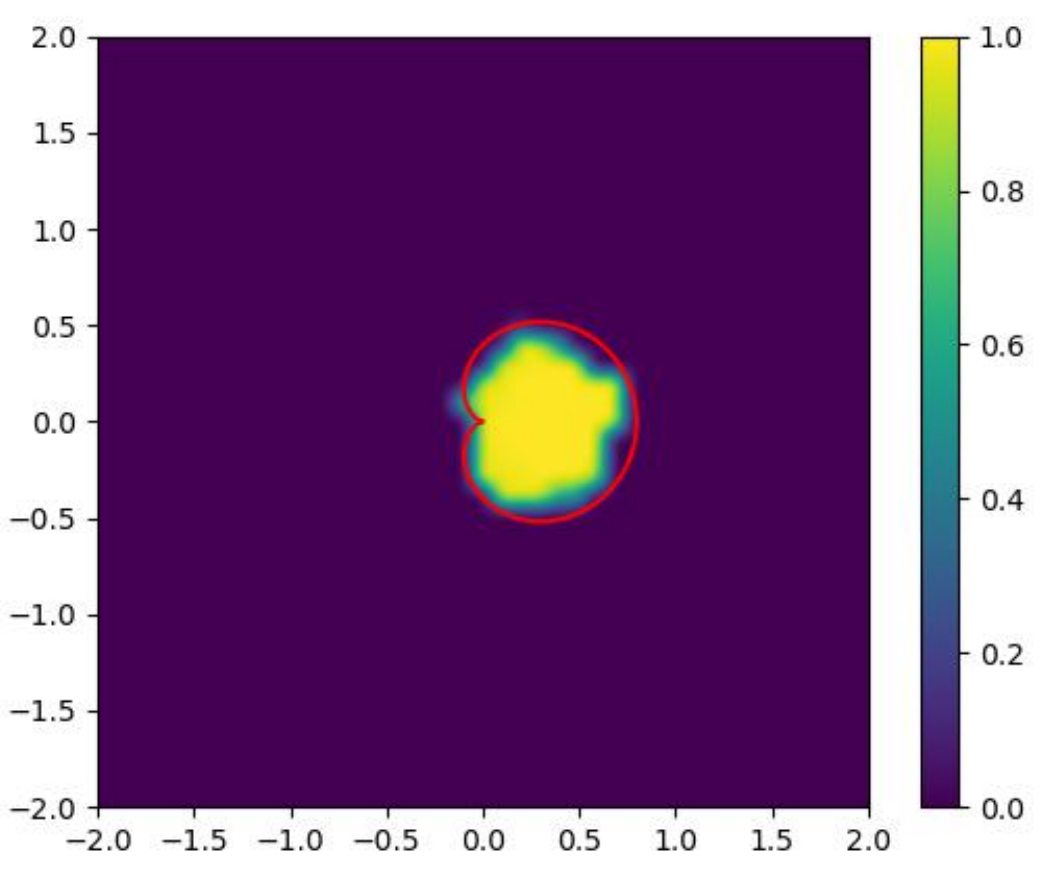}
	\end{minipage}
  }
   \subfigure[0.0071, $0.9232$]
  {
	\begin{minipage}{3cm}
 	\centering
	\includegraphics[scale=0.18]{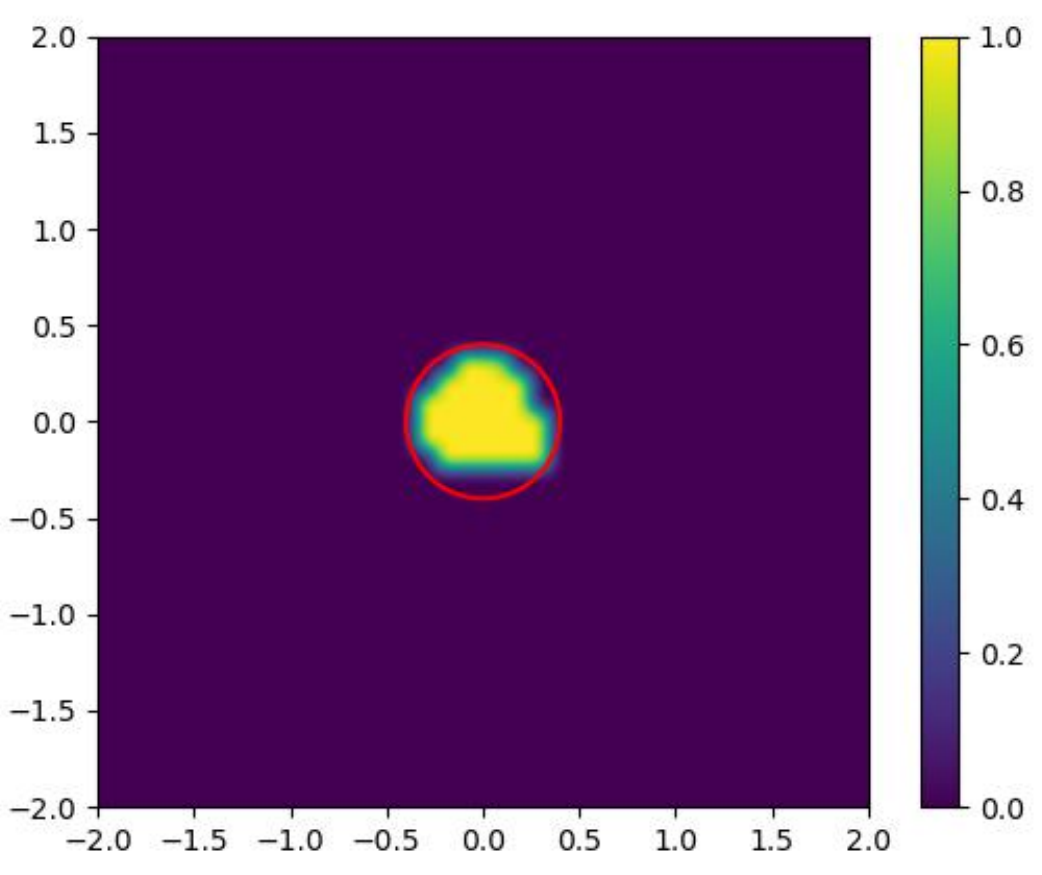}
	\end{minipage}
  }

   \subfigure[0.0303, $0.8091$]
  {
	\begin{minipage}{3cm}
 	\centering
	\includegraphics[scale=0.18]{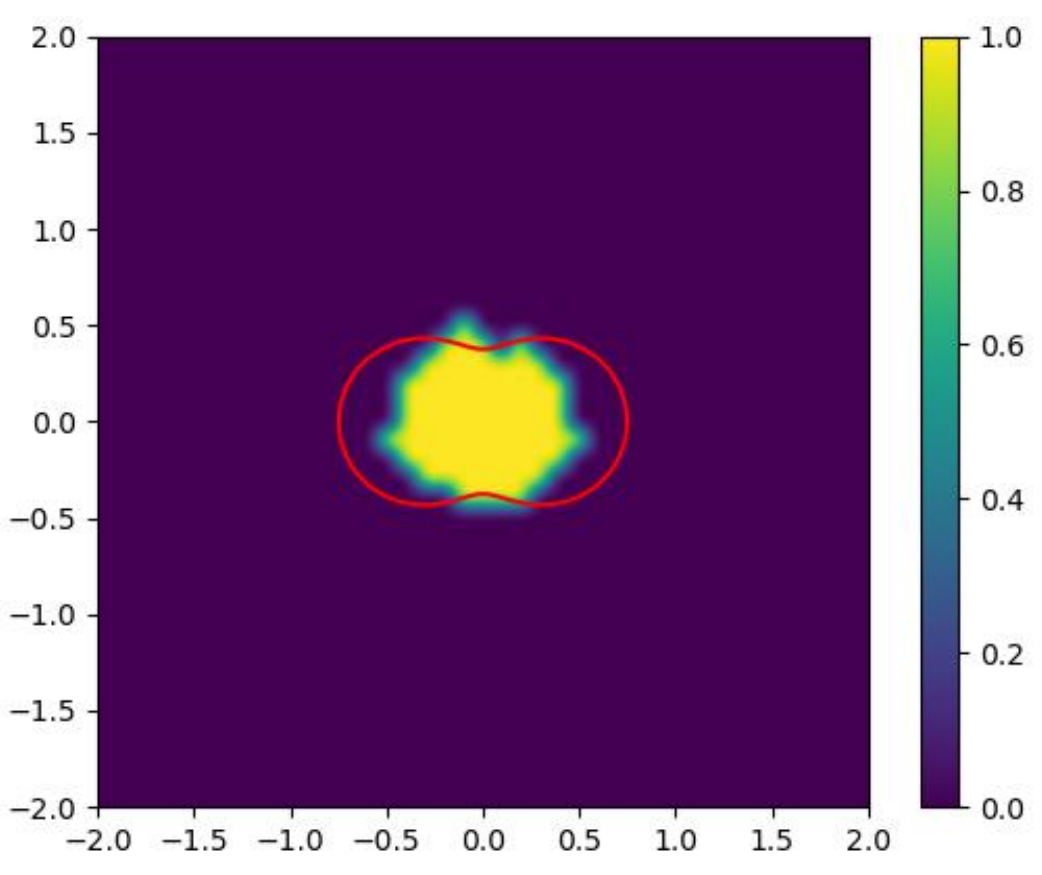}
	\end{minipage}
  }
  \subfigure[0.0234, $0.8015$]
  {
	\begin{minipage}{3cm}
 	
	\includegraphics[scale=0.18]{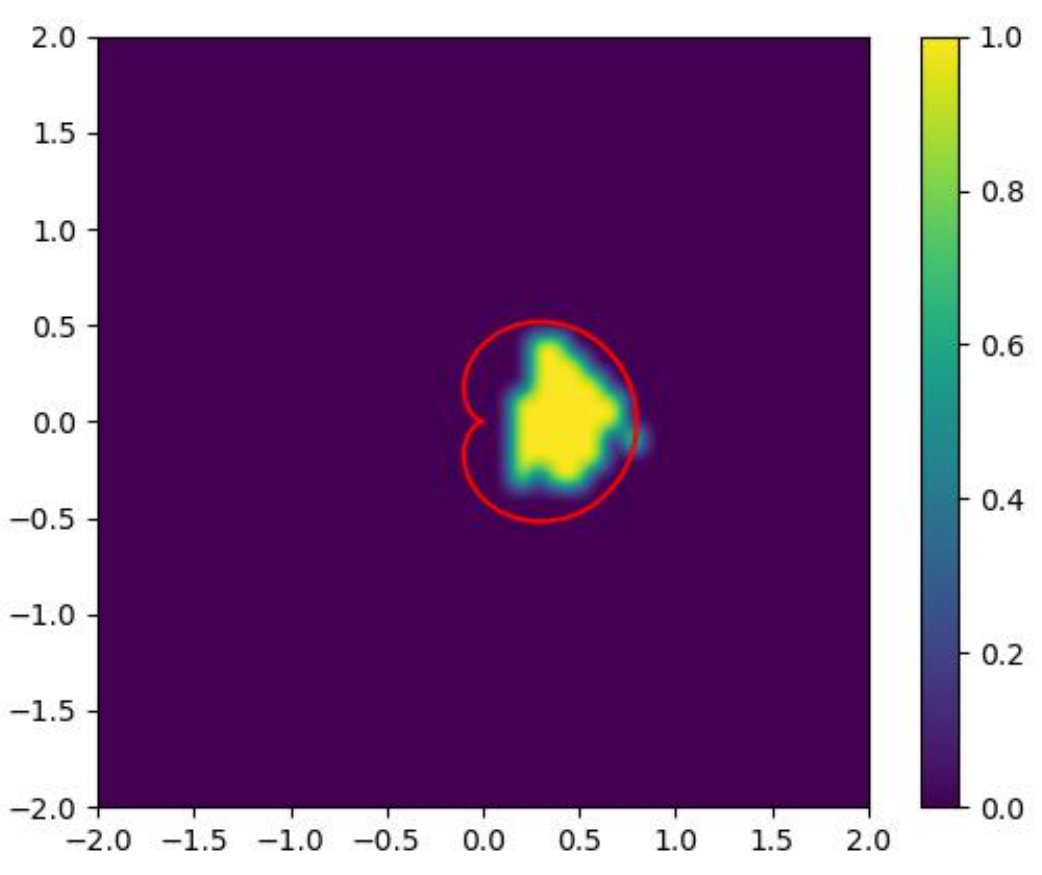}
	\end{minipage}
  }
   \subfigure[0.0101, $0.9343$]
  {
	\begin{minipage}{3cm}
 	\centering
	\includegraphics[scale=0.18]{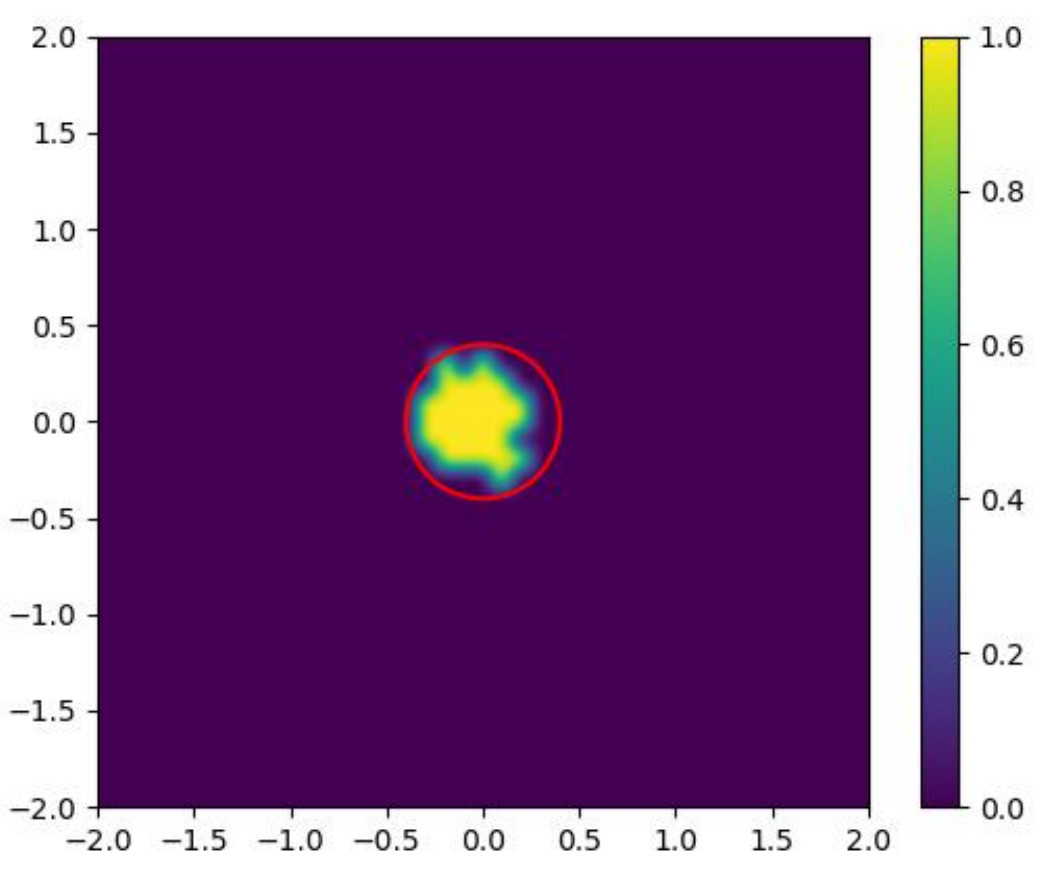}
	\end{minipage}
  }
 \caption{The numerical reconstructions of Scenario 3. The three pictures in the first row correspond to the noise-free case, while the three pictures in the second row correspond to the $3\%$ noise case. }
 \label{Results-Generalization}
\end{figure}

This subsection is devoted to various numerical experiments. We recall from $\S$\ref{S-TSS} that the final reconstruction obtained through the LRT is given by
\[
\Phi_{\Theta^*}(y^\delta(B)) = \sum_{\ell=1}^3 \Lambda_{\tilde\theta^*}(y^\delta(B))_\ell \,\Phi_{\theta_\ell^*}(y^\delta(B)),
\]
where $\Theta^*=(\theta_1^*,\dots,\theta_L^*,\tilde\theta^*)$. In all the experiments, the quality of the reconstruction is measured by the mean squared error 
\begin{equation}\label{eq:MSE}
   {\rm MSE}(B) = \frac{1}{N}\Vert\indic(B)- \Phi_{\Theta^*}(y^\delta(B))\Vert_2^2,
\end{equation}
which  measures the distance between the actual contrast $\indic(B)$ and the predicted result $\Phi_{\Theta^*}(y^\delta(B))$ of the LRT, for a domain $B$ in the test set. To further quantify the reconstruction quality, we employ the Structural Similarity Index Measure (SSIM)—an established metric used to gauge the similarity of two images \cite{Wang-2004-13}. This dimensionless quantity exhibits a theoretical upper bound of $1$ ($\text{SSIM}\in [-1,1]$) and demonstrates superior correlation to human visual perception, where values approaching $1$ indicate higher structural fidelity between reconstructed and ground-truth inclusions. In all the figures, the red polygons or curves refer to our reconstruction targets, and the corresponding MSE (left) and SSIM (right) are shown below each picture.
\begin{figure}%[htbp]
\centering
 \subfigure[0.0059, $0.9345$]
  {
	\begin{minipage}{3cm}
 	\centering
	\includegraphics[scale=0.18]{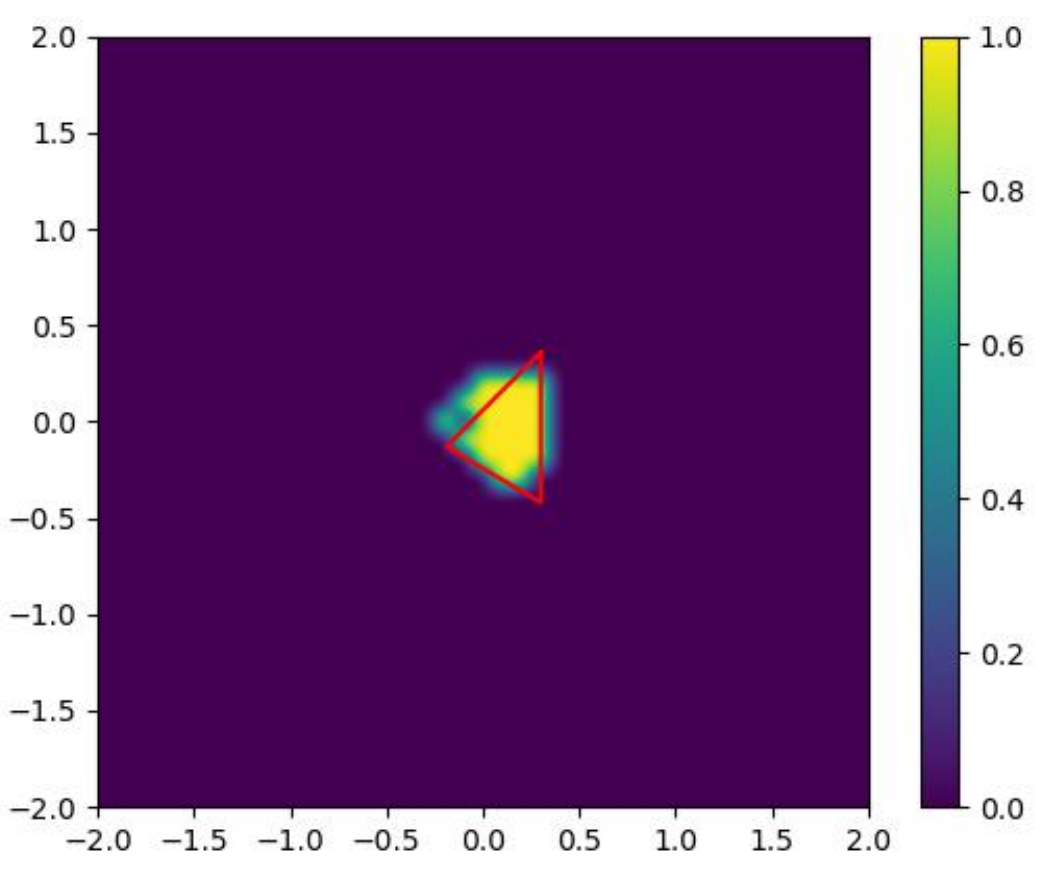}
	\end{minipage}
  }
  \subfigure[0.0089, $0.9070$]
  {
	\begin{minipage}{3cm}
 	
	\includegraphics[scale=0.18]{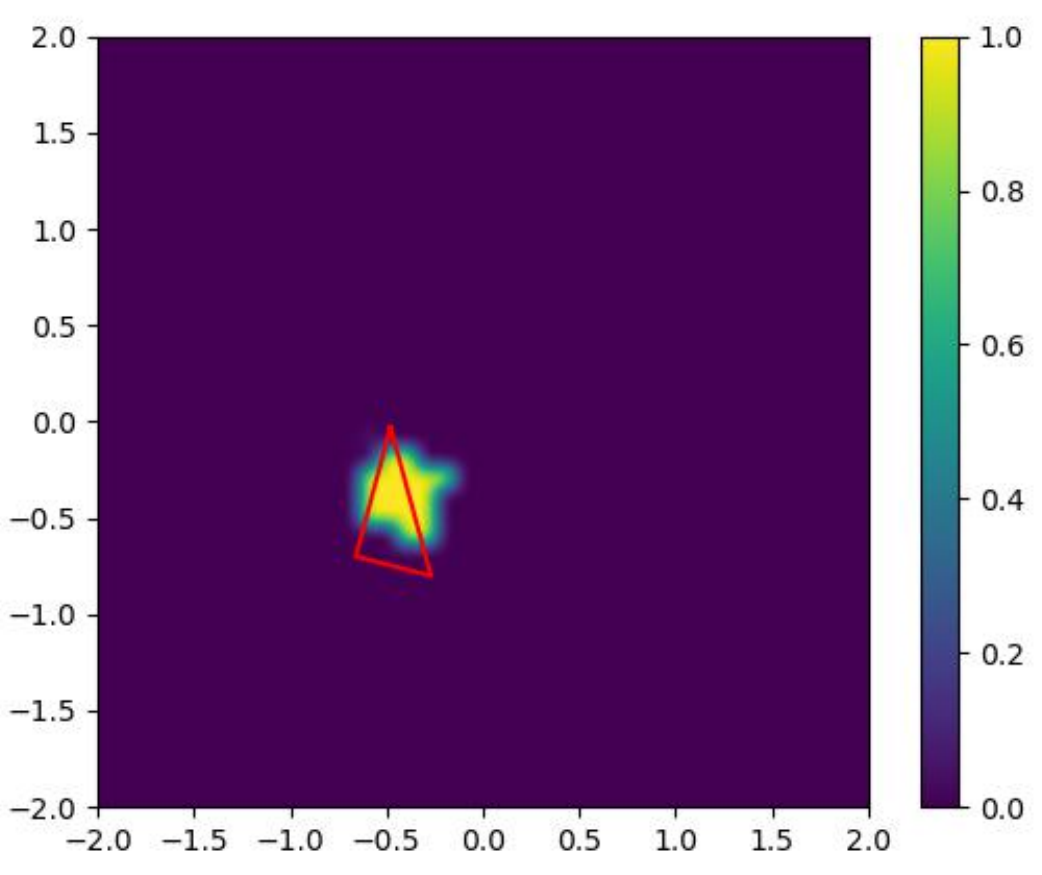}
	\end{minipage}
  }
   \subfigure[0.0071, $0.9218$]
  {
	\begin{minipage}{3cm}
 	\centering
	\includegraphics[scale=0.18]{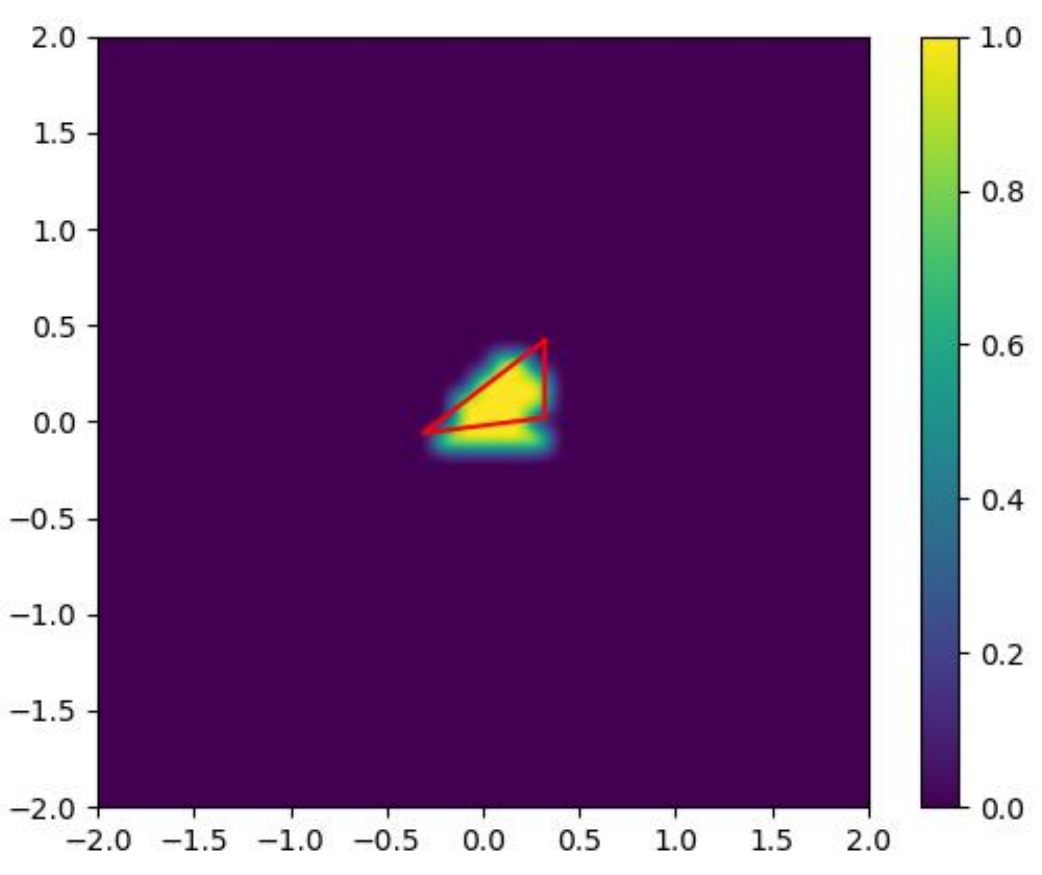}
	\end{minipage}
  \label{Results-Noise-Freec}
  }

   \subfigure[0.0111, $0.8862$]
  {
	\begin{minipage}{3cm}
 	\centering
	\includegraphics[scale=0.18]{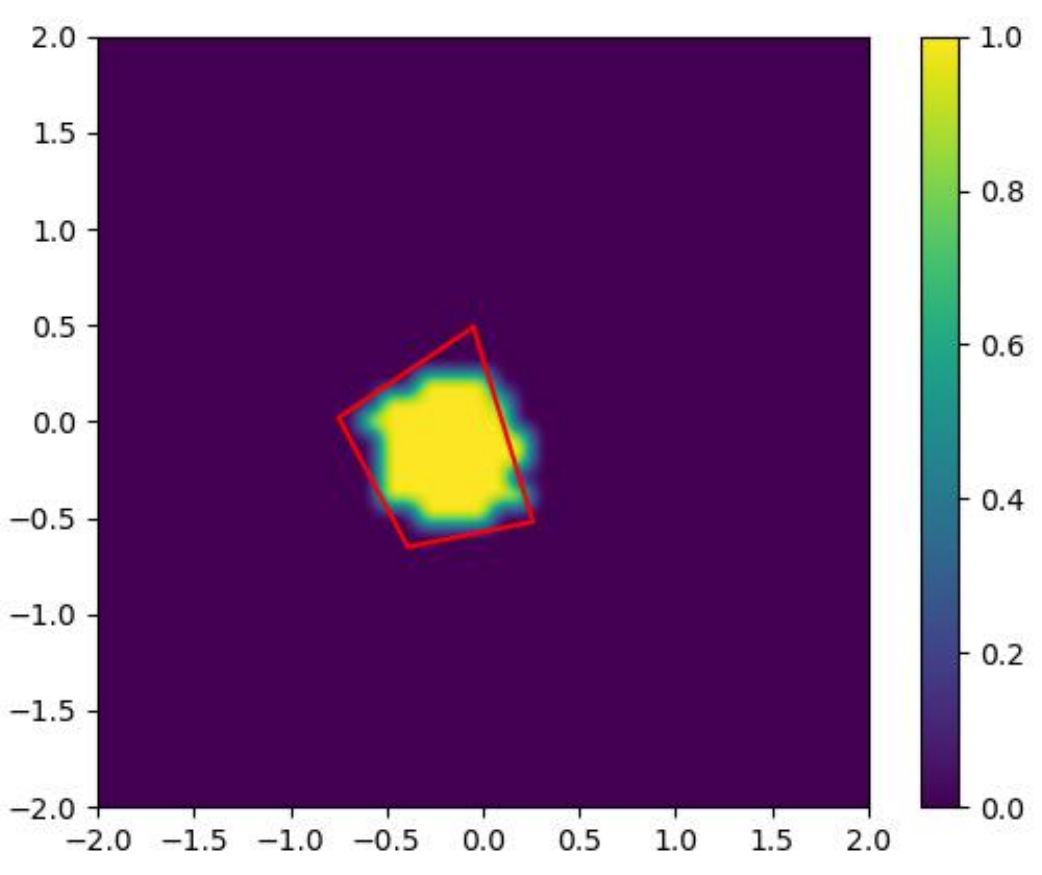}
	\end{minipage}
  }
  \subfigure[0.0158, $0.8640$]
  {
	\begin{minipage}{3cm}
 	
	\includegraphics[scale=0.18]{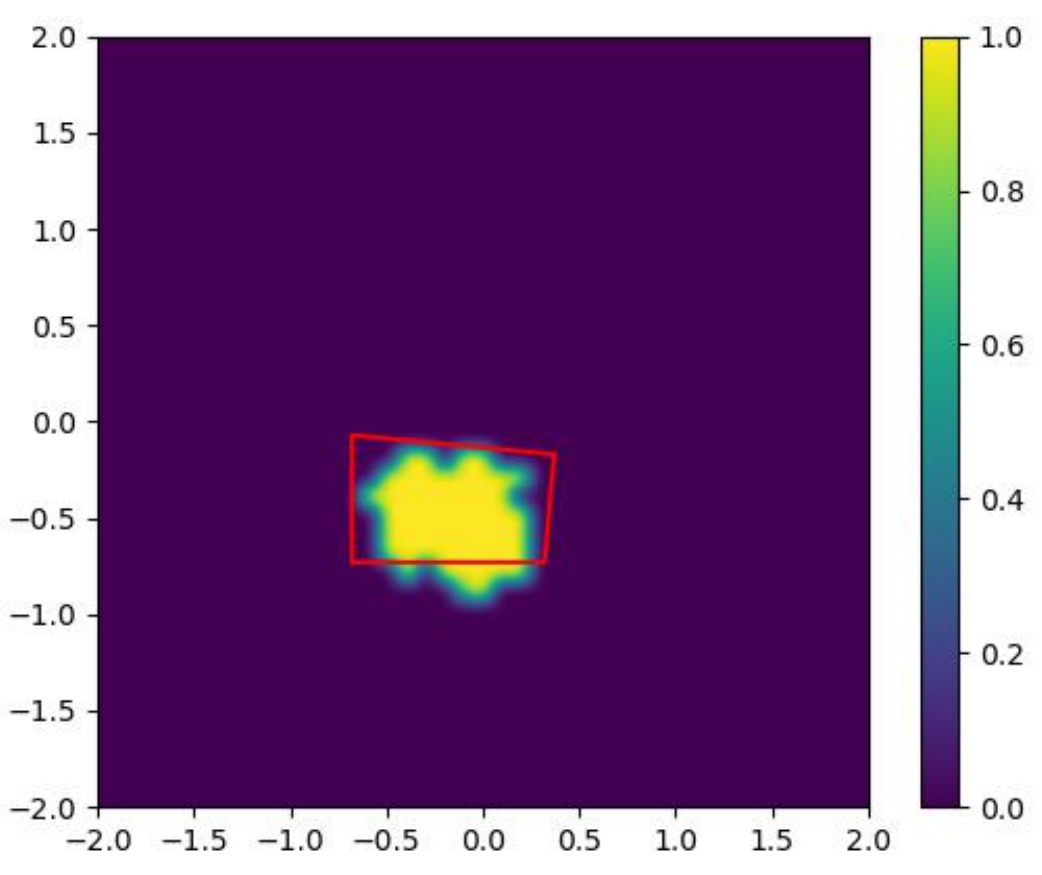}
	\end{minipage}
  }
   \subfigure[0.0073, $0.9223$]
  {
	\begin{minipage}{3cm}
 	\centering
	\includegraphics[scale=0.18]{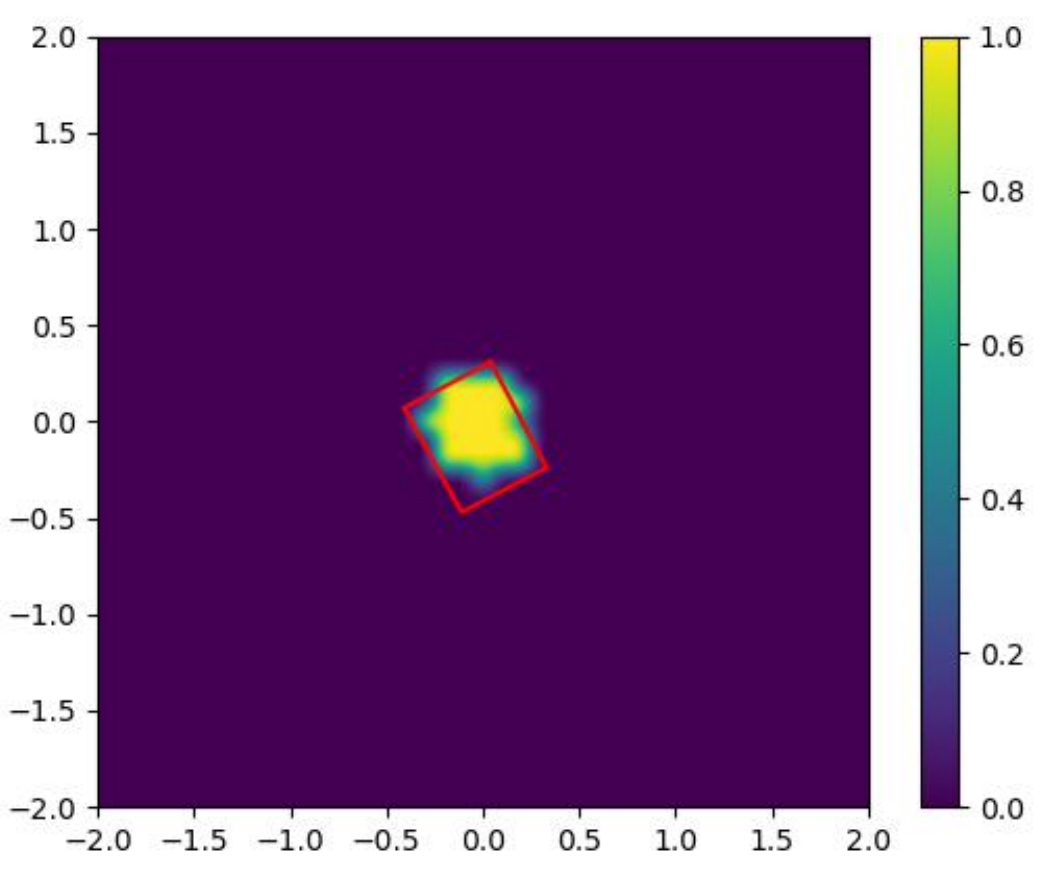}
	\end{minipage}
  }

   \subfigure[0.0060, $0.9320$]
  {
	\begin{minipage}{3cm}
 	\centering
	\includegraphics[scale=0.18]{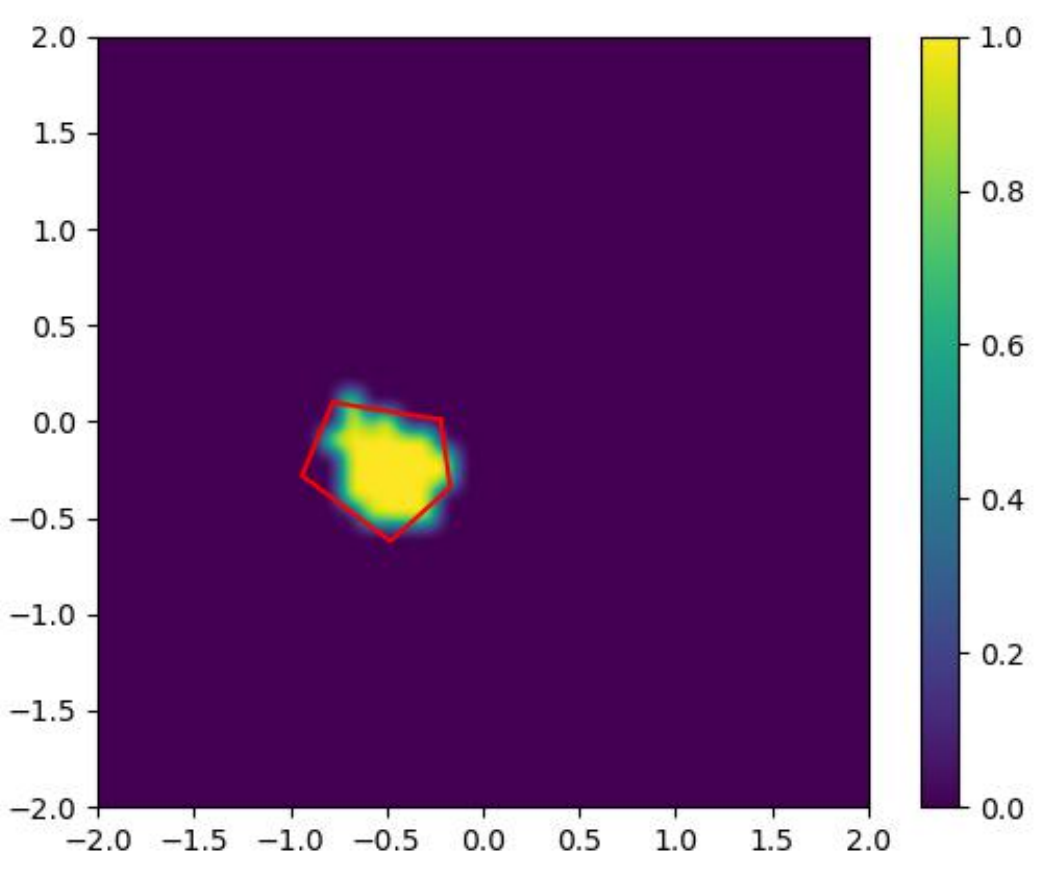}
	\end{minipage}
  }
  \subfigure[0.0092, $0.9032$]
  {
	\begin{minipage}{3cm}
 	
	\includegraphics[scale=0.18]{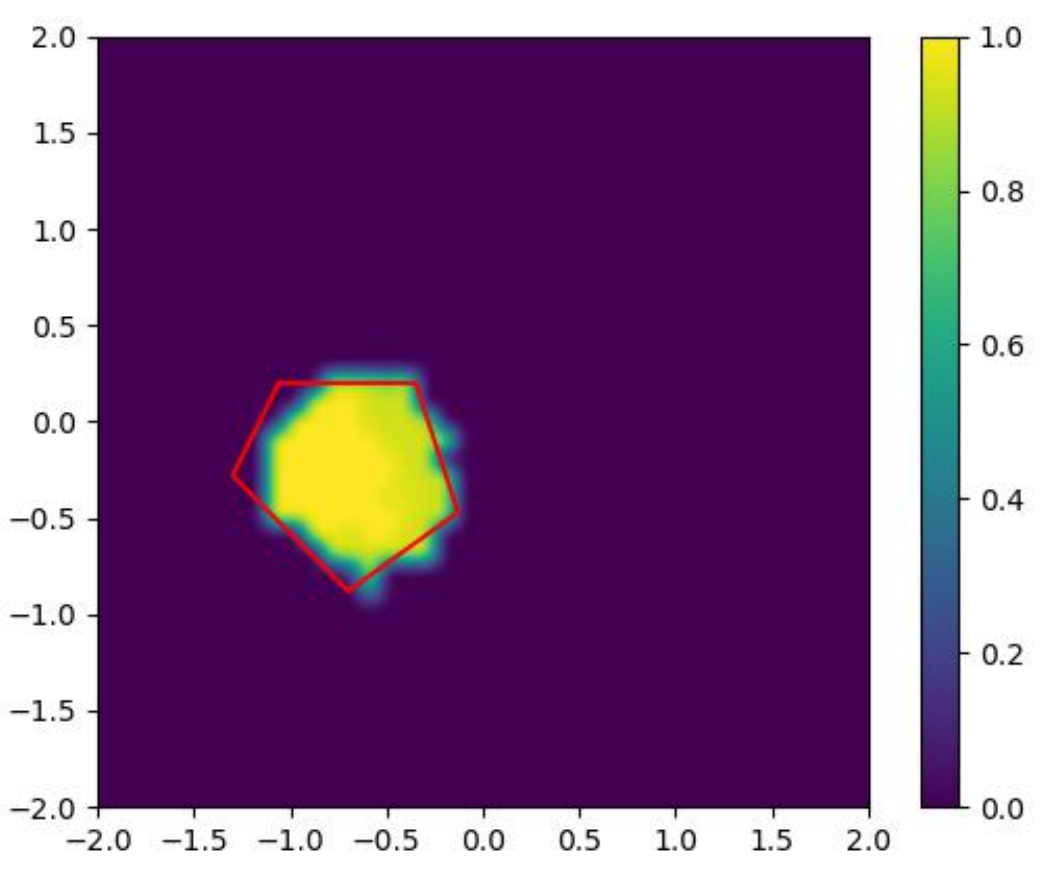}
	\end{minipage}
  }
   \subfigure[0.0093, $0.9160$]
  {
	\begin{minipage}{3cm}
 	\centering
	\includegraphics[scale=0.18]{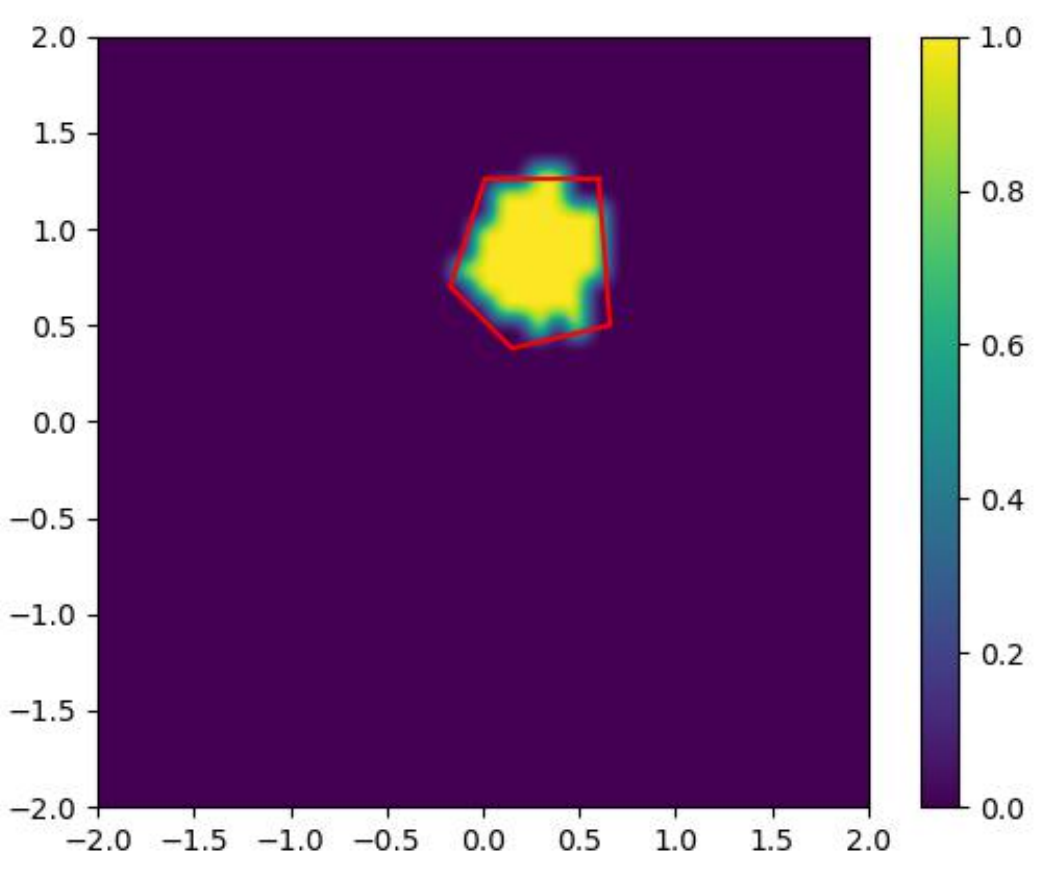}
	\end{minipage}
  }

   \subfigure[0.0059, $0.9241$]
  {
	\begin{minipage}{3cm}
 	\centering
	\includegraphics[scale=0.18]{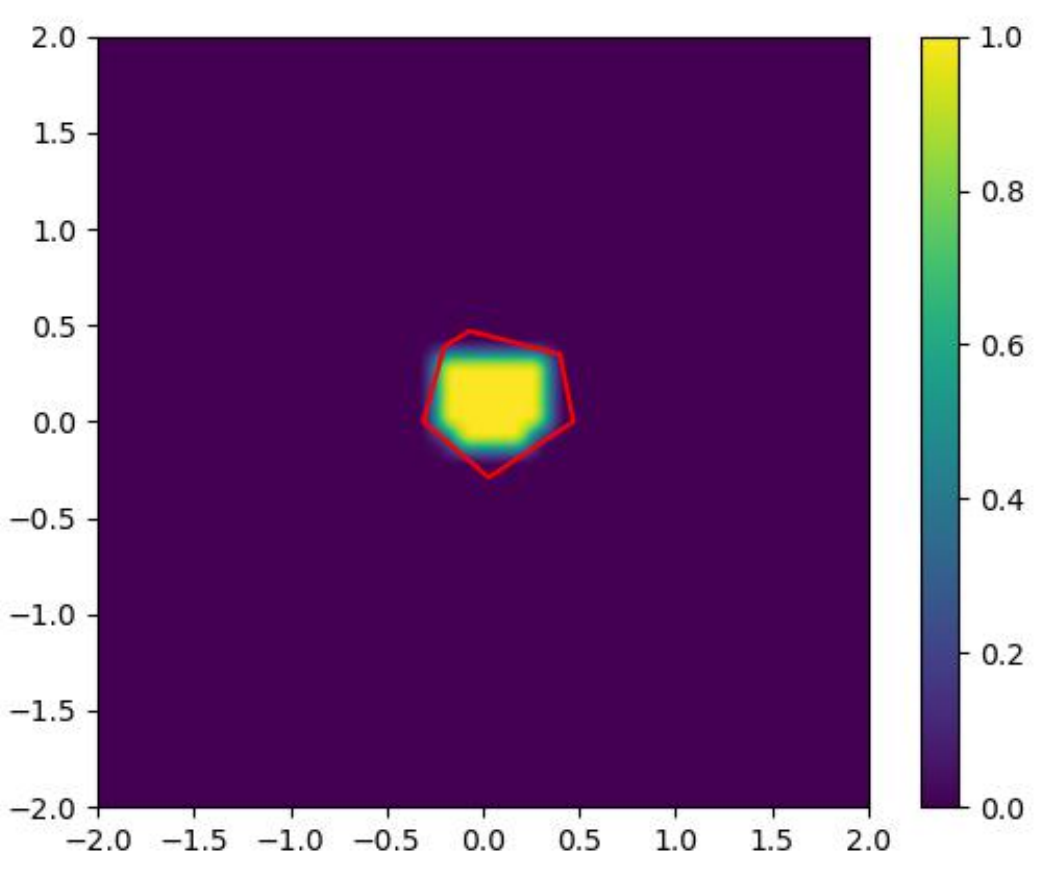}
	\end{minipage}
  }
  \subfigure[0.0120, $0.8984$]
  {
	\begin{minipage}{3cm}
 	
	\includegraphics[scale=0.18]{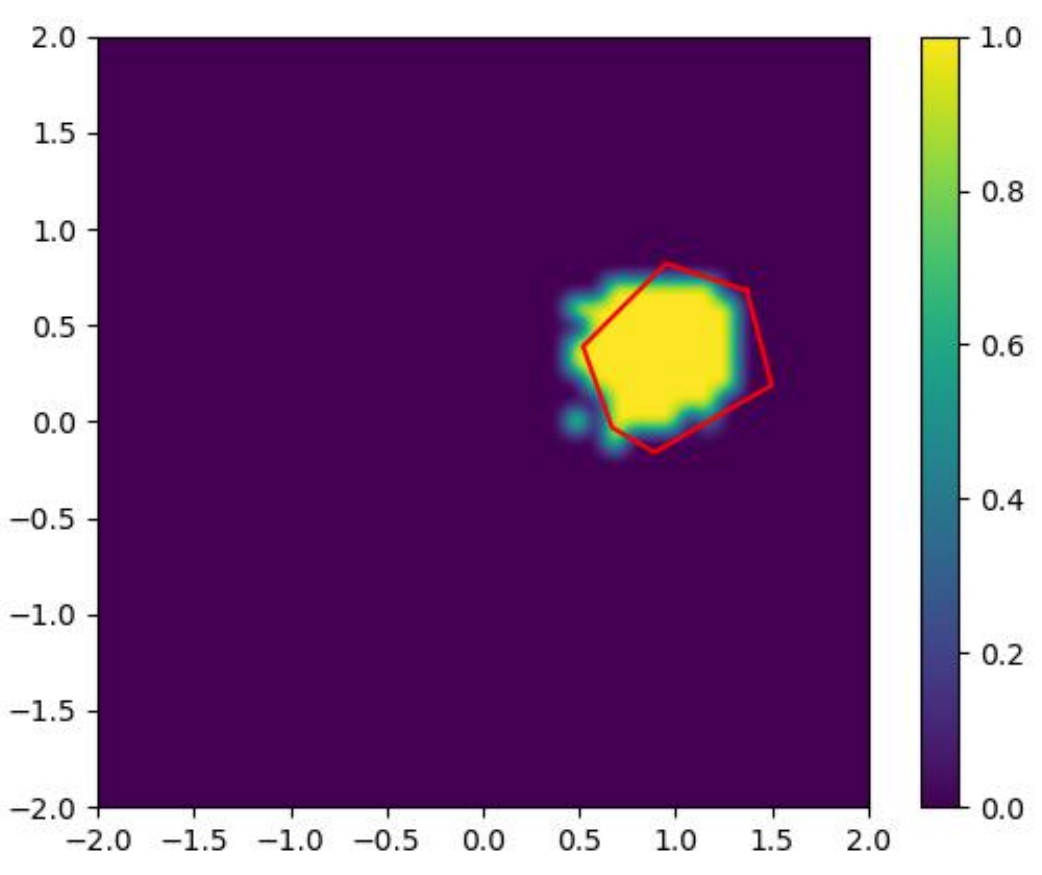}
	\end{minipage}
  }
   \subfigure[0.0089, $0.9098$]
  {
	\begin{minipage}{3cm}
 	\centering
	\includegraphics[scale=0.18]{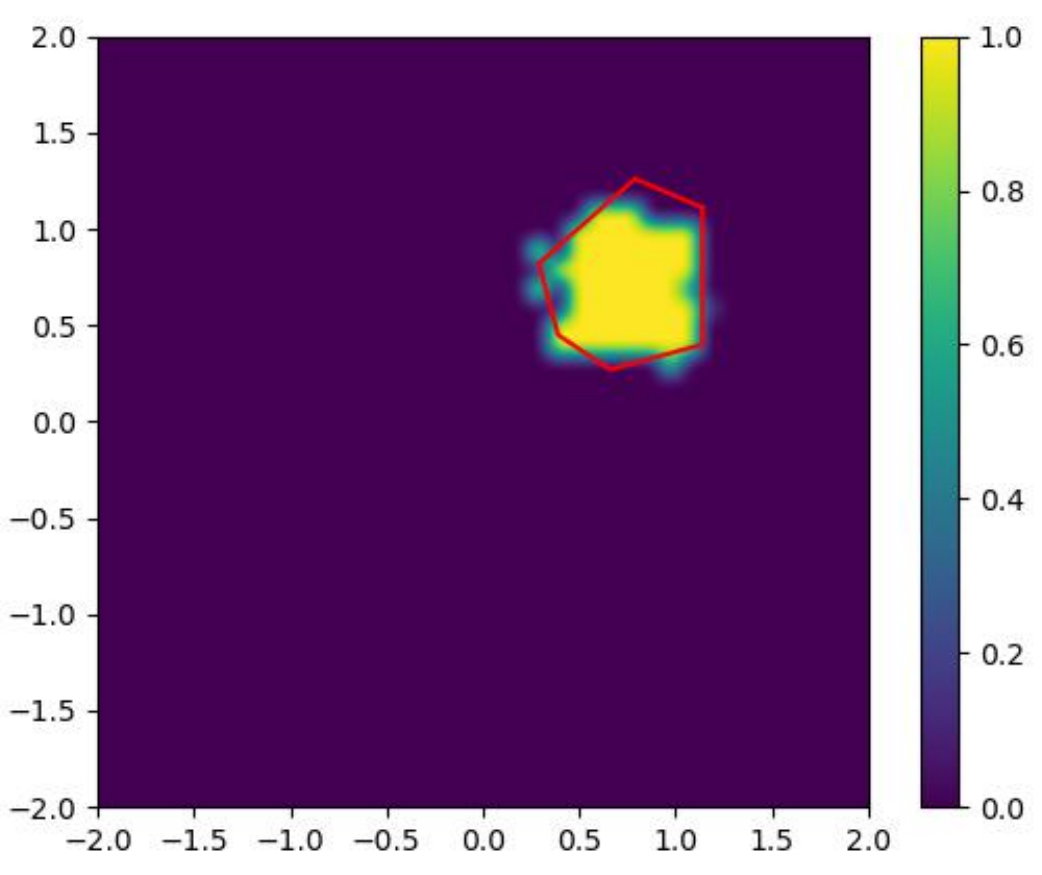}
	\end{minipage}
  }

 \caption{Some samples of the reconstructions of Scenario 1. The MSE (left) and SSIM (right) are recorded below each image.}
 \label{Results-Noise-Free}
\end{figure}
We consider three different scenarios.
\begin{itemize}
    \item \textbf{Scenario 1: Evaluate the performance of the LRT with  noise-free measurements.}  
We consider the LRT method in the noiseless case, namely, with 
 $\delta = 0$. Some reconstructions taken from the test set are shown in Figure~\ref{Results-Noise-Free}. The distribution of the mean squared errors for all the samples in the test set is shown in Figure~\ref{Histograma}.

     \item \textbf{Scenario 2: Evaluate the performance of the LRT with  noisy measurements.}
We consider the same setting of Scenario 1, but with noisy measurements, with noise level $\delta = 3\%$.
     Some reconstructions are shown in Figure~\ref{Results-Noise}, and the distribution of the mean squared errors in Figure~\ref{Histogramb}.

     \item \textbf{Scenario 3: Evaluate the generalization performance of the LRT.}
    We consider $3$ inclusions that are not within the class of our training set. More specifically, the samples in the training set are all convex polygons (with the maximum number of sides being $8$), while these three inclusions either are non-convex or have $C^2$ boundaries. Notably, the diameter of the peanut-shaped inclusion is $1.5$, exceeding the maximum diameter of all the samples in the training set, which is 1.2 according to \eqref{eq:Sample}. To test the generalization performance of the LRT, these three inclusions are used as the test samples in Scenario 1 and Scenario 2. Here, we highlight that the selected samples are mainly used to test the generalization of LRT in terms of different shapes instead of positions. The inclusions and the reconstructions are depicted in Figure~\ref{Results-Generalization}.
\end{itemize}

From Figures~\ref{Results-Noise-Free} and \ref{Results-Noise}, we can observe that the LRT can reconstruct the inclusion with high accuracy in both the noise-free and noisy measurement cases. In particular, it recovers the shape and the position of the inclusion even when the target is small and far away from the boundary; see, for instance, Figure \ref{Results-Noise-Freec} and Figure \ref{Results-Noisec}.  It is worth recalling that the recovery of inclusions from boundary measurements is a severely ill-posed problem, for which only logarithmic-type stability holds \cite{rondi-1999,Alessandrini-2001-176,dicristo-rondi-2003}. As a consequence, the reconstructions obtained typically have very low resolution, unless specific methods designed for polygonal inclusions are used \cite{beretta-etal-2018}. Thus, in general, reconstructing the exact shape of the inclusion is very difficult, and the reconstructions shown in Figures~\ref{Results-Noise-Free} and \ref{Results-Noise} are to be considered very accurate in this context. 

As shown in Figures~\ref{Histograma} and \ref{Histogramb}, the loss values for the samples in the test set are primarily distributed within the interval $[0, 0.02]$, which further illustrates the effectiveness of the LRT. Furthermore, the numerical results shown in Figures \ref{Results-Generalization} suggest that the LRT exhibits excellent generalization property. It is noteworthy that only one boundary measurement is employed to realize the reconstruction. 

\subsection{Comparisons with other schemes}\label{sub:other-strategy}
To further illustrate the effectiveness of the LRT, we conduct some numerical experiments to compare the LRT with other numerical strategies. %For brevity, only the noise-free measurements are considered in this subsection. For completeness, the corresponding simulations in the noisy case are included in Appendix~\ref{app:noisy}.

\subsubsection{Deterministic RT}\label{subsub:RT}
\begin{figure}%[htbp]
\centering
 \subfigure[0.0059, $0.9345$]
  {
	\begin{minipage}{3cm}
 	\centering
	\includegraphics[scale=0.18]{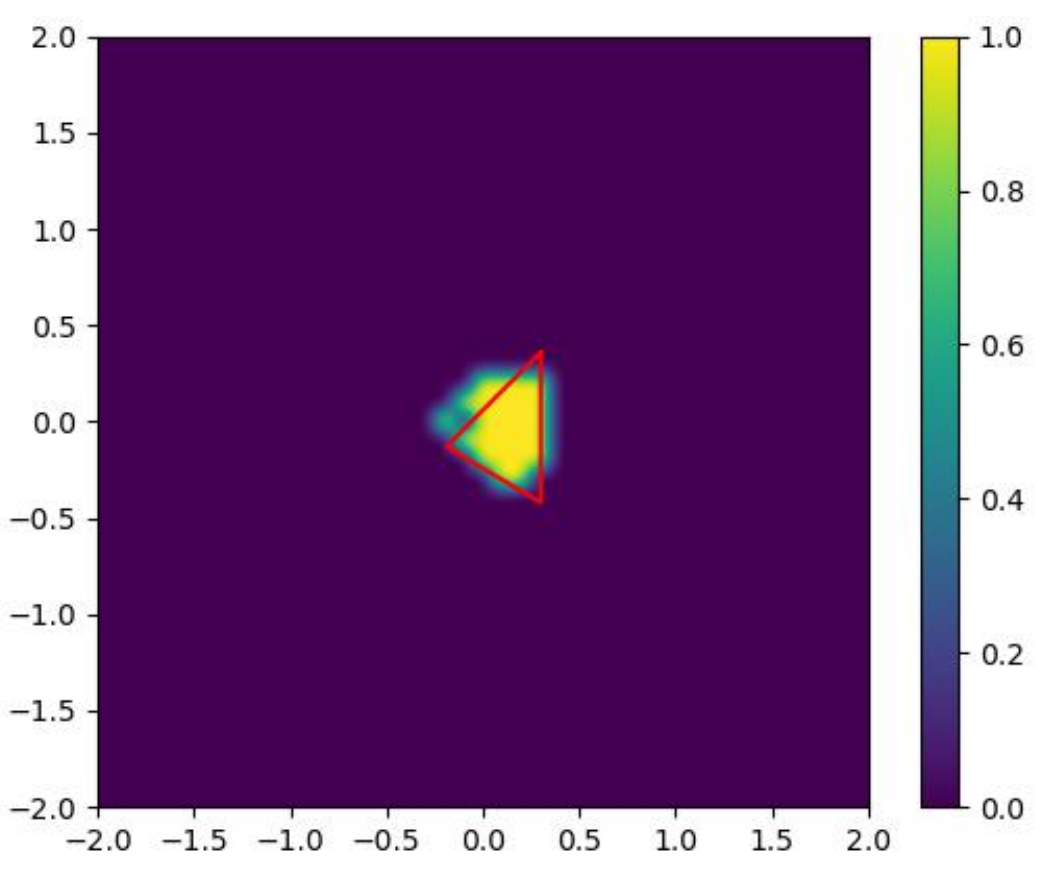}
	\end{minipage}
  }
  \subfigure[0.0158, $0.8640$]
  {
	\begin{minipage}{3cm}
 	
	\includegraphics[scale=0.18]{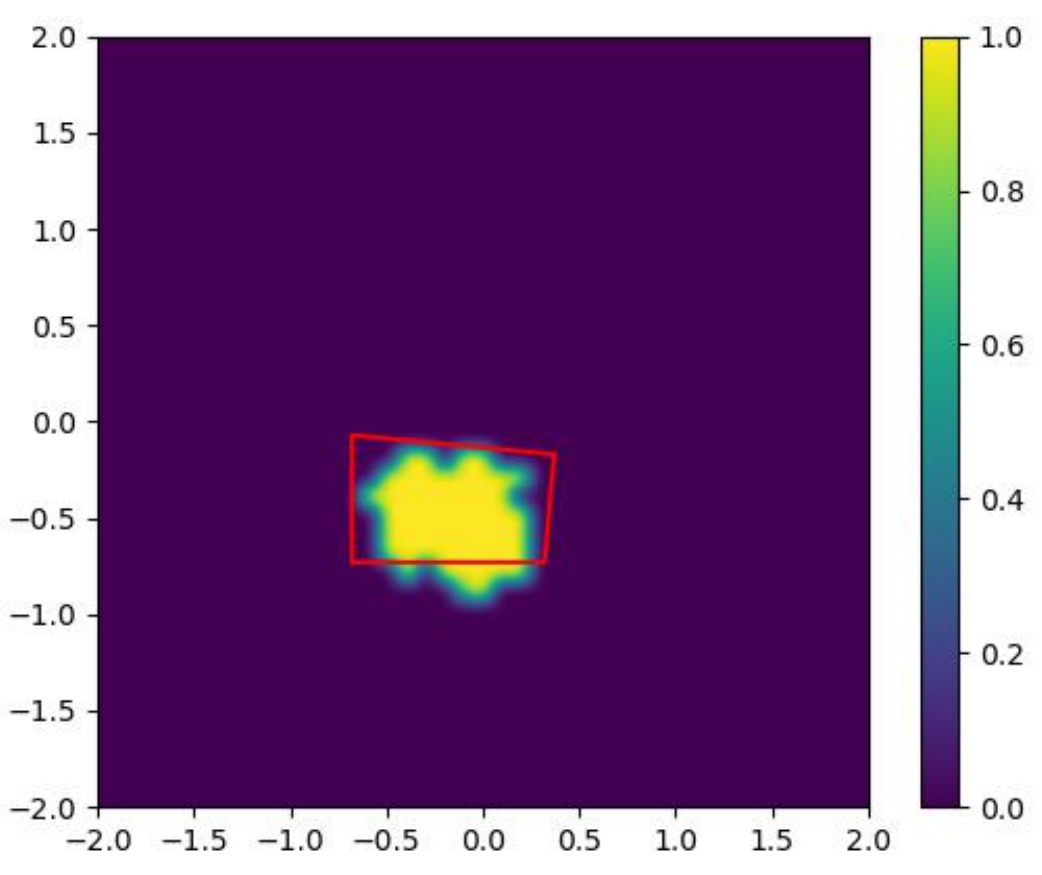}
	\end{minipage}
  }
   \subfigure[0.0071, $0.9241$]
  {
	\begin{minipage}{3cm}
 	\centering
	\includegraphics[scale=0.18]{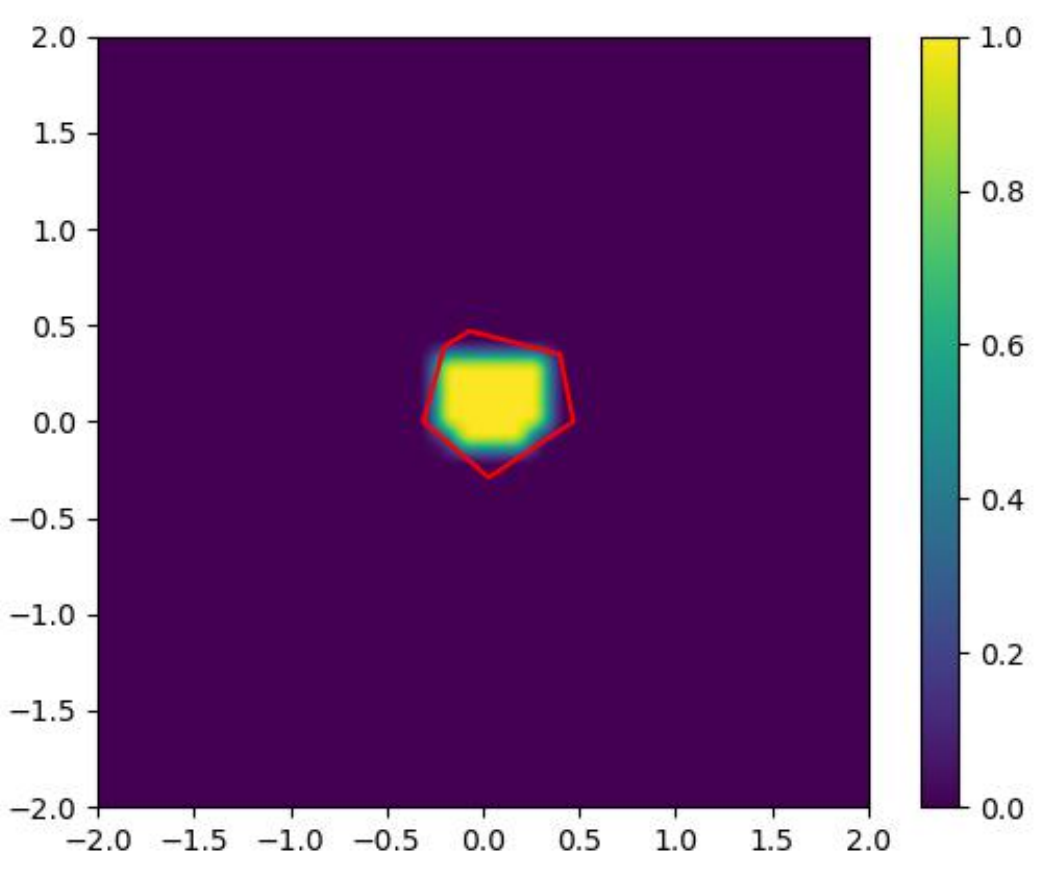}
	\end{minipage}
  }

   \subfigure[0.0357, $0.8249$]
  {
	\begin{minipage}{3cm}
 	\centering
	\includegraphics[scale=0.18]{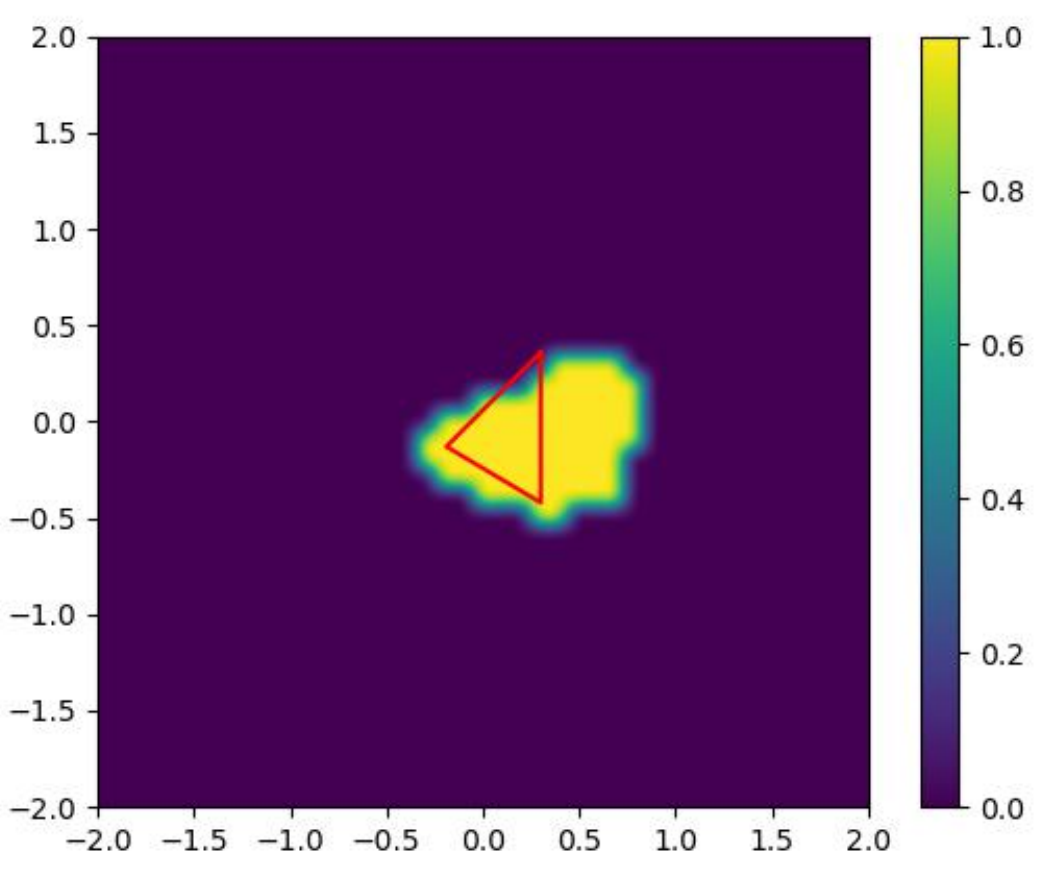}
	\end{minipage}
  }
  \subfigure[0.0619, $0.7383$]
  {
	\begin{minipage}{3cm}
 	
	\includegraphics[scale=0.18]{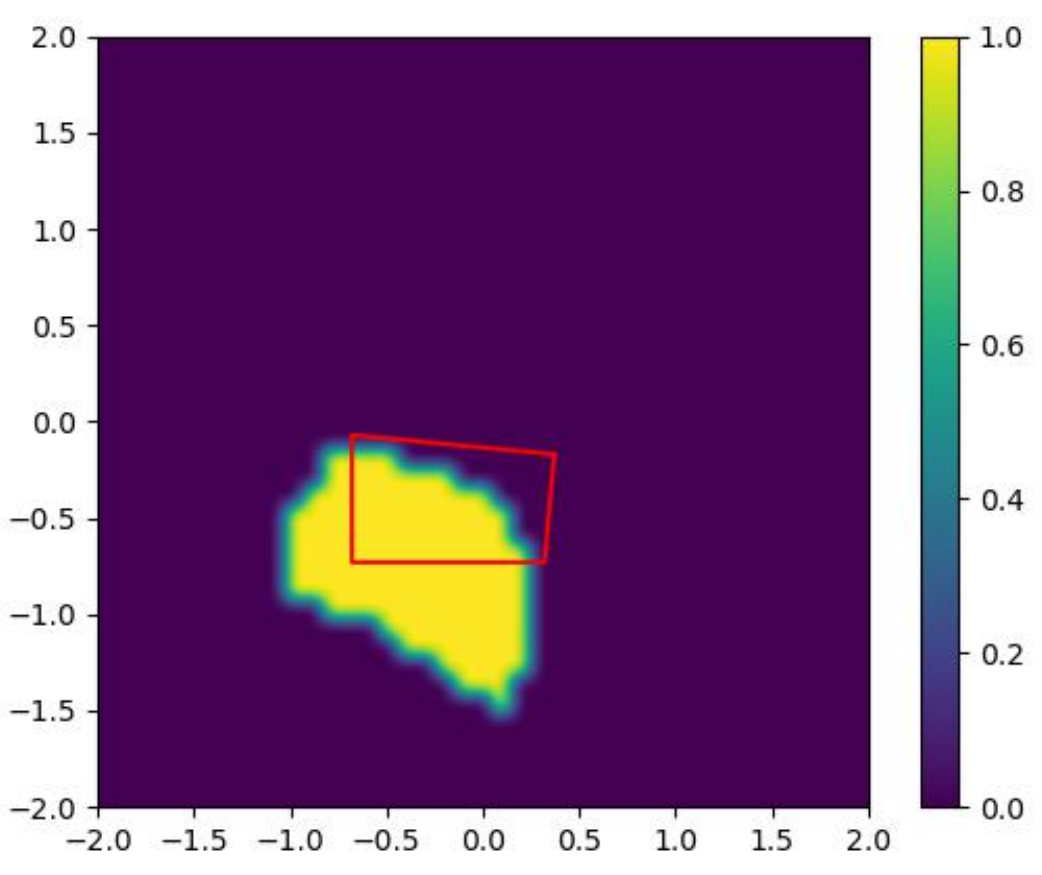}
	\end{minipage}
  }
   \subfigure[0.0553, $0.7855$]
  {
	\begin{minipage}{3cm}
 	\centering
	\includegraphics[scale=0.18]{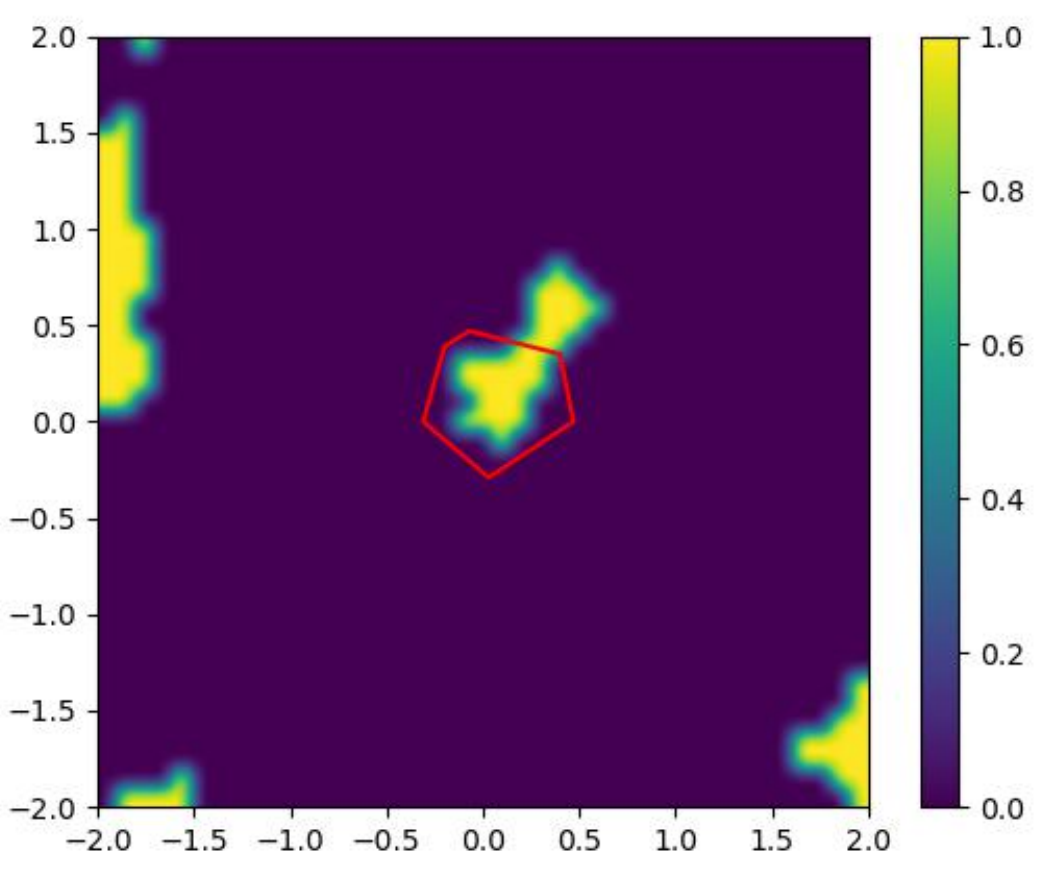}
	\end{minipage}
  }

  \subfigure[0.0080, $0.8860$]
  {
	\begin{minipage}{3cm}
 	\centering
	\includegraphics[scale=0.18]{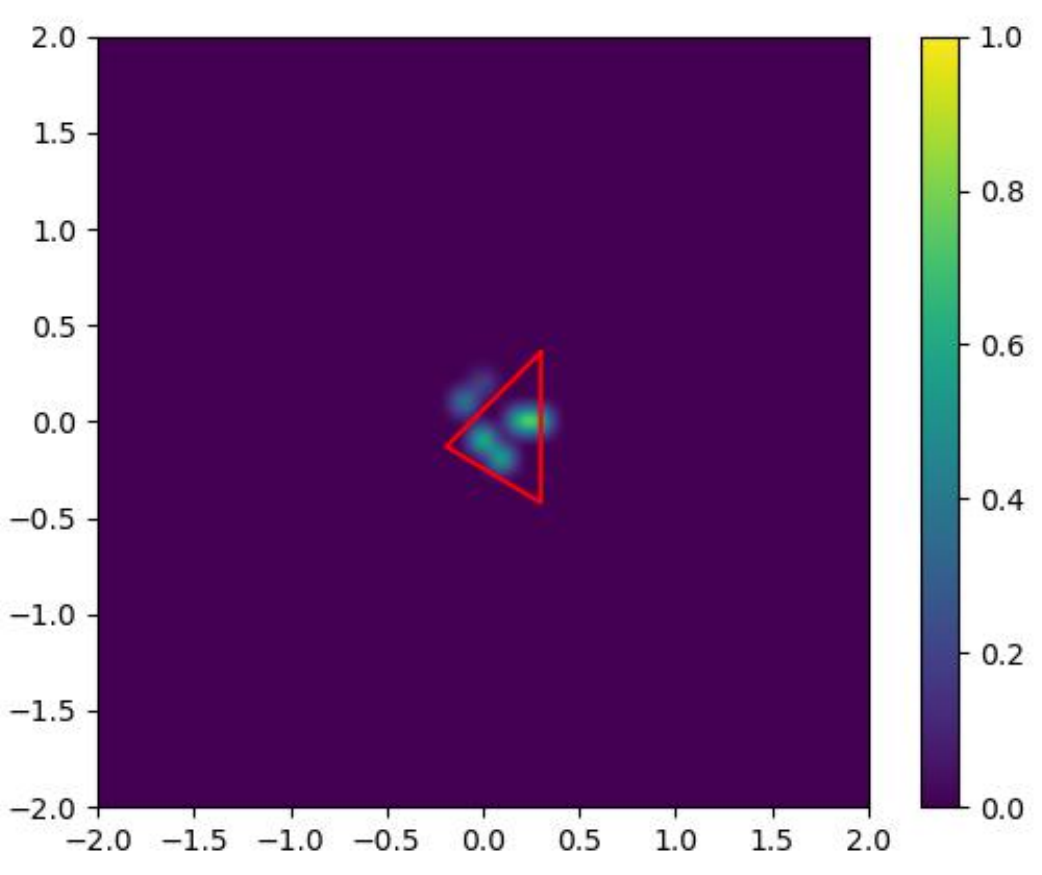}
	\end{minipage}
  }
  \subfigure[0.0388, $0.7909$]
  {
	\begin{minipage}{3cm}
 	
	\includegraphics[scale=0.18]{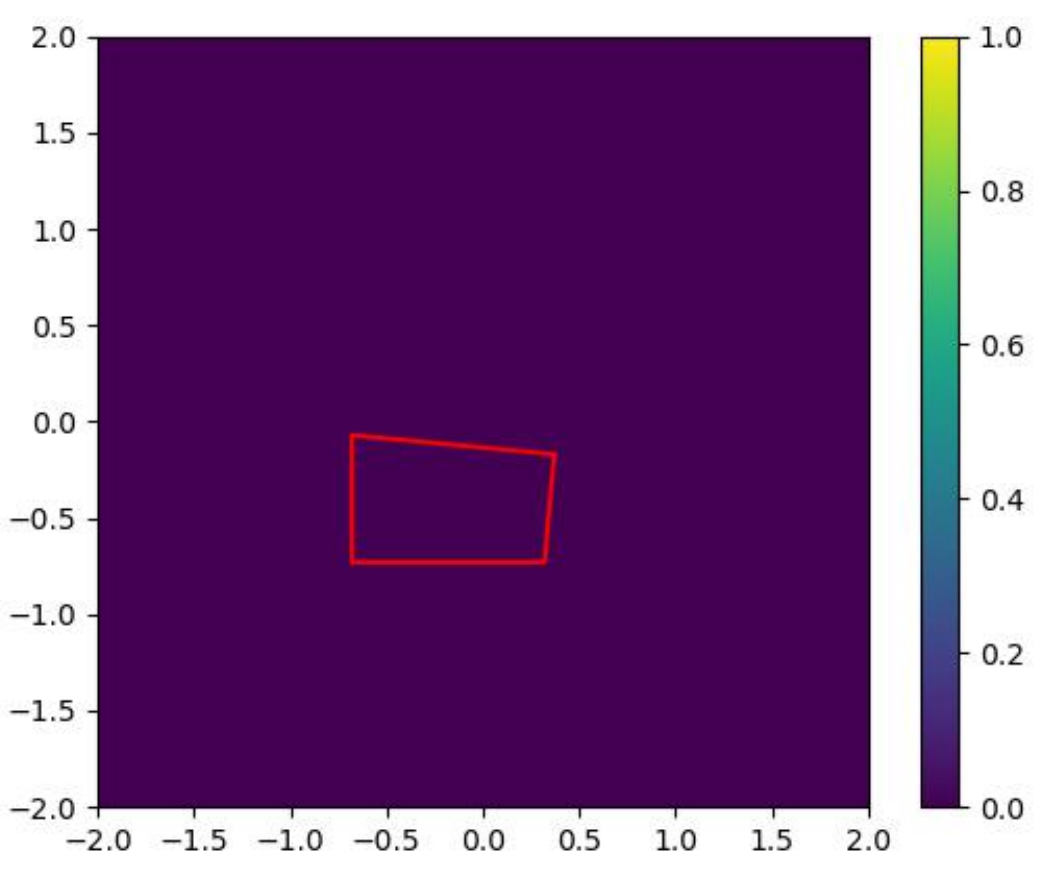}
	\end{minipage}
  }
   \subfigure[0.0185, $0.8516$]
  {
	\begin{minipage}{3cm}
 	\centering
	\includegraphics[scale=0.18]{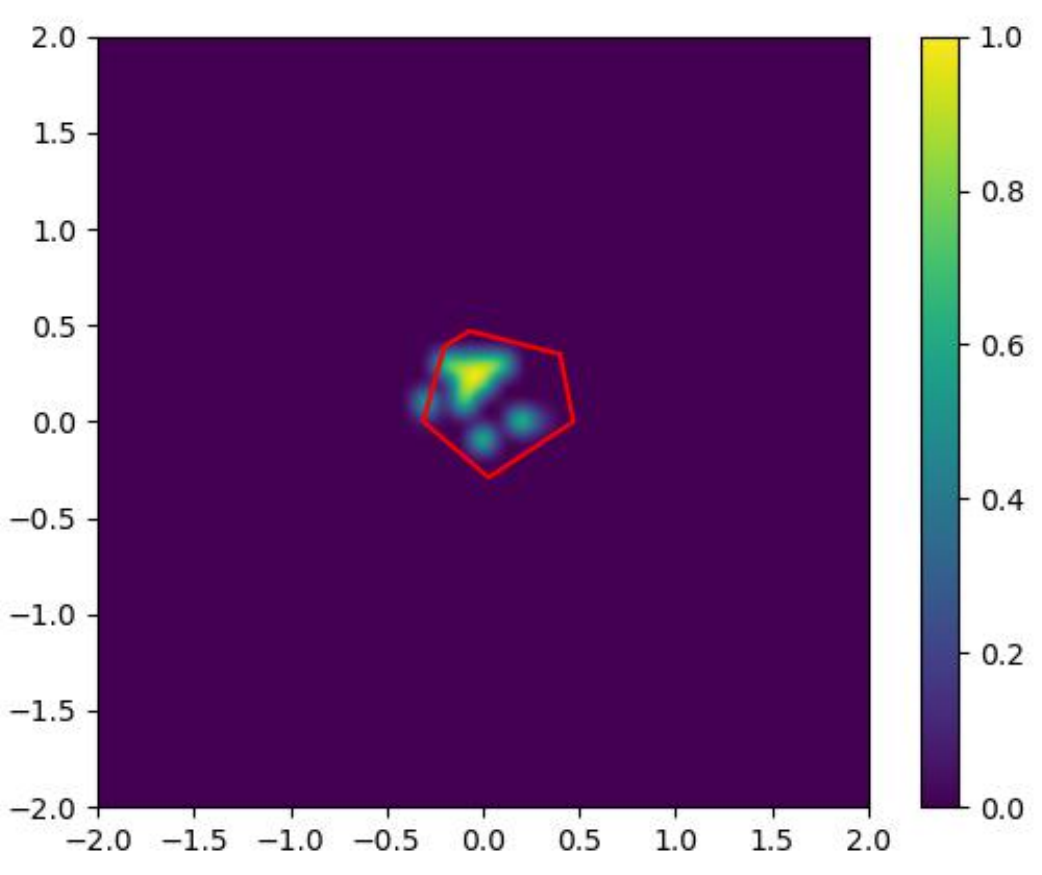}
	\end{minipage}
  }

  \subfigure[0.0101, $0.9032$]
  {
	\begin{minipage}{3cm}
 	\centering
	\includegraphics[scale=0.18]{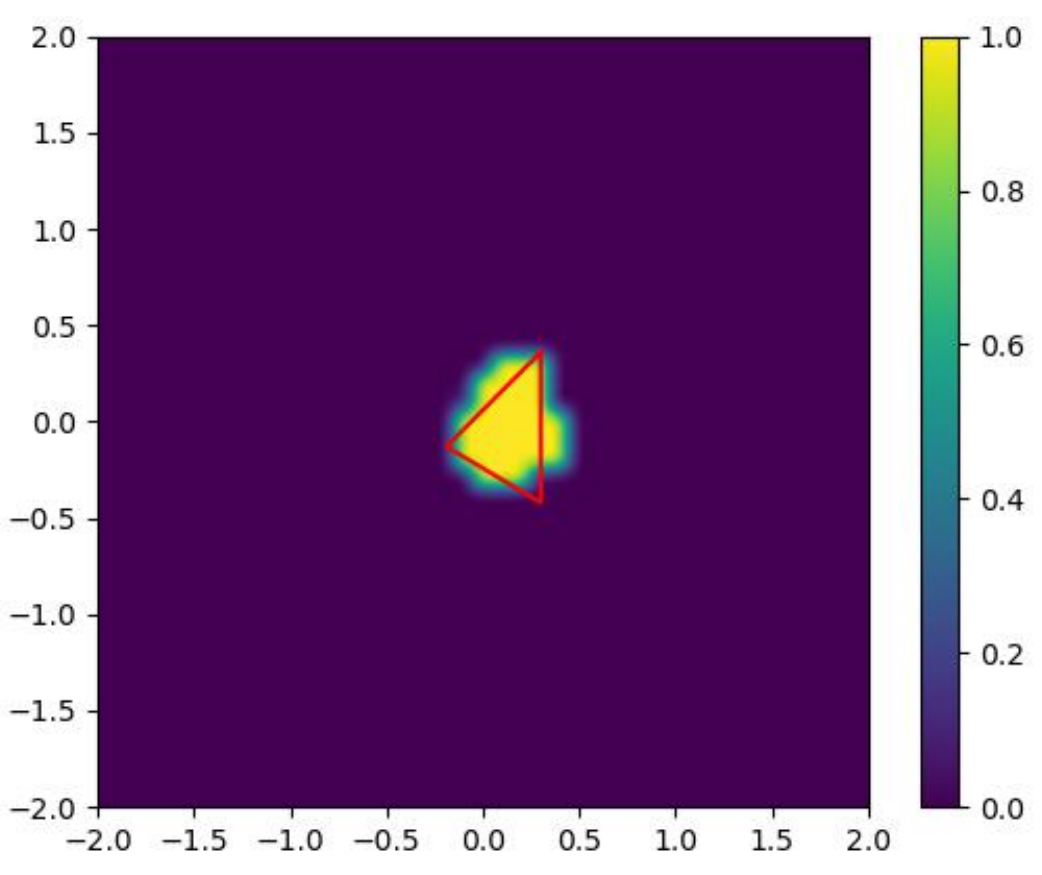}
	\end{minipage}
  }
  \subfigure[0.0170, $0.8568$]
  {
	\begin{minipage}{3cm}
 	
	\includegraphics[scale=0.18]{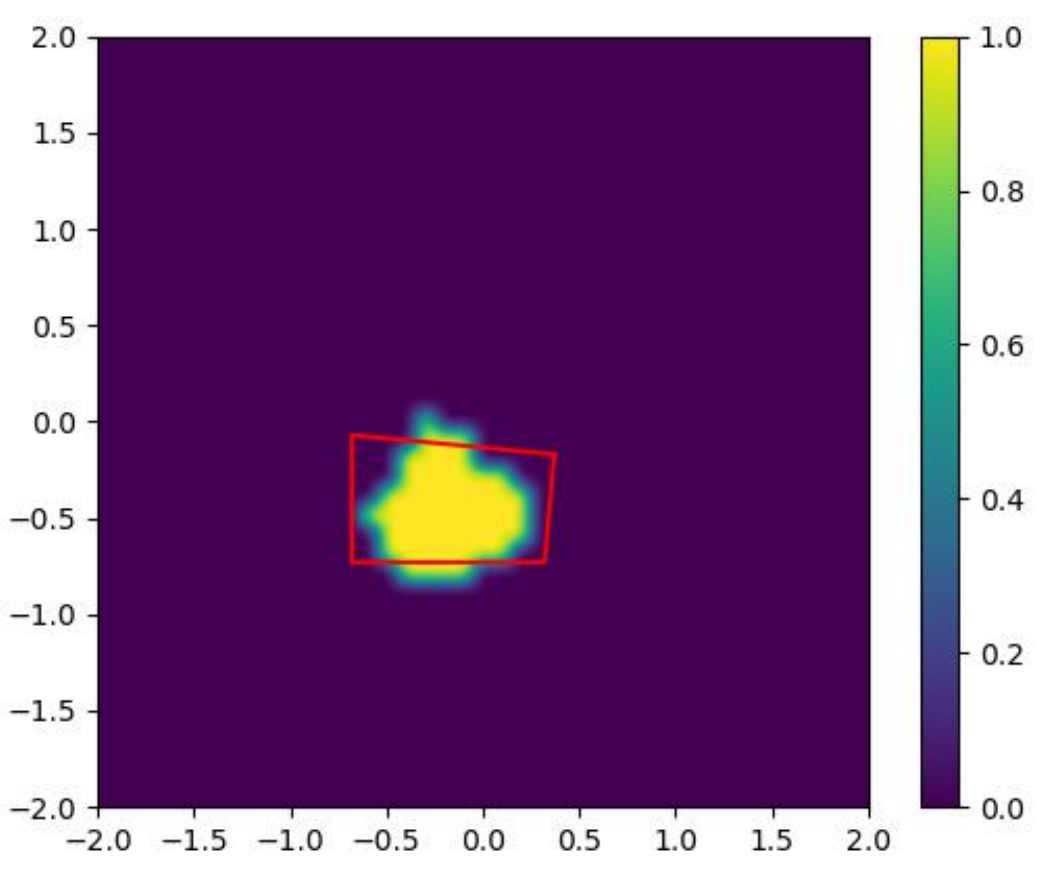}
	\end{minipage}
  }
   \subfigure[0.0059, $0.9410$]
  {
	\begin{minipage}{3cm}
 	\centering
	\includegraphics[scale=0.18]{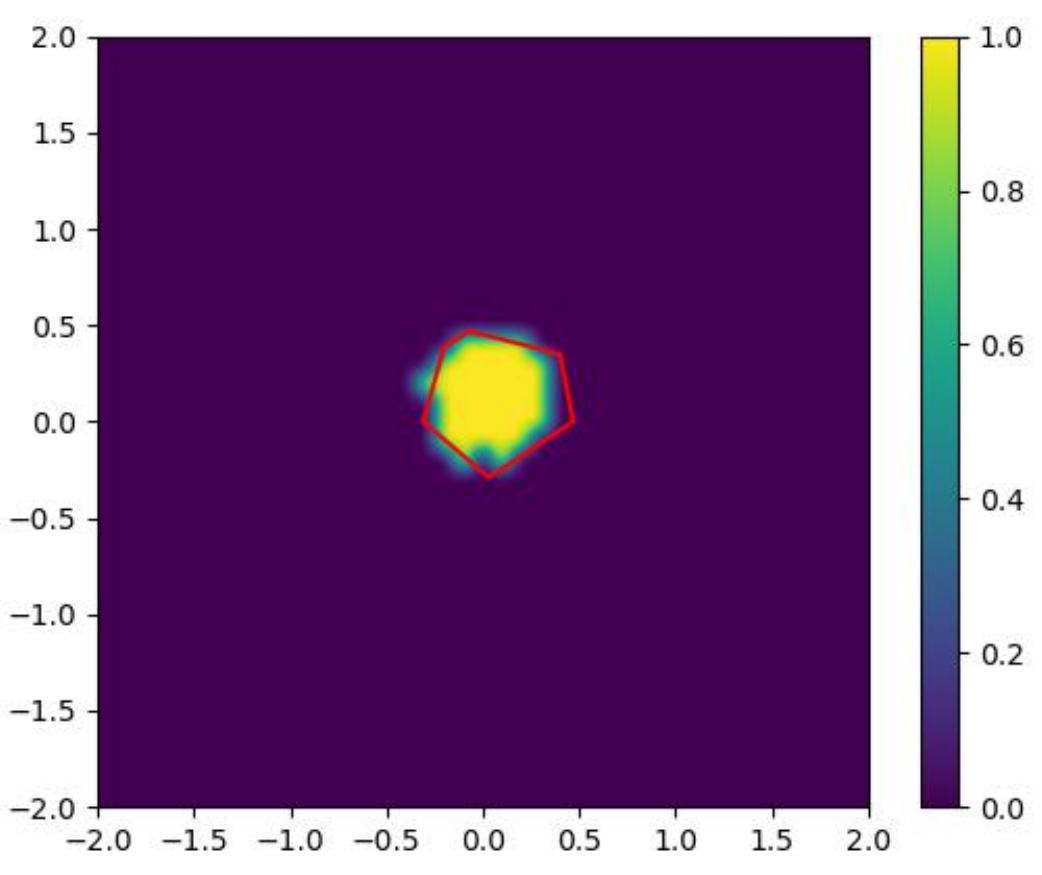}
	\end{minipage}
  }
 \caption{The comparisons between the LRT and other strategies in the noiseless case ($\delta=0$). First row: the LRT. Second row: the RT (deterministic, no learning). Third row: end-to-end fully connected NNs. Fourth row: the LRT with random initialization during training.}
 \label{Comparison}
\end{figure}

We first concentrate on the comparison with the deterministic strategy, i.e., the RT without learning described in $\S$\ref{Algorithm-RT} and $\S$\ref{sub:LRT}. Three samples in the test set are selected, and the reconstructions with the two strategies are shown in the first two rows of Figure~\ref{Comparison}. Moreover, we employ the deterministic strategy to reconstruct the inclusion in the test set with different parameters, or more precisely, the optimal parameters for each $\mathcal{T}_l$, $l=1,2,3$. The corresponding error distributions for the noise-free and noisy cases are recorded in Figures~\ref{Histogramc} and \ref{Histogramd}, respectively.

\begin{figure}%[htbp]
\centering
 \subfigure[0.0077, $0.9168$]
  {
	\begin{minipage}{3cm}
 	\centering
	\includegraphics[scale=0.18]{Figures/Noise/P3-1-Err.pdf}
	\end{minipage}
  }
  \subfigure[0.0160, $0.8485$]
  {
	\begin{minipage}{3cm}
 	
	\includegraphics[scale=0.18]{Figures/Noise/P4-2-Err.pdf}
	\end{minipage}
  }
   \subfigure[0.0089, $0.9021$]
  {
	\begin{minipage}{3cm}
 	\centering
	\includegraphics[scale=0.18]{Figures/Noise/P6-1-Err.pdf}
	\end{minipage}
  }

   \subfigure[0.0161, $0.8142$]
  {
	\begin{minipage}{3cm}
 	\centering
	\includegraphics[scale=0.18]{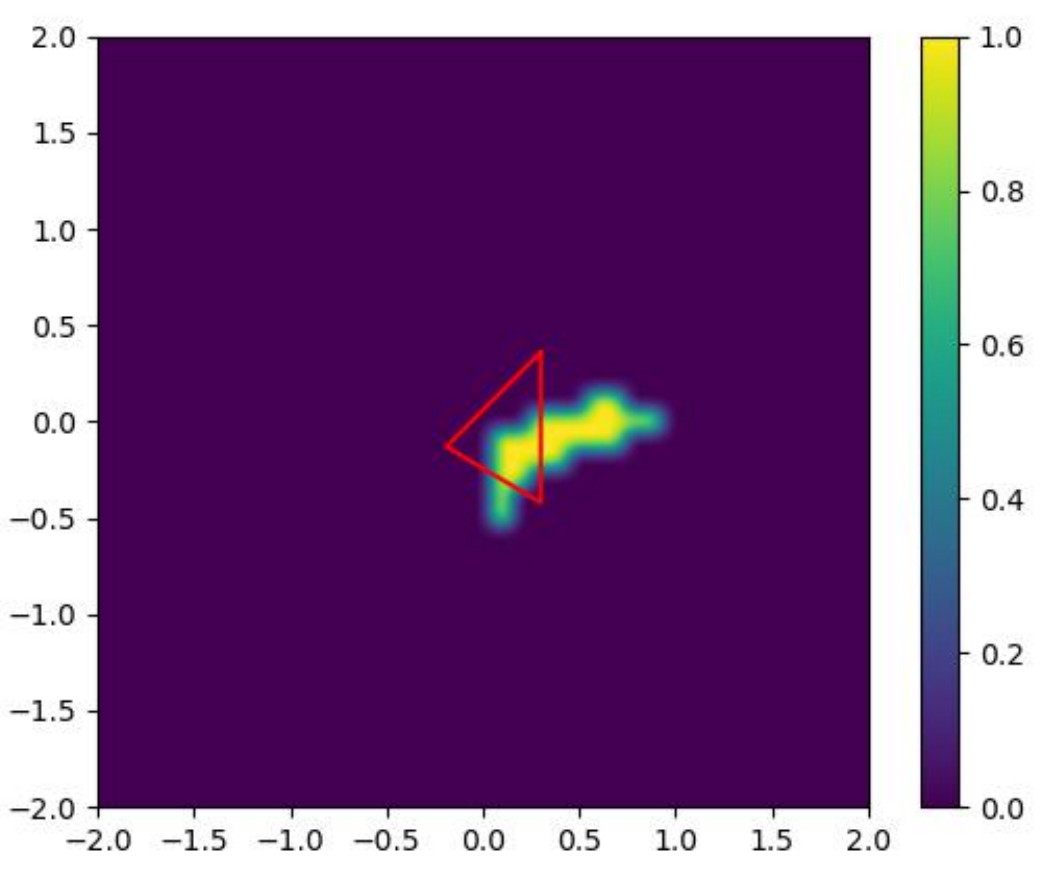}
	\end{minipage}
  }
  \subfigure[0.0333, $0.7647$]
  {
	\begin{minipage}{3cm}
 	
	\includegraphics[scale=0.18]{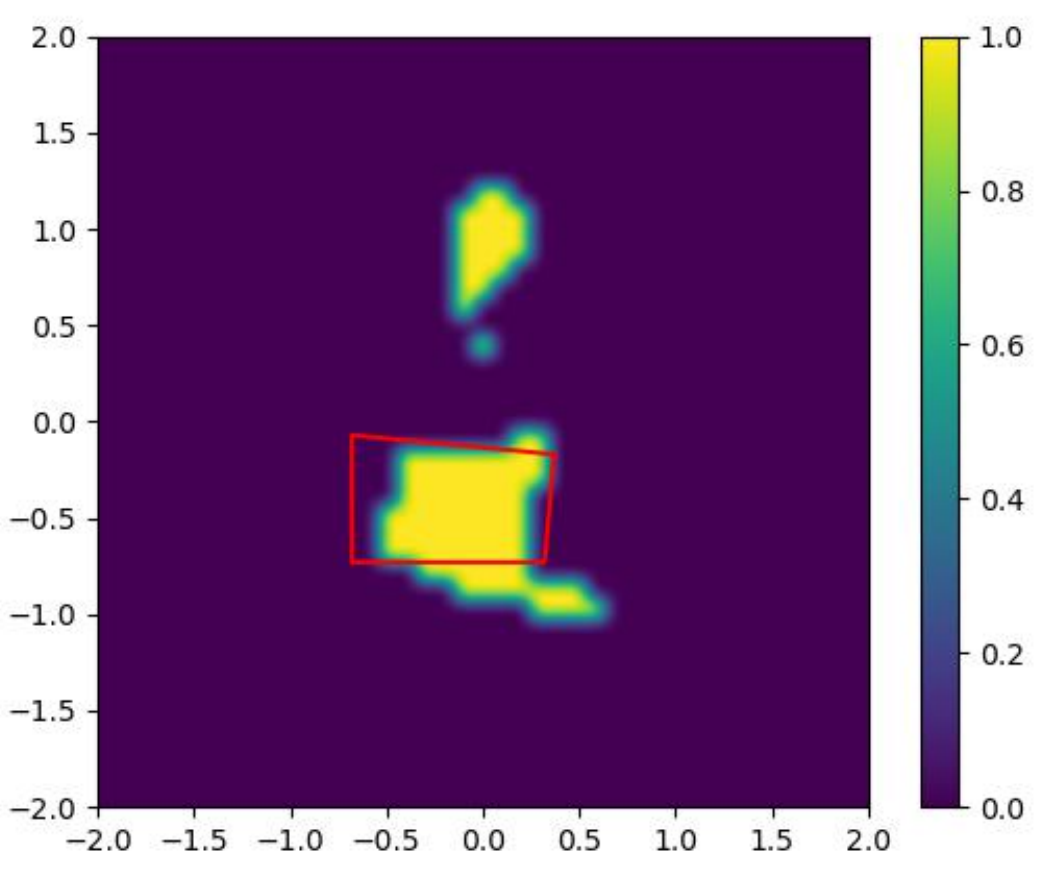}
	\end{minipage}
  }
   \subfigure[0.0321, $0.8262$]
  {
	\begin{minipage}{3cm}
 	\centering
	\includegraphics[scale=0.18]{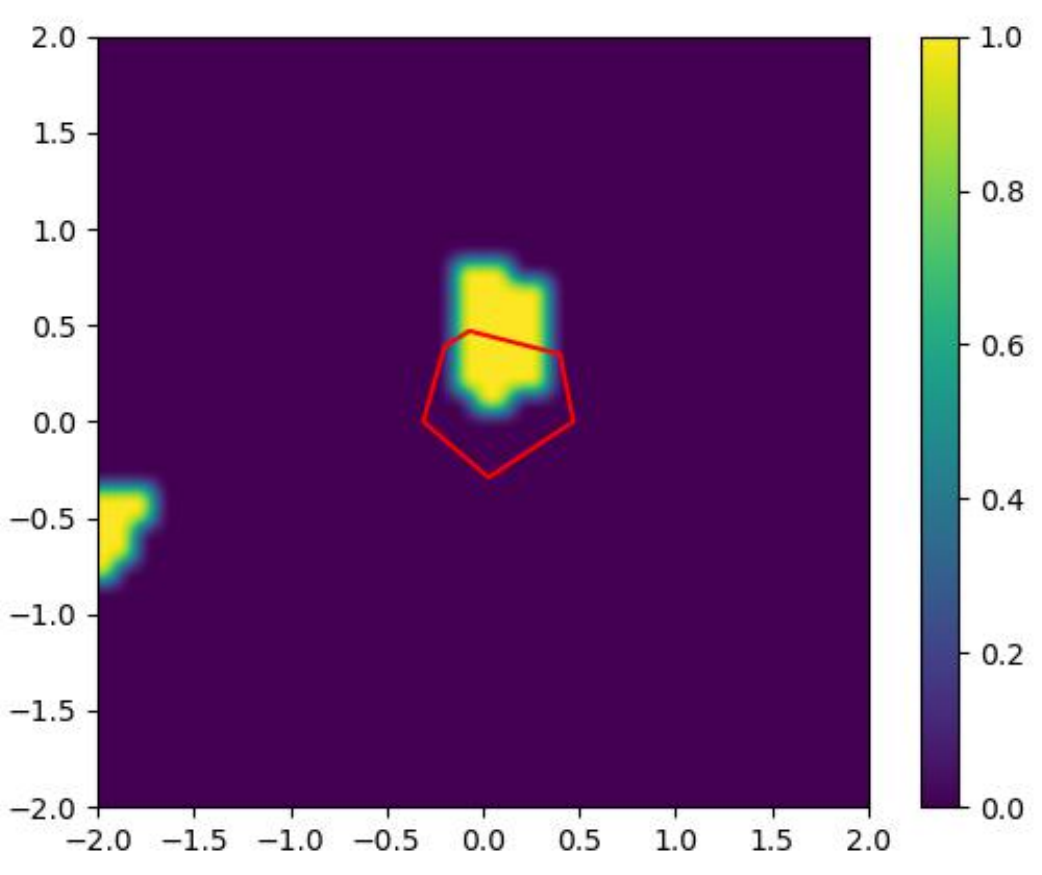}
	\end{minipage}
  }

  \subfigure[0.0087, $0.8880$]
  {
	\begin{minipage}{3cm}
 	\centering
	\includegraphics[scale=0.18]{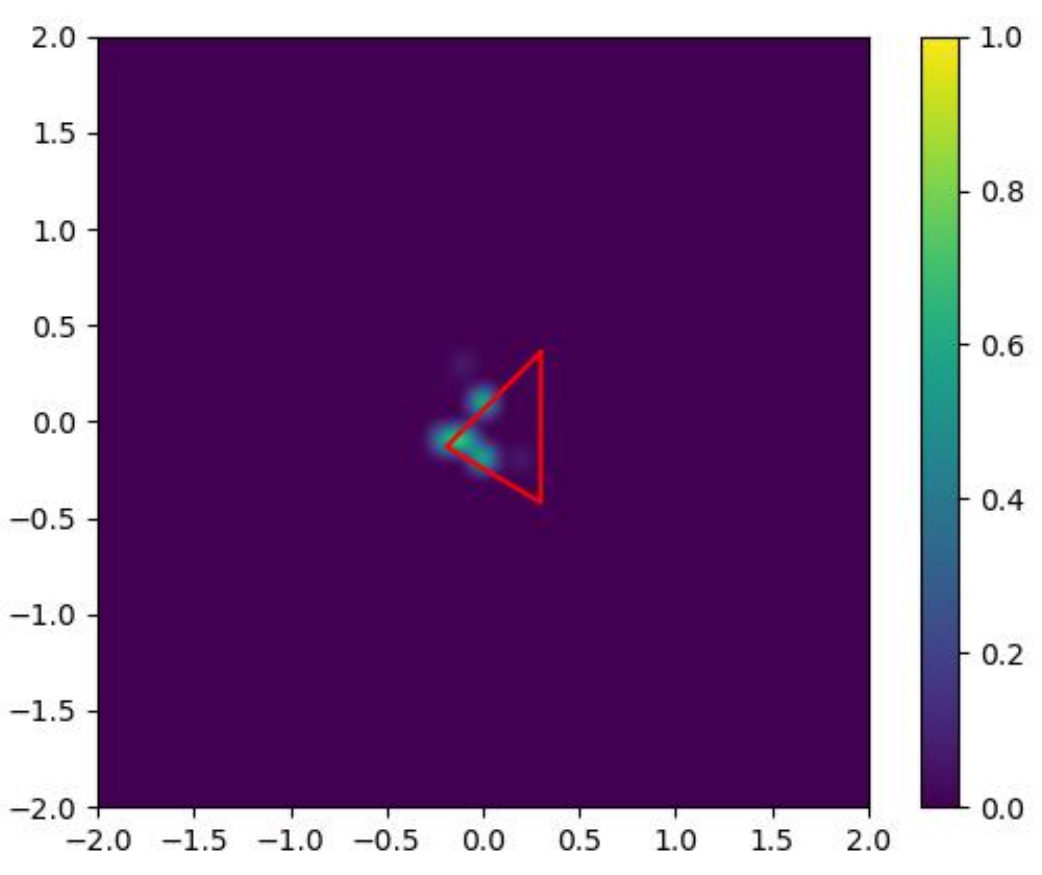}
	\end{minipage}
  }
  \subfigure[0.0375, $0.7909$]
  {
	\begin{minipage}{3cm}
 	
	\includegraphics[scale=0.18]{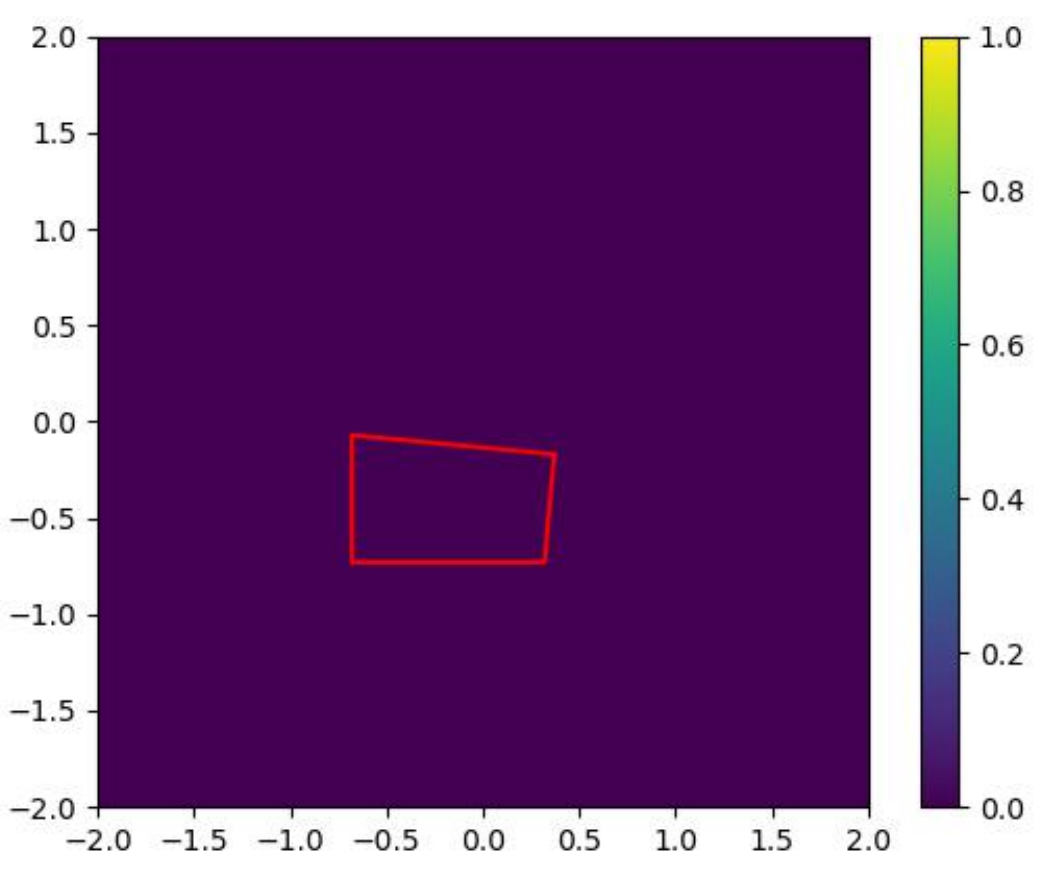}
	\end{minipage}
  }
   \subfigure[0.0191, $0.8465$]
  {
	\begin{minipage}{3cm}
 	\centering
	\includegraphics[scale=0.18]{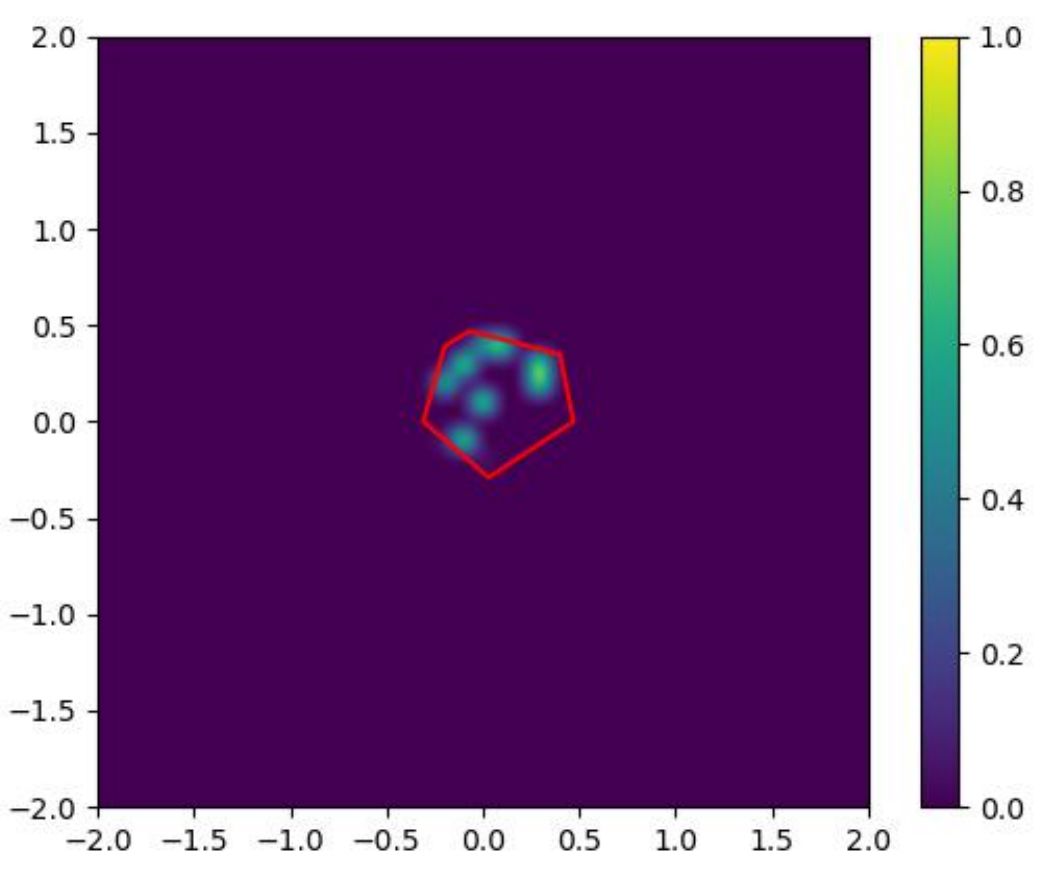}
	\end{minipage}
  }

  \subfigure[0.0131, $0.8774$]
  {
	\begin{minipage}{3cm}
 	\centering
	\includegraphics[scale=0.18]{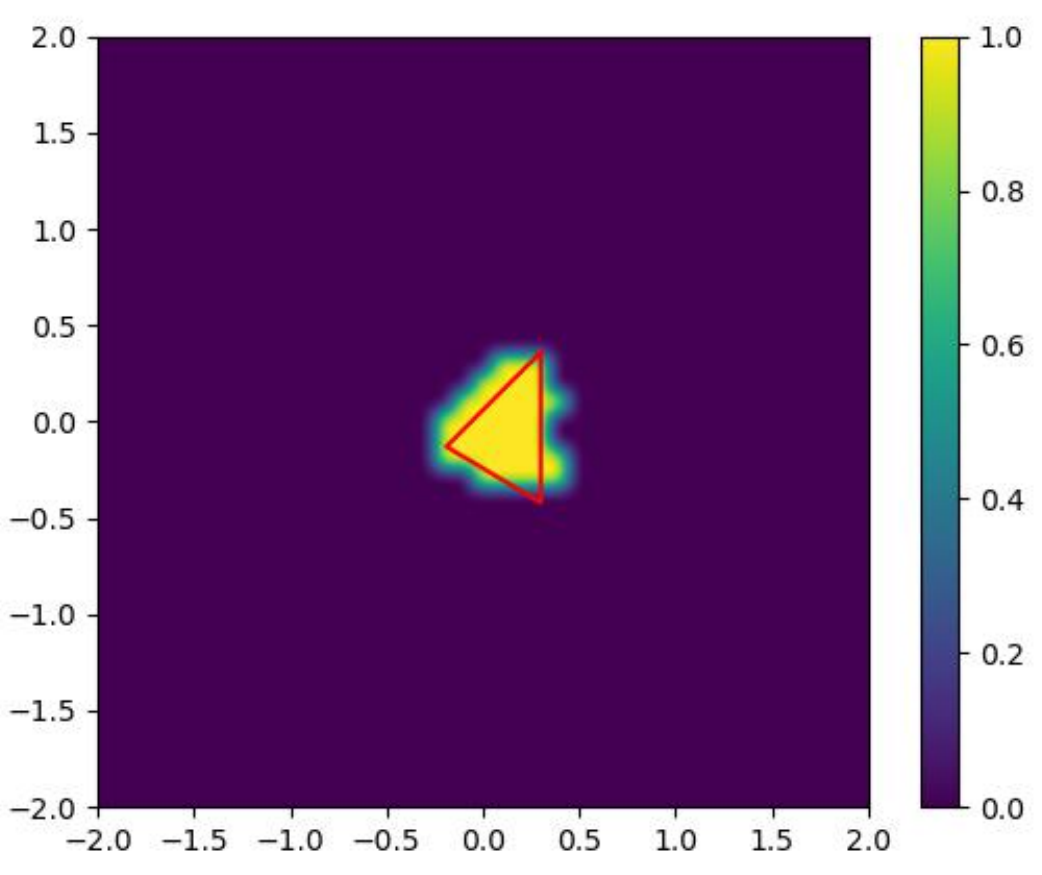}
	\end{minipage}
  }
  \subfigure[0.0190, $0.8465$]
  {
	\begin{minipage}{3cm}
 	
	\includegraphics[scale=0.18]{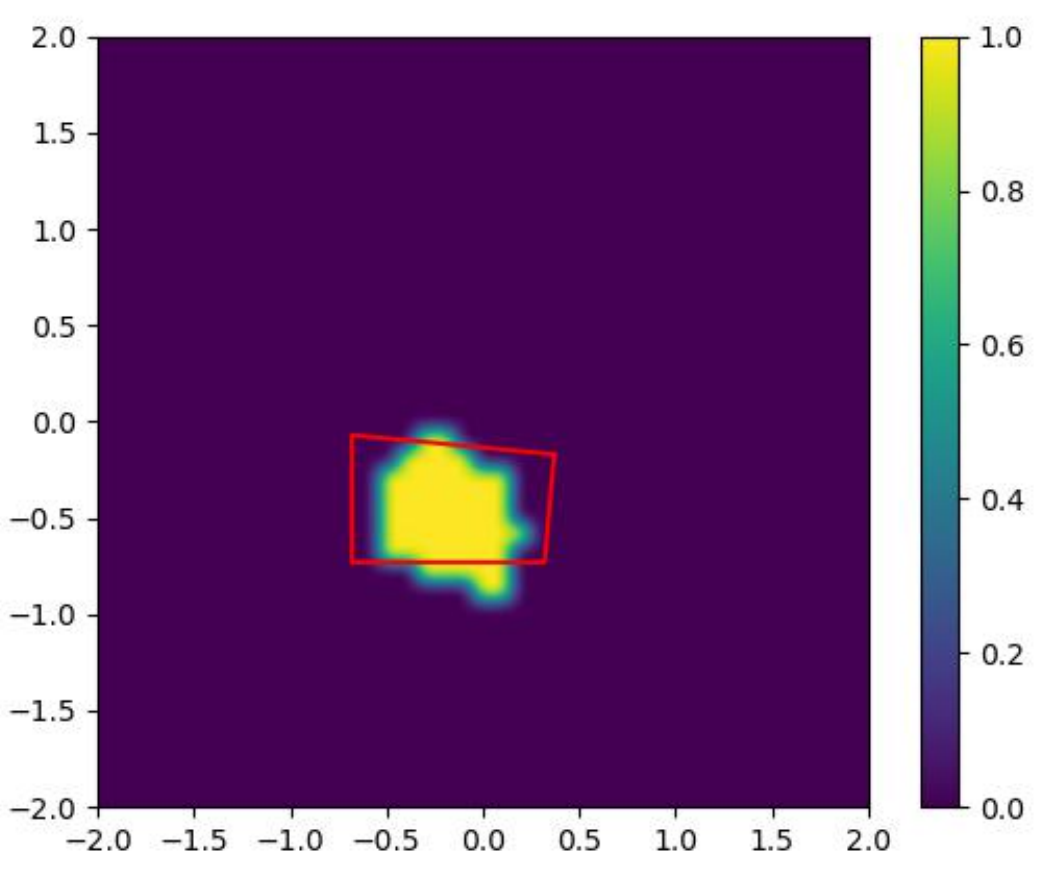}
	\end{minipage}
  }
   \subfigure[0.0270, $0.9351$]
  {
	\begin{minipage}{3cm}
 	\centering
	\includegraphics[scale=0.18]{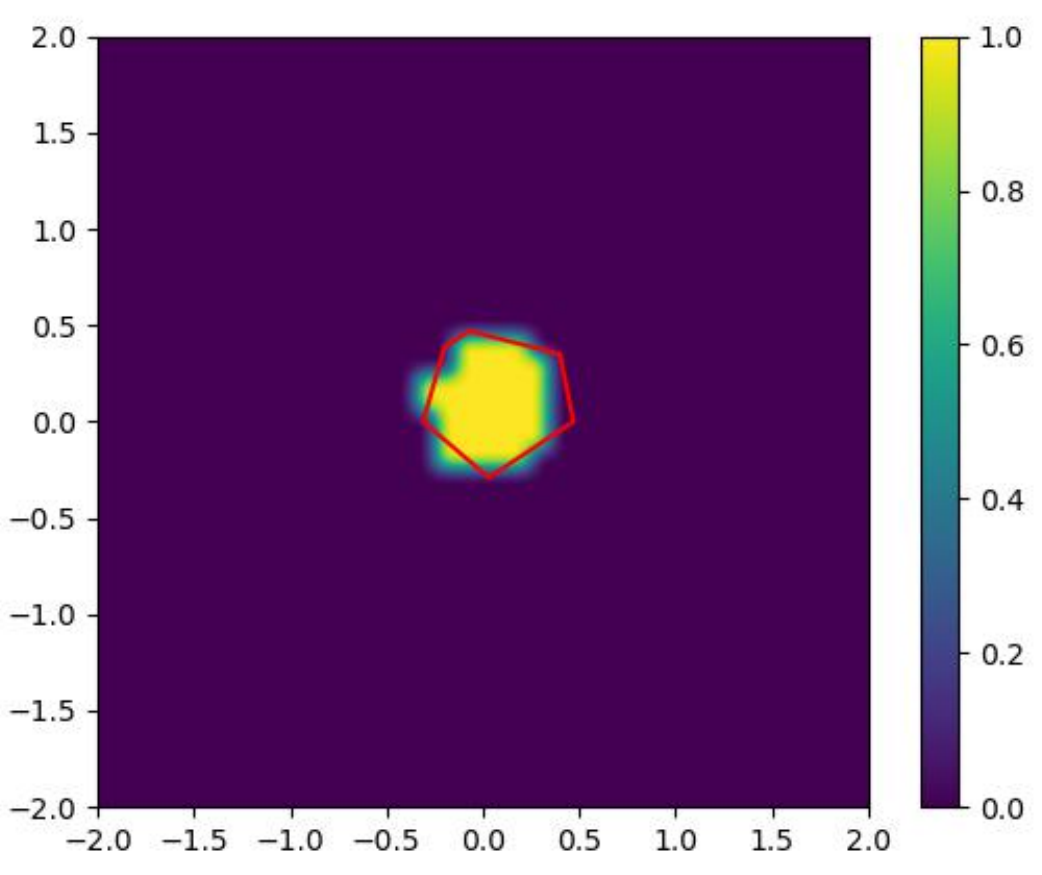}
	\end{minipage}
  }
 \caption{The comparisons between the LRT and other strategies in the noisy case ($\delta=3\%$). First row: the LRT. Second row: the RT (deterministic, no learning). Third row: end-to-end fully connected NNs. Fourth row: the LRT with random initialization during training.}
 \label{Comparison-Noise}
\end{figure}

We can conclude from Figures~\ref{Histogram}, \ref{Comparison} and \ref{Comparison-Noise} that the LRT substantially improves the reconstructions obtained by the RT, which illustrates the validity of the learning strategy. More precisely, compared with the RT, we highlight that the advantages of the LRT come from the following two aspects. First, the RT requires different parameters $\alpha$ and $\mathcal C$ for different inclusions to generate reasonable reconstructions, while this problem is addressed automatically by the LRT. Second, the computational time for the RT to finish a reconstruction is about $0.5$h, while the LRT can make it within $1$s after training, which greatly improves the computational efficiency.

\subsubsection{End-to-end fully connected neural network}
It is natural to compare the LRT with  a standard end-to-end fully connected neural network. For a more faithful comparison, we use the same three-step strategy of the LRT, namely, we design $L=3$ fully connected NNs with four layers, where the number of the neurons are $400$, $1200$, $1200$ and $1681$, respectively. The activation functions between the four layers are  the ReLUs. 
The NNs take as input the measurements $y^\delta(B)$ and output an approximation of $\chi(B)$, for $B\in\mathcal{T}_\ell$. We use the same reconstruction algorithm discussed in $\S$\ref{S-TSS}, in which the learned classifier is combined with the learned fully connected NNs (which replace the $\Phi_{\theta_\ell}$). 
 The numerical results of the same inclusions considered in the previous case ($\S$\ref{subsub:RT}) are shown in Figures~\ref{Comparison} and \ref{Comparison-Noise}, third row. We also test other fully connected NNs with different depths and widths,  as well as CNNs because of their relevance for imaging tasks, but they turn out to give similar results, which we omit for the sake of conciseness. As it is evident, the results are much worse than those of the LRT, and the NNs tend to favor the value $0$ in the reconstructions. This indicates that the architecture of the networks used in the LRT, which takes into account the physical model and the RT method, is a key aspect for a successful recovery.

\subsubsection{The LRT without the initialization given by the deterministc RT}

\begin{figure}%[htbp]
\centering
 \subfigure[Initialization with deterministic RT, $\delta=0$.]
  {
	\begin{minipage}{.47\textwidth}
 	\centering
	\includegraphics[width=\columnwidth]{Figures/Noise-free/Histogram-Noise-Free.pdf}
	\end{minipage}
  }
   \subfigure[Initialization with random weights, $\delta=0$.]
  {
	\begin{minipage}{.47\textwidth}
 	\centering
	\includegraphics[width=\columnwidth]{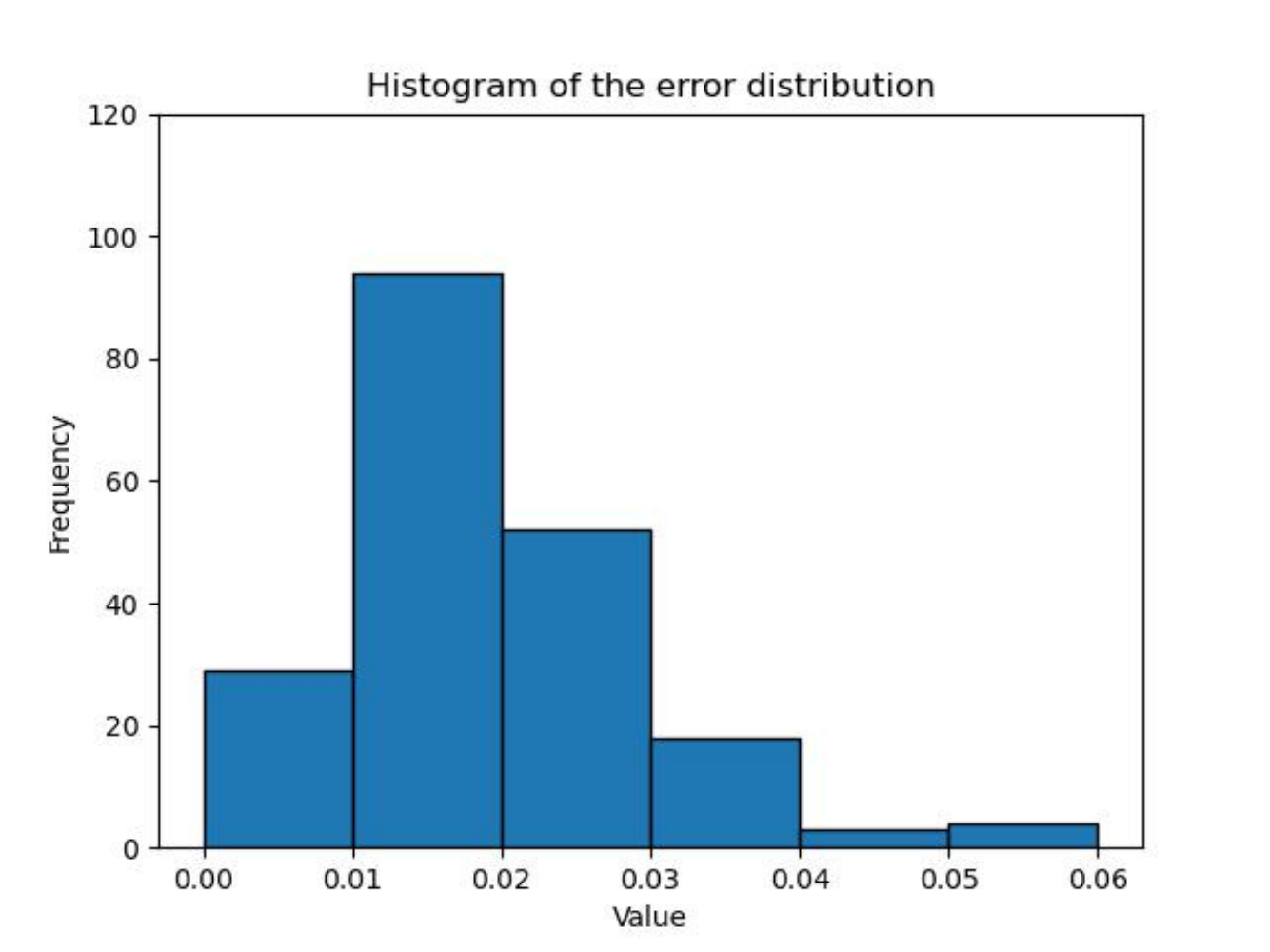}
	\end{minipage}
  }

  \subfigure[Initialization with deterministic RT, $\delta=3\%$.]
  {
	\begin{minipage}{.47\textwidth}
 	\centering
	\includegraphics[width=\columnwidth]{Figures/Noise/Histogram-Noise.pdf}
	\end{minipage}
  }
   \subfigure[Initialization with random weights,  $\delta=3\%$.]
  {
	\begin{minipage}{.47\textwidth}
 	\centering
	\includegraphics[width=\columnwidth]{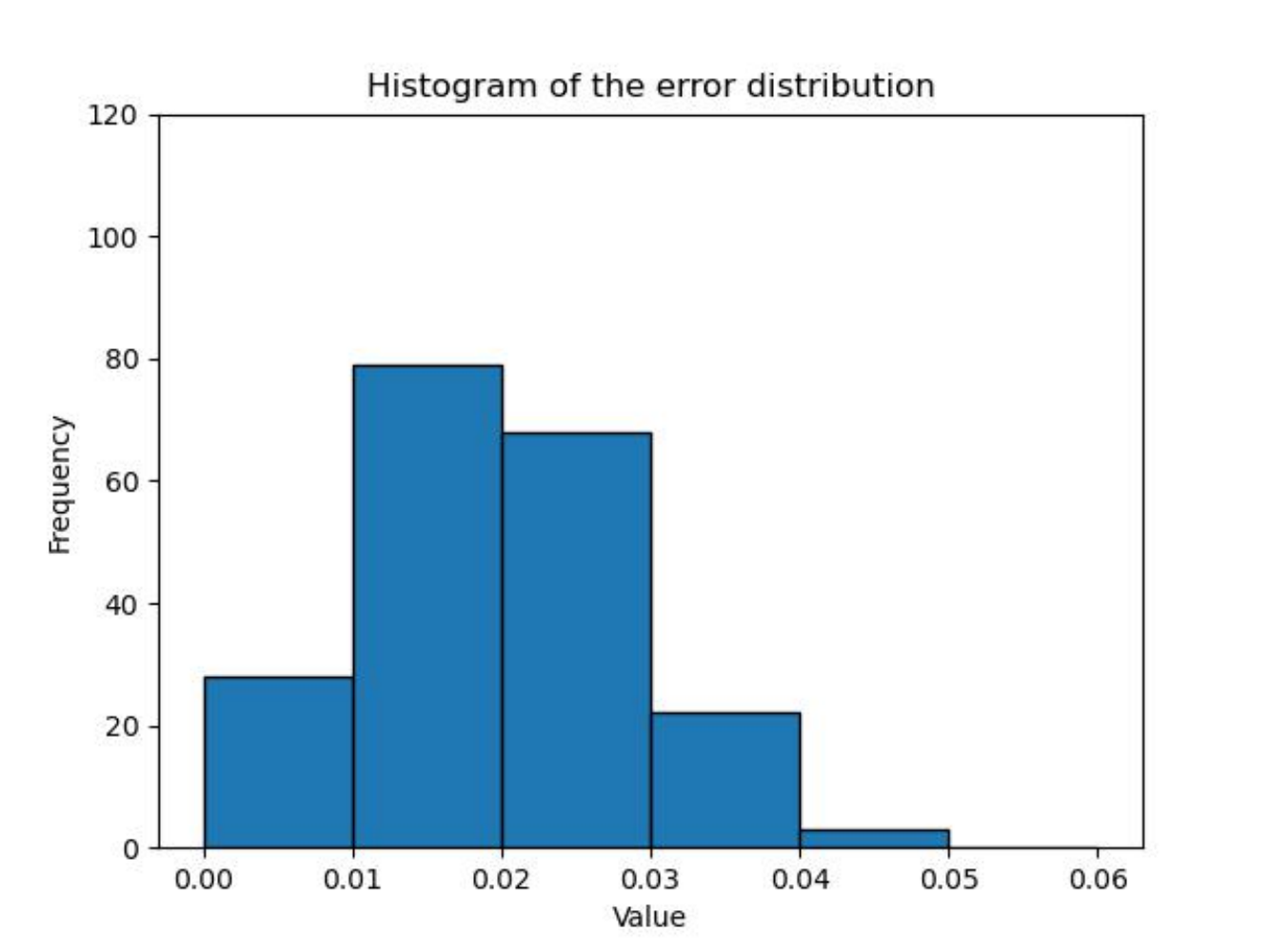}
	\end{minipage}
  }
 \caption{The error distribution of the mean squared errors for all the samples in the test set with the reconstruction obtained with the LRT with different initializations during training.}
 \label{Histogram-LRT-RTLNI}
\end{figure}

In the LRT, during training we initialize the weights of the neural networks $\Phi_{\theta_\ell}$ with the weights obtained with the deterministic RT of $\S$\ref{sub:LRT}, as explained in $\S$\ref{sub:trainingnet}.
 Next, we give some numerical results to verify the validity of our initialization strategy. Specifically, we train the networks $\Phi_{\theta_\ell}$ with random initialization, and then we use the LRT as discussed above. The same three samples of the test set are used for comparison, and the results can be seen in Figures~\ref{Comparison} and \ref{Comparison-Noise}, fourth row. In addition, we also record the distributions of the mean squared errors for reconstructing all the samples in the test set by using the LRT with these two possible initialization strategies during training, which are depicted in Figure~\ref{Histogram-LRT-RTLNI}. These numerical results show that carefully initializing the weights by using the physical model and the deterministic RT has a small positive effect.

\section{Conclusions}\label{Sec:Conclusion}
In this paper, we have considered an inverse boundary value problem for the Laplace equation, in which an insulated inclusion has to be reconstructed from one boundary measurement. We have shown that the reconstruction algorithm based on the RT can be written as a neural network. This makes it possible to learn the weights of this network with a fixed architecture. This is further combined with a pre-trained classifier.

The numerical results suggest that our method is superior to both the standard RT method (fully deterministic and without learning) and to an end-to-end deep neural network (fully data-driven, and not taking into account the physical model). This provides yet another indication that, very often, the most accurate reconstruction algorithms for inverse problems are obtained as a combination of a model-based approach and a data-driven part.

Let us discuss some of the limitations of the current work. The method discussed in this paper is adapted to one single measurement and one single inclusion to detect. It would be natural to extend this approach to more complicated disconnected inclusions, perhaps by using multiple boundary measurements. Similarly, our method is primarily designed for the reconstruction of convex polygonal inclusions, based on the properties of the RT. When applied to inclusions with more complex geometries, the current approach may not achieve high accuracy as it does now. However, it can serve as an effective initial guess for iteration-based methods to further refine the reconstruction. Furthermore,  the proposed strategy currently relies on training two separate neural networks, which increases the complexity of the implementation. A promising direction for improvement would be to integrate the two networks into a unified architecture, potentially through the design of a more sophisticated network structure.

It would be interesting to investigate whether the method discussed in this paper, based on a ``write as a NN and then learn'' approach, can be extended to other qualitative methods for inverse problems. Finally, more extensive numerical results, possibly in real-world scenarios, are needed to fully evaluate the capabilities of this technique. These investigations are left for future research.

\section*{Acknowledgements}
This work has been carried out at the Machine Learning Genoa (MaLGa) center, University of Genoa. S.S.\ gratefully acknowledges the support provided by the China Scholarship Council program (Project ID: 202306090232), the National Natural Science Foundation of China (Nos. 12471395, 12071072) and the Jiangsu Provincial Scientific Research Center of Applied Mathematics under Grant No. BK20233002. G.S.A.\ is supported by the Air Force Office of Scientific Research under award number FA8655-23-1-7083 and by the European Union (ERC, SAMPDE, 101041040). Views and opinions expressed are however those of the authors only and do not necessarily reflect those of the European Union or the European Research Council. Neither the European Union nor the granting authority can be held responsible for them. G.S.A.\ is a member of the ``Gruppo Nazionale per l’Analisi Matematica, la Probabilità e le loro Applicazioni'', of the ``Istituto Nazionale di Alta Matematica''. The research was supported in part by the MIUR Excellence Department Project awarded to Dipartimento di Matematica, Università di Genova, CUP D33C23001110001. Co-funded by European Union – Next Generation EU, Missione 4 Componente 1 CUP D53D23005770006 and CUP D53D23016180001.

\bibliographystyle{abbrv}
\bibliography{biblio}

\begin{thebibliography}{10}

\bibitem{Akduman-2002-18}
I.~Akduman and R.~Kress.
\newblock Electrostatic imaging via conformal mapping.
\newblock {\em Inverse Problems}, 18(6):1659--1672, 2002.

\bibitem{2023-alberti-arroyo-santacesaria}
G.~S. Alberti, A.~Arroyo, and M.~Santacesaria.
\newblock Inverse problems on low-dimensional manifolds.
\newblock {\em Nonlinearity}, 36(1):734--808, 2023.

\bibitem{alberti2024learning}
G.~S. Alberti, L.~Ratti, M.~Santacesaria, and S.~Sciutto.
\newblock Learning a {G}aussian mixture for sparsity regularization in inverse
  problems.
\newblock {\em IMA Journal of Numerical Analysis}, to appear.

\bibitem{Alessandrini-2001-176}
G.~Alessandrini and L.~Rondi.
\newblock Optimal stability for the inverse problem of multiple cavities.
\newblock {\em J. Differential Equations}, 176(2):356--386, 2001.

\bibitem{alessandrini-scapin-2017}
G.~Alessandrini and A.~Scapin.
\newblock Depth dependent resolution in electrical impedance tomography.
\newblock {\em J. Inverse Ill-Posed Probl.}, 25(3):391--402, 2017.

\bibitem{Alves-2009-25}
C.~Alves, R.~Kress, and P.~Serranho.
\newblock Iterative and range test methods for an inverse source problem for
  acoustic waves.
\newblock {\em Inverse Problems}, 25(5):055005, 17, 2009.

\bibitem{Ammari-2008-108}
H.~Ammari, H.~Kang, E.~Kim, K.~Louati, and M.~S. Vogelius.
\newblock A {MUSIC}-type algorithm for detecting internal corrosion from
  electrostatic boundary measurements.
\newblock {\em Numer. Math.}, 108(4):501--528, 2008.

\bibitem{Aparicio-1996-12}
N.~D. Aparicio and M.~K. Pidcock.
\newblock The boundary inverse problem for the {L}aplace equation in two
  dimensions.
\newblock {\em Inverse Problems}, 12(5):565--577, 1996.

\bibitem{2019-arridge-etal-actanumerica}
S.~Arridge, P.~Maass, O.~\"Oktem, and C.-B. Sch\"onlieb.
\newblock Solving inverse problems using data-driven models.
\newblock {\em Acta Numer.}, 28:1--174, 2019.

\bibitem{Lorenzo-2024}
L.~Audibert and S.~Meng.
\newblock Shape and parameter identification by the linear sampling method for
  a restricted {F}ourier integral operator.
\newblock {\em Inverse Problems}, 40(9):Paper No. 095007, 30, 2024.

\bibitem{2022-beretta-francini}
E.~Beretta and E.~Francini.
\newblock Global {L}ipschitz stability estimates for polygonal conductivity
  inclusions from boundary measurements.
\newblock {\em Appl. Anal.}, 101(10):3536--3549, 2022.

\bibitem{beretta-etal-2018}
E.~Beretta, S.~Micheletti, S.~Perotto, and M.~Santacesaria.
\newblock Reconstruction of a piecewise constant conductivity on a polygonal
  partition via shape optimization in {EIT}.
\newblock {\em J. Comput. Phys.}, 353:264--280, 2018.

\bibitem{Beretta-1999-30}
E.~Beretta and S.~Vessella.
\newblock Stable determination of boundaries from {C}auchy data.
\newblock {\em SIAM J. Math. Anal.}, 30(1):220--232, 1999.

\bibitem{Bonnet-2008-24}
M.~Bonnet.
\newblock Inverse acoustic scattering by small-obstacle expansion of a misfit
  function.
\newblock {\em Inverse Problems}, 24(3):035022, 27, 2008.

\bibitem{borcea-2002}
L.~Borcea.
\newblock Electrical impedance tomography.
\newblock {\em Inverse Problems}, 18(6):R99--R136, 2002.

\bibitem{Bourgeois-2010-4}
L.~Bourgeois and J.~Dard\'e.
\newblock A quasi-reversibility approach to solve the inverse obstacle problem.
\newblock {\em Inverse Probl. Imaging}, 4(3):351--377, 2010.

\bibitem{2023bubba}
T.~A. Bubba, M.~Burger, T.~Helin, and L.~Ratti.
\newblock Convex regularization in statistical inverse learning problems.
\newblock {\em Inverse Problems and Imaging}, 17(6):1193--1225, 2023.

\bibitem{Bubba-2021-14}
T.~A. Bubba, M.~Galinier, M.~Lassas, M.~Prato, L.~Ratti, and S.~Siltanen.
\newblock Deep neural networks for inverse problems with pseudodifferential
  operators: an application to limited-angle tomography.
\newblock {\em SIAM J. Imaging Sci.}, 14(2):470--505, 2021.

\bibitem{Bukhgeim-1999-15}
A.~L. Bukhgeim, J.~Cheng, and M.~Yamamoto.
\newblock Stability for an inverse boundary problem of determining a part of a
  boundary.
\newblock {\em Inverse Problems}, 15(4):1021--1032, 1999.

\bibitem{Burger-2001-17}
M.~Burger.
\newblock A level set method for inverse problems.
\newblock {\em Inverse Problems}, 17(5):1327--1355, 2001.

\bibitem{Cakoni-2014}
F.~Cakoni and D.~Colton.
\newblock {\em A qualitative approach to inverse scattering theory}, volume 188
  of {\em Applied Mathematical Sciences}.
\newblock Springer, New York, 2014.

\bibitem{cakoni2011linear}
F.~Cakoni, D.~Colton, and P.~Monk.
\newblock {\em The linear sampling method in inverse electromagnetic
  scattering}.
\newblock SIAM, 2011.

\bibitem{Cakoni-2007-1}
F.~Cakoni and R.~Kress.
\newblock Integral equations for inverse problems in corrosion detection from
  partial {C}auchy data.
\newblock {\em Inverse Probl. Imaging}, 1(2):229--245, 2007.

\bibitem{Cen-2023-493}
S.~Cen, B.~Jin, K.~Shin, and Z.~Zhou.
\newblock Electrical impedance tomography with deep {C}alder\'on method.
\newblock {\em J. Comput. Phys.}, 493:Paper No. 112427, 14, 2023.

\bibitem{Chow-2014-30}
Y.~T. Chow, K.~Ito, and J.~Zou.
\newblock A direct sampling method for electrical impedance tomography.
\newblock {\em Inverse Problems}, 30(9):095003, 25, 2014.

\bibitem{Chung-2024}
J.~Chung and S.~Gazzola.
\newblock Computational methods for large-scale inverse problems: a survey on
  hybrid projection methods.
\newblock {\em SIAM Rev.}, 66(2):205--284, 2024.

\bibitem{dehoop-etal-2021}
M.~V. de~Hoop, M.~Lassas, and C.~A. Wong.
\newblock Deep learning architectures for nonlinear operator functions and
  nonlinear inverse problems.
\newblock {\em Math. Stat. Learn.}, 4(1-2):1--86, 2021.

\bibitem{dicristo-rondi-2003}
M.~Di~Cristo and L.~Rondi.
\newblock Examples of exponential instability for inverse inclusion and
  scattering problems.
\newblock {\em Inverse Problems}, 19(3):685--701, 2003.

\bibitem{Engl}
H.~W. Engl, M.~Hanke, and A.~Neubauer.
\newblock {\em Regularization of inverse problems}, volume 375.
\newblock Springer Science \& Business Media, 1996.

\bibitem{Evans-2010}
L.~C. Evans.
\newblock {\em Partial differential equations}, volume~19 of {\em Graduate
  Studies in Mathematics}.
\newblock American Mathematical Society, Providence, RI, second edition, 2010.

\bibitem{1989-friedman-isakov}
A.~Friedman and V.~Isakov.
\newblock On the uniqueness in the inverse conductivity problem with one
  measurement.
\newblock {\em Indiana Univ. Math. J.}, 38(3):563--579, 1989.

\bibitem{garde-hyvonen-2020}
H.~Garde and N.~Hyv\"onen.
\newblock Optimal depth-dependent distinguishability bounds for electrical
  impedance tomography in arbitrary dimension.
\newblock {\em SIAM J. Appl. Math.}, 80(1):20--43, 2020.

\bibitem{gregor-lecun-2010}
K.~Gregor and Y.~LeCun.
\newblock Learning fast approximations of sparse coding.
\newblock In {\em Proceedings of the 27th International Conference on
  International Conference on Machine Learning}, ICML'10, page 399–406,
  Madison, WI, USA, 2010. Omnipress.

\bibitem{Guo-2021-43}
R.~Guo and J.~Jiang.
\newblock Construct deep neural networks based on direct sampling methods for
  solving electrical impedance tomography.
\newblock {\em SIAM J. Sci. Comput.}, 43(3):B678--B711, 2021.

\bibitem{Haddar-2015-21}
H.~Haddar and R.~Kress.
\newblock Conformal mappings and inverse boundary value problems.
\newblock {\em Inverse Problems}, 21(3):935--953, 2005.

\bibitem{Hanke_2024}
M.~Hanke.
\newblock Lipschitz stability of an inverse conductivity problem with two
  cauchy data pairs.
\newblock {\em Inverse Problems}, 40(10):105015, 2024.

\bibitem{Hu-2020}
Y.~Hu, B.~Wang, and T.~Li.
\newblock Corrosion shape reconstruction of the mixed boundary in electrostatic
  imaging.
\newblock {\em Adv. Difference Equ.}, pages Paper No. 405, 19, 2020.

\bibitem{Ikehata-1998-23}
M.~Ikehata.
\newblock Reconstruction of the shape of the inclusion by boundary
  measurements.
\newblock {\em Comm. Partial Differential Equations}, 23(7-8):1459--1474, 1998.

\bibitem{Ikehata-1999-15}
M.~Ikehata.
\newblock Enclosing a polygonal cavity in a two-dimensional bounded domain from
  {C}auchy data.
\newblock {\em Inverse Problems}, 15(5):1231--1241, 1999.

\bibitem{Ikehata-2002-18}
M.~Ikehata and T.~Ohe.
\newblock A numerical method for finding the convex hull of polygonal cavities
  using the enclosure method.
\newblock {\em Inverse Problems}, 18(1):111--124, 2002.

\bibitem{Jakubik-2008-58}
P.~Jakubik and R.~Potthast.
\newblock Testing the integrity of some cavity---the {C}auchy problem and the
  range test.
\newblock {\em Appl. Numer. Math.}, 58(6):899--914, 2008.

\bibitem{Karageorghis-2009-17}
A.~Karageorghis and D.~Lesnic.
\newblock Detection of cavities using the method of fundamental solutions.
\newblock {\em Inverse Probl. Sci. Eng.}, 17(6):803--820, 2009.

\bibitem{Kim-2002-62}
S.~Kim, O.~Kwon, and J.~K. Seo.
\newblock Location search techniques for a grounded conductor.
\newblock {\em SIAM J. Appl. Math.}, 62(4):1383--1393, 2002.

\bibitem{Kirsch-2007}
A.~Kirsch and N.~Grinberg.
\newblock {\em The Factorization Method for Inverse Problems}.
\newblock Oxford University Press, 12 2007.

\bibitem{Kress-2014}
R.~Kress.
\newblock {\em Linear integral equations}, volume~82 of {\em Applied
  Mathematical Sciences}.
\newblock Springer, New York, third edition, 2014.

\bibitem{Kress-2005-21}
R.~Kress and W.~Rundell.
\newblock Nonlinear integral equations and the iterative solution for an
  inverse boundary value problem.
\newblock {\em Inverse Problems}, 21(4):1207--1223, 2005.

\bibitem{Le-2023-784}
T.~Le, D.-L. Nguyen, V.~Nguyen, and T.~Truong.
\newblock Sampling type method combined with deep learning for inverse
  scattering with one incident wave.
\newblock In {\em Advances in inverse problems for partial differential
  equations}, volume 784 of {\em Contemp. Math.}, pages 63--80. Amer. Math.
  Soc., 2023.

\bibitem{Li-2008-30}
J.~Li, H.~Liu, and J.~Zou.
\newblock Multilevel linear sampling method for inverse scattering problems.
\newblock {\em SIAM J. Sci. Comput.}, 30(3):1228--1250, 2008.

\bibitem{Li-2024-40}
K.~Li, B.~Zhang, and H.~Zhang.
\newblock Reconstruction of inhomogeneous media by an iteration algorithm with
  a learned projector.
\newblock {\em Inverse Problems}, 40(7):075008, 2024.

\bibitem{Lin-2021-15}
Y.-H. Lin, G.~Nakamura, R.~Potthast, and H.~Wang.
\newblock Duality between range and no-response tests and its application for
  inverse problems.
\newblock {\em Inverse Probl. Imaging}, 15(2):367--386, 2021.

\bibitem{Meng-2024-SIAM}
S.~Meng and B.~Zhang.
\newblock A kernel machine learning for inverse source and scattering problems.
\newblock {\em SIAM J. Numer. Anal.}, 62(3):1443--1464, 2024.

\bibitem{mukherjee-etal-2023}
S.~Mukherjee, A.~Hauptmann, O.~Öktem, M.~Pereyra, and C.-B. Schönlieb.
\newblock Learned reconstruction methods with convergence guarantees: A survey
  of concepts and applications.
\newblock {\em IEEE Signal Processing Magazine}, 40(1):164--182, 2023.

\bibitem{Murtagh-1991-2}
F.~Murtagh.
\newblock Multilayer perceptrons for classification and regression.
\newblock {\em Neurocomputing}, 2(5-6):183--197, 1991.

\bibitem{Nakamura-2015}
G.~Nakamura and R.~Potthast.
\newblock {\em Inverse modeling}.
\newblock IOP Expanding Physics. IOP Publishing, Bristol, 2015.
\newblock An introduction to the theory and methods of inverse problems and
  data assimilation.

\bibitem{Ning-2024-40}
J.~Ning, F.~Han, and J.~Zou.
\newblock A direct sampling-based deep learning approach for inverse medium
  scattering problems.
\newblock {\em Inverse Problems}, 40(1):015005, 2024.

\bibitem{Ongie-2020}
G.~Ongie, A.~Jalal, C.~A. Metzler, R.~G. Baraniuk, A.~G. Dimakis, and
  R.~Willett.
\newblock Deep learning techniques for inverse problems in imaging.
\newblock {\em IEEE Journal on Selected Areas in Information Theory},
  1(1):39--56, 2020.

\bibitem{Potthast-2005-75}
R.~Potthast.
\newblock Sampling and probe methods---an algorithmical view.
\newblock {\em Computing}, 75(2-3):215--235, 2005.

\bibitem{Potthast-2006-22}
R.~Potthast.
\newblock A survey on sampling and probe methods for inverse problems.
\newblock {\em Inverse Problems}, 22(2):R1--R47, 2006.

\bibitem{Potthast-2003-19}
R.~Potthast, J.~Sylvester, and S.~Kusiak.
\newblock A `range test' for determining scatterers with unknown physical
  properties.
\newblock {\em Inverse Problems}, 19(3):533--547, 2003.

\bibitem{Raissi-2019-378}
M.~Raissi, P.~Perdikaris, and G.~E. Karniadakis.
\newblock Physics-informed neural networks: a deep learning framework for
  solving forward and inverse problems involving nonlinear partial differential
  equations.
\newblock {\em J. Comput. Phys.}, 378:686--707, 2019.

\bibitem{rondi-1999}
L.~Rondi.
\newblock Optimal stability estimates for the determination of defects by
  electrostatic measurements.
\newblock {\em Inverse Problems}, 15(5):1193--1212, 1999.

\bibitem{1996-seo}
J.~K. Seo.
\newblock On the uniqueness in the inverse conductivity problem.
\newblock {\em J. Fourier Anal. Appl.}, 2(3):227--235, 1996.

\bibitem{Sun-2023-485}
S.~Sun, G.~Nakamura, and H.~Wang.
\newblock Numerical studies of domain sampling methods for inverse boundary
  value problems by one measurement.
\newblock {\em J. Comput. Phys.}, 485:112099, 2023.

\bibitem{sun-2024}
S.~Sun, G.~Nakamura, and H.~Wang.
\newblock Domain sampling methods for an inverse boundary value problem of the
  heat equation.
\newblock {\em Inverse Problems}, 40(12):Paper No. 125009, 21, 2024.

\bibitem{Tan-2019}
H.~H. Tan and K.~H. Lim.
\newblock Vanishing gradient mitigation with deep learning neural network
  optimization.
\newblock In {\em 2019 7th International Conference on Smart Computing \&
  Communications (ICSCC)}, pages 1--4, 2019.

\bibitem{Tanyu-2023-39}
D.~N. Tanyu, J.~Ning, T.~Freudenberg, N.~Heilenk\"otter, A.~Rademacher,
  U.~Iben, and P.~Maass.
\newblock Deep learning methods for partial differential equations and related
  parameter identification problems.
\newblock {\em Inverse Problems}, 39(10):103001, 2023.

\bibitem{uhlmann-2009}
G.~Uhlmann.
\newblock Electrical impedance tomography and {C}alder\'on's problem.
\newblock {\em Inverse Problems}, 25(12):123011, 39, 2009.

\bibitem{Wang-2004-13}
Z.~Wang, A.~Bovik, H.~Sheikh, and E.~Simoncelli.
\newblock Image quality assessment: from error visibility to structural
  similarity.
\newblock {\em IEEE Transactions on Image Processing}, 13(4):600--612, 2004.

\bibitem{Wu-2021-7}
Y.~Wu, B.~Chen, K.~Liu, C.~Zhu, H.~Pan, J.~Jia, H.~Wu, and J.~Yao.
\newblock Shape reconstruction with multiphase conductivity for electrical
  impedance tomography using improved convolutional neural network method.
\newblock {\em IEEE Sensors Journal}, 21(7):9277--9287, 2021.

\bibitem{Ye-2018}
J.~C. Ye, Y.~Han, and E.~Cha.
\newblock Deep convolutional framelets: a general deep learning framework for
  inverse problems.
\newblock {\em SIAM J. Imaging Sci.}, 11(2):991--1048, 2018.

\bibitem{Zia-2016-304}
Q.~Zia and R.~Potthast.
\newblock The range test and the no response test for {O}seen problems:
  theoretical foundation.
\newblock {\em J. Comput. Appl. Math.}, 304:201--211, 2016.

\end{thebibliography}

\appendix

\section{Further details on the computational aspects} \label{app:Training}

In this section, we provide additional details on the training algorithms employed for the neural networks used in the LRT, and on the corresponding computational cost.

\subsection{The training of the networks \texorpdfstring{$\Phi_{\theta_\ell}$}{Phi Theta}}

In order to train $\Phi_{\theta_\ell}$, we minimize the loss function given in \eqref{eq:lossthetal}, namely
\begin{equation*}
\mathcal L_\ell (\theta_\ell) = \frac{1}{\# \mathcal T_\ell} \sum_{B\in\mathcal T_\ell} d(\indic (B), \Phi_{\theta_\ell}(y^\delta(B))),
\end{equation*}
where
\begin{equation*}
  d(\indic (B), \Phi_{\theta_\ell}(y^\delta(B))) = \frac1N\Vert \indic(B)- \Phi_{\theta_\ell}(y^\delta(B))\Vert_2^2.
\end{equation*}
Our goal is to find the optimal $\theta_\ell^*$ such that 
\begin{equation*}
    \theta_\ell^* = \argmin_{\theta_\ell} \mathcal L_{\ell}(\theta_\ell),
\end{equation*}
which can be solved using the gradient descent algorithm. To improve the computational efficiency when $\# \mathcal T_\ell$ is large, we split $\mathcal T_\ell$ into $Q$ mini-batches, denoted as $\mathscr Q_1,\,\mathscr Q_2,\,\cdots, \mathscr Q_Q$. Then, we can iteratively update $\theta_\ell$ using these mini-batches. More precisely, one iteration is given by
\begin{equation}\label{eq:iterationmini}
    \theta_\ell \leftarrow \theta_\ell - r_q\frac{1}{\# \mathscr Q_q}\sum_{B\in \mathscr Q_q}\nabla_{\theta} d(\indic (B), \Phi_{\theta}(y^\delta(B))),
\end{equation}
where  $r_q>0$ is the learning rate. 

 We describe the detailed training process in Algorithm~\ref{Al-RTLNN}. For the simulations of this paper, we use $\mathcal E=80$ epochs and we choose the batch size equal to $256$. The learning rate begins with $r_1 = 0.5$ and decreases by a factor $2$ every $3$ epochs.  
\begin{algorithm}[ht]
 \caption{The training process for the networks $\Phi_{\theta_\ell}$}
 \label{Al-RTLNN}
    \begin{algorithmic}[1]
    \Input Training set: $\mathcal T_\ell$; The initial parameter set: $\theta_\ell^0$; The number of epochs: $\mathcal E$; The number of mini-batches: $Q$; The learning rate: $\{r_q\}^Q_{q=1}$.
    \Output The optimal parameter set $\theta_\ell^*$.

    \State Initialize  $\Phi_{\theta_\ell}$ by letting $\theta_\ell:=\theta_\ell^0$;
    \For{$\tau = 1,2,\cdots,\mathcal E$}
        \For{$q = 1,2,\cdots,Q$}    
        \State Update $\theta_l$ as in \eqref{eq:iterationmini};
        \EndFor
    \EndFor
    \State Set $\theta_\ell^*=\theta_\ell$.
    \end{algorithmic}
\end{algorithm}

\subsection{The training process of the classifier}
In this subsection, we devote ourselves to the description of the classifier's training process. In this case the training set is
$\mathcal{T}=\mathcal{T}_1\cup\dots\cup\mathcal{T}_L$, and we minimize the  cross entropy loss given by \eqref{eq:lossclass}, namely 
\begin{equation*}
\widetilde{\mathcal{L}}(\tilde\theta) = -\frac{1}{\#\mathcal{T}}\sum_{B\in\mathcal T}\tilde d(\eta(B),\Lambda_{\tilde\theta}(y^\delta(B)),
\end{equation*}
where
\[
\tilde d\bigl(\eta(B),\Lambda_{\tilde\theta}(y^\delta(B))\bigr) = \sum_{\ell=1}^3 \eta_\ell(B) \log\bigl(\Lambda_{\tilde\theta}(y^\delta(B))_\ell\bigr)
\]
and $\eta_\ell(B)=1$ if $B\in\mathcal{T}_\ell$ and $\eta_\ell(B)=0$ otherwise. 
We aim to find an optimal $\tilde\theta^*$ such that 
 \begin{equation*}
    \tilde\theta^* = \argmin_{\tilde\theta}\widetilde{\mathcal{L}}(\tilde\theta).
\end{equation*}
We adopt a training method that is similar to Algorithm~\ref{Al-RTLNN}, and we illustrate the details in Algorithm~\ref{Al-Class}  without further discussion.
For the simulations of this paper, the number of epochs $\tilde {\mathcal E}$ is set to  $80$ and the batch size is fixed as $128$. Moreover, the initial learning rate $\gamma_1$ is $0.1$ and decreases by a factor $2$ every $3$ epochs. 
\begin{algorithm}[ht]
 \caption{The training process of the classifier}
 \label{Al-Class}
    \begin{algorithmic}[1]
    \Input Training set: $\mathcal T$;  The initial (random) parameter set: $\tilde \theta$; The number of epochs: $\tilde {\mathcal E}$; The number of mini-batches: $\widetilde Q$; The learning rate: $\{\gamma_q\}^{\widetilde Q}_{q=1}$. 
    \Output The optimal parameter set $\tilde\theta^*$;

    \For{$\tau = 1,2,\cdots,\tilde {\mathcal E}$}
        \For{$q = 1,2,\cdots,\widetilde Q$}
        
        \State Update $\tilde\theta \leftarrow \tilde\theta - \gamma_q\frac{1}{\#\widetilde{\mathscr Q_q}}\sum_{B\in \widetilde{\mathscr Q_q}}\nabla_{\tilde\theta}\tilde d\bigl(\eta(B),\Lambda_{\tilde\theta}(y^\delta(B))\bigr)$, where $\widetilde{\mathscr Q_q}$ is a mini-batch of $\mathcal T$; 
        \EndFor
    \EndFor
    \State Set $\tilde \theta^*=\tilde\theta$.
    \end{algorithmic}
\end{algorithm}

\subsection{Computational cost}

During the training processes, we use the ADAM optimizer to update the parameters. All the computations are completed using Pytorch in a PC with 64GB of RAM and $11$th Gen Intel(R) Core(TM) i7-11700F CPU. Without using any parallelization technique, the training process of the classifier can be finished within $20$s, and every epoch in Algorithm~\ref{Al-RTLNN} takes about $1$min. As we have anticipated, after the training processes of the classifier and of the $\Phi_{\theta_\ell}$,  the reconstruction of one inclusion using the LRT takes less than $1$s.

\section{Histograms of SSIM distribution} \label{app:SSIM}

\begin{figure}[htbp]
\centering
 \subfigure[LRT, $\delta=0$.]
  {
	\begin{minipage}{.47\textwidth}
 	\centering
	\includegraphics[width=\columnwidth]{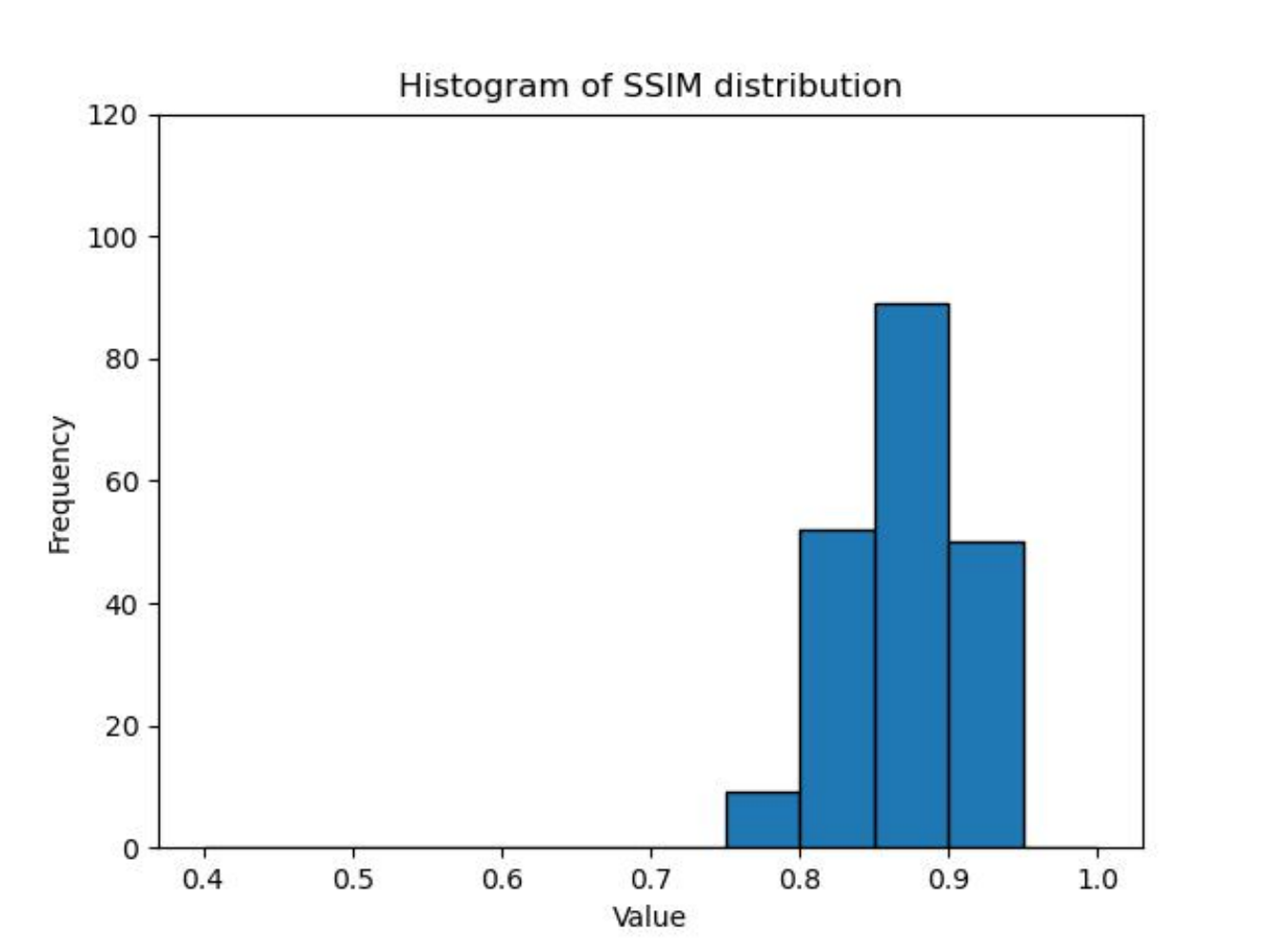}
	\end{minipage}
  }
  \subfigure[LRT, $\delta=3\%$.]
  {
	\begin{minipage}{.47\textwidth}
 	\centering
	\includegraphics[width=\columnwidth]{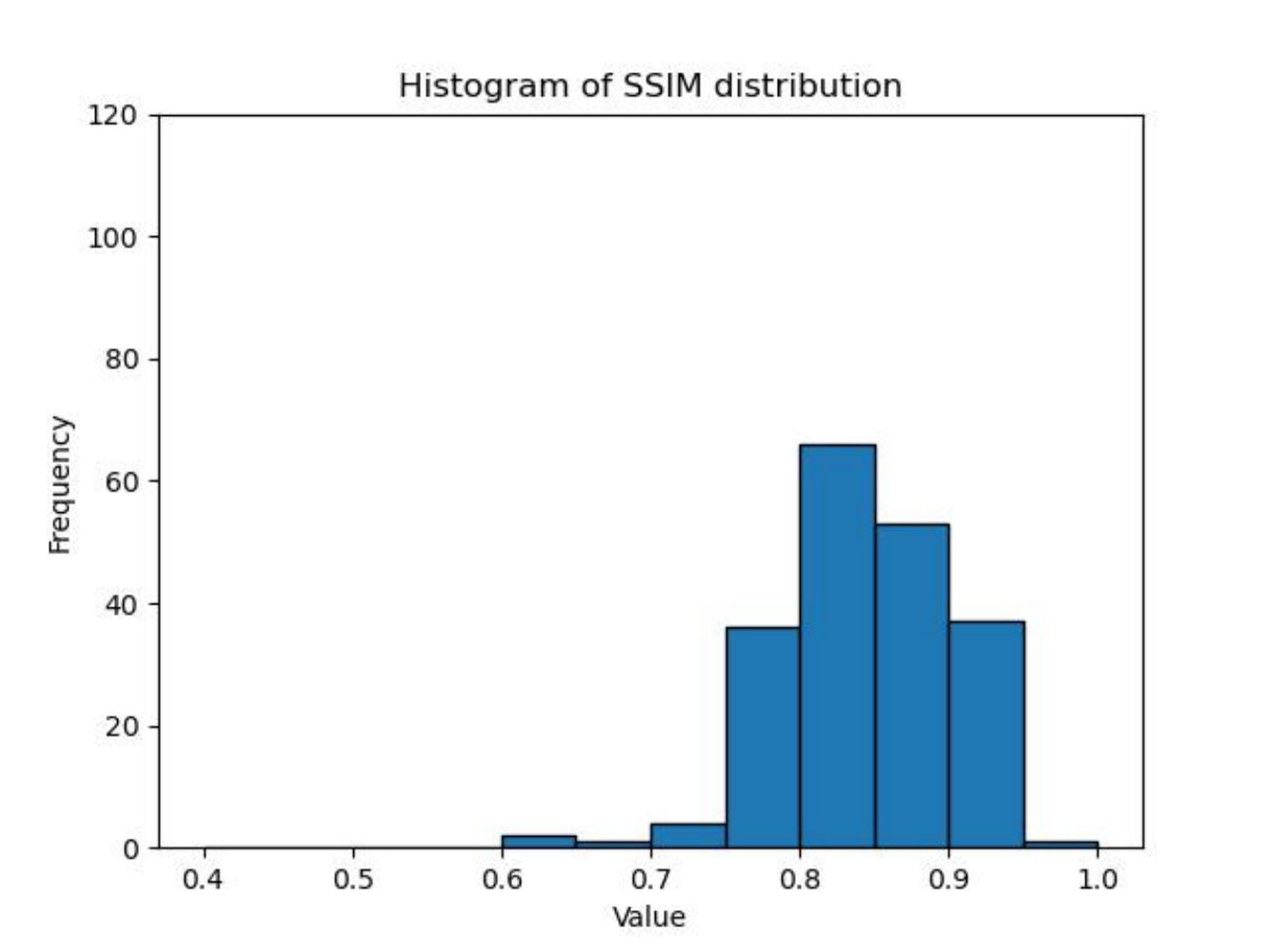}
	\end{minipage}
  }
  
   \subfigure[LRT (Random initialization), $\delta=0$.]
  {
	\begin{minipage}{.47\textwidth}
 	\centering
	\includegraphics[width=\columnwidth]{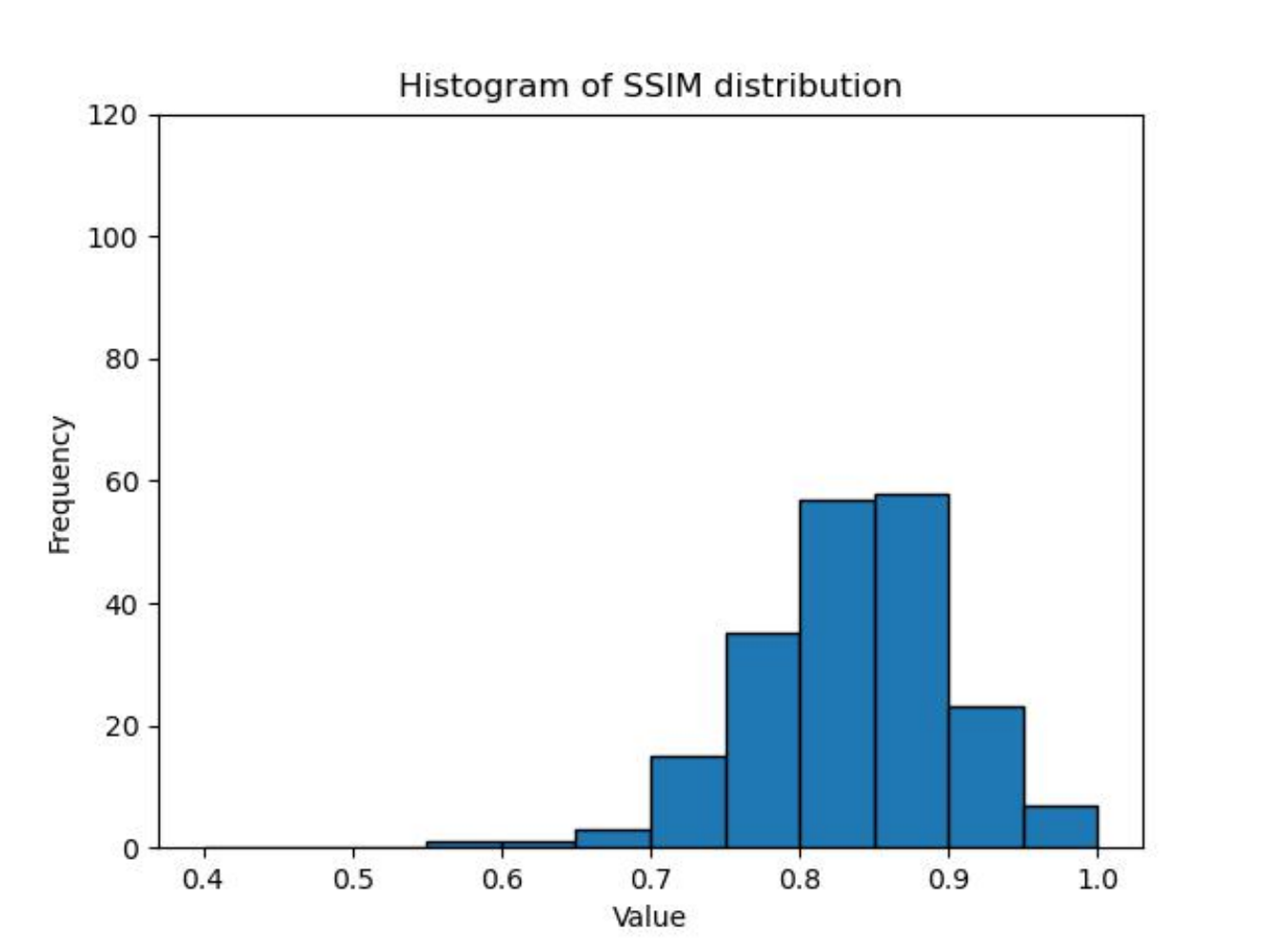}
	\end{minipage}
  }
   \subfigure[LRT (Random initialization),  $\delta=3\%$.]
  {
	\begin{minipage}{.47\textwidth}
 	\centering
	\includegraphics[width=\columnwidth]{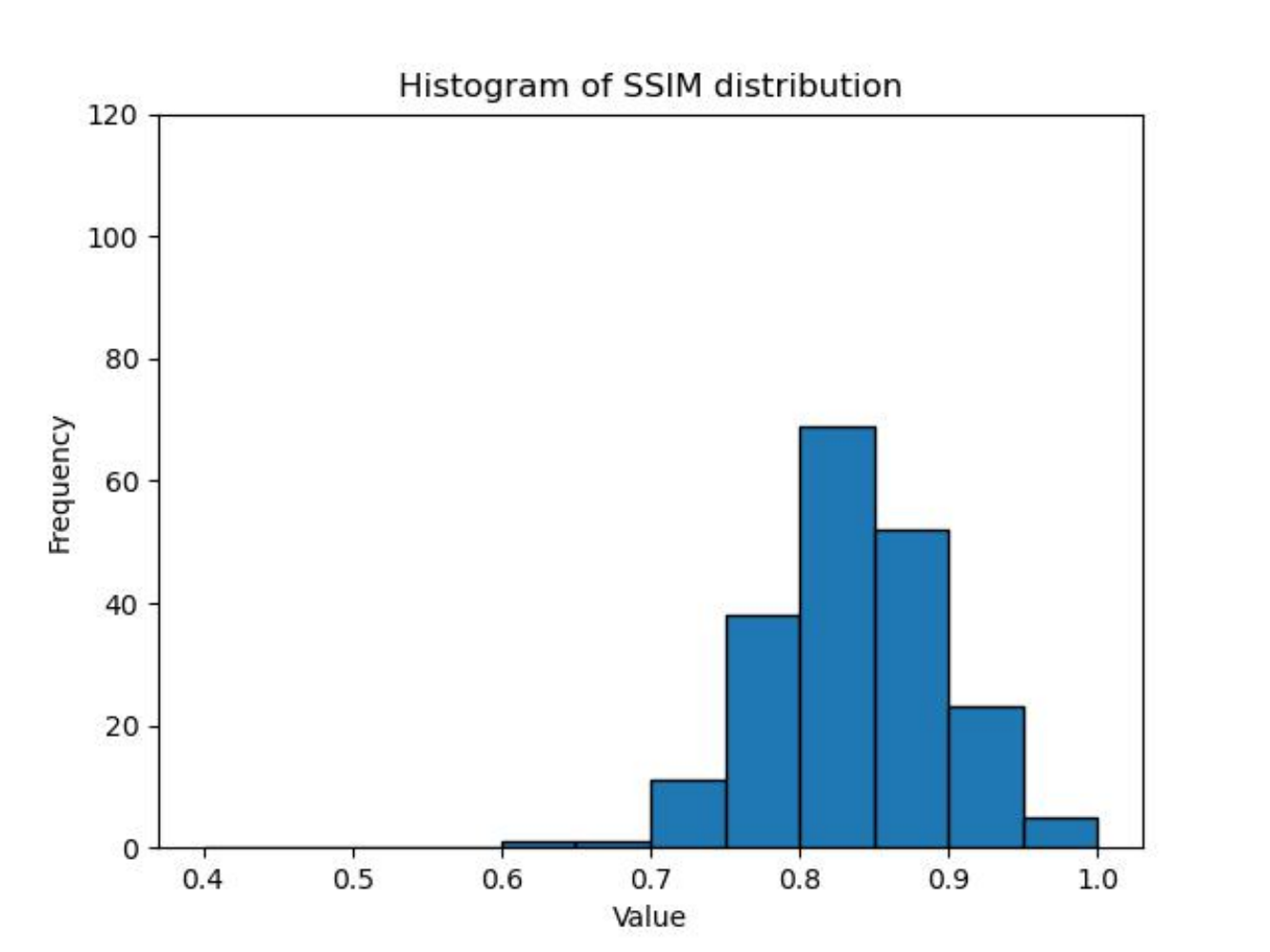}
	\end{minipage}
  }

  \subfigure[Deterministic RT, $\delta = 0$.]
  {
	\begin{minipage}{.47\textwidth}
 	\centering
	\includegraphics[width=\columnwidth]{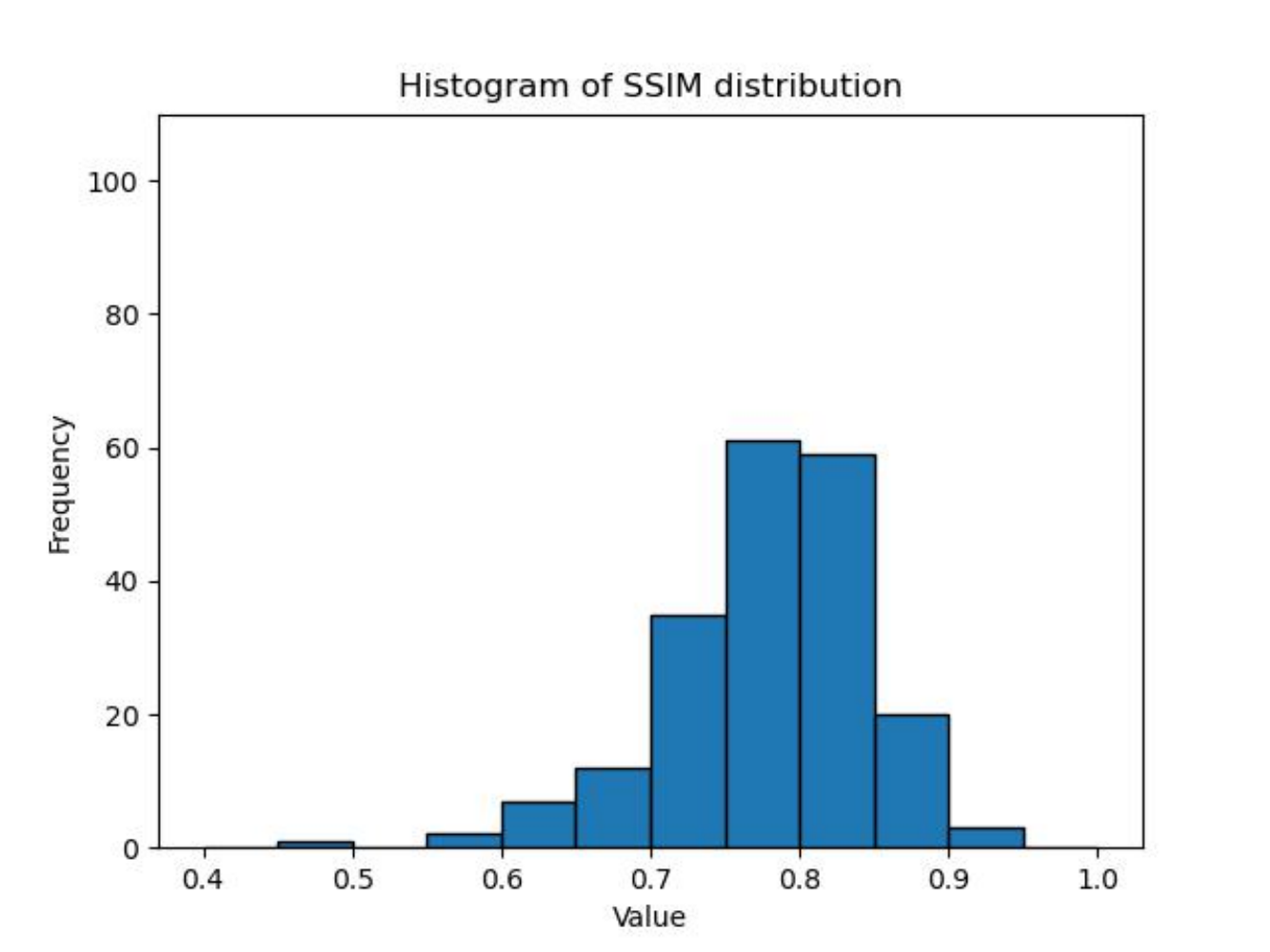}
	\end{minipage}
  }
   \subfigure[Deterministic RT, $\delta=3\%$.]
  {
	\begin{minipage}{.47\textwidth}
 	\centering
	\includegraphics[width=\columnwidth]{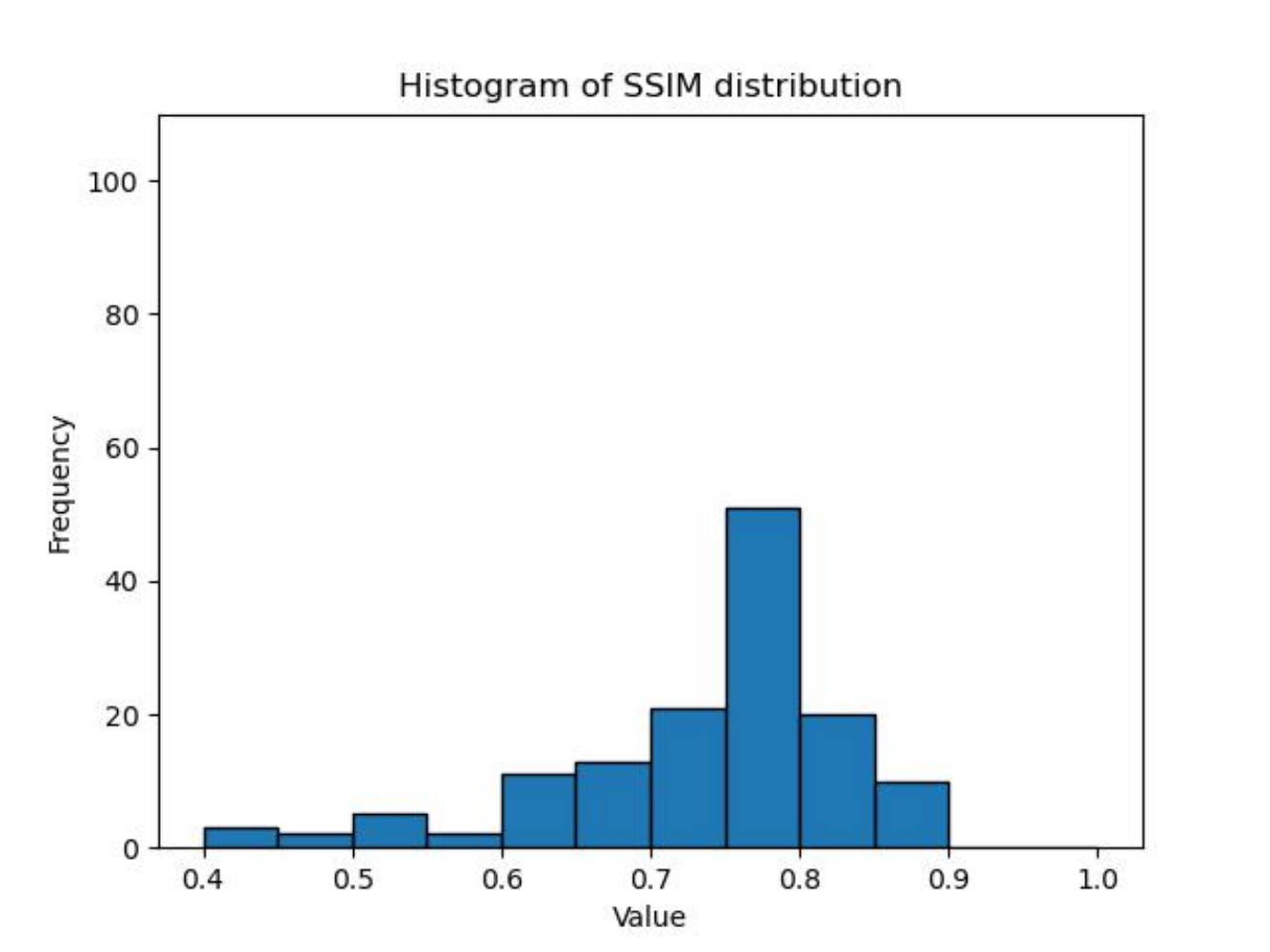}
	\end{minipage}
  }
 \caption{The SSIM distribution for all the samples in the test set with the reconstruction obtained with different methods. First row: the LRT. Second row: the LRT with random initialization during training. Third row: the RT (deterministic, no learning).}
 \label{Histogram-SSIM}
\end{figure}

In Figure~\ref{Histogram-SSIM}, we show the histograms of the SSIM distribution for the reconstructions obtained by the LRT, the LRT with random initialization and the RT with without learning. The results obtained are similar to those shown in Figures~\ref{Histogram} and \ref{Histogram-LRT-RTLNI} for the MSE.

\end{document}